\input amstex

\loadeufm
\loadmsbm
\loadeufm

\documentstyle{amsppt}
\input amstex
\catcode `\@=11
\def\logo@{}
\catcode `\@=12
\magnification \magstep1
\NoRunningHeads
\NoBlackBoxes
\TagsOnLeft

\pagewidth{5in}
\pageheight{8in}

\def \={\ = \ }
\def \+{\ +\ }
\def \-{\ - \ }

\def \b|{\big |}

\def \g1{\Gamma_1}

\def \nfp{\demo\nofrills{Proof:\usualspace\usualspace }}

\def\rarr#1#2{\smash{\mathop{\hbox to .5in{\rightarrowfill}}
 	 \limits^{\scriptstyle#1}_{\scriptstyle#2}}}

\def\larr#1#2{\smash{\mathop{\hbox to .5in{\leftarrowfill}}
	  \limits^{\scriptstyle#1}_{\scriptstyle#2}}}

\def\swarr#1#2 {\llap{$\scriptstyle #1$}  \swarrow
  	\vcenter to .5in{}\rlap{$\scriptstyle #2$}}

\topmatter
\title Diophantine Geometry over Groups IX: 
\centerline{Envelopes and Imaginaries}
\endtitle
\author
\centerline{ 
Z. Sela${}^{1,2}$}
\endauthor
\footnote""{${}^1$Hebrew University, Jerusalem 91904, Israel.}
\footnote""{${}^2$Partially supported by an Israel academy of sciences fellowship.}
\abstract\nofrills{}
This paper is the ninth in a sequence 
on the structure of
sets of solutions to systems of equations in  free and hyperbolic groups, projections
of such sets (Diophantine sets), and the structure of definable sets over  free and hyperbolic groups.
In the ninth paper we associate a Diophantine set with a definable set, and view it as the Diophantine envelope of
the definable set. We use the envelope and  duo limit groups that were used in proving  stability of the theory
of free and torsion-free hyperbolic groups [Se9], to study definable equivalence relations (imaginaries)
in these groups.
\endabstract
\endtopmatter

\document

\baselineskip 12pt

In the first 6 papers in the sequence on Diophantine geometry over groups we
studied sets of solutions to systems of equations in a free group, and
developed basic techniques and objects that are required for the analysis of sentences and
definable sets in a free group.  The techniques that we developed 
enabled us to present an iterative procedure that analyzes $EAE$ sets that are defined
in a free group (i.e., sets that are defined using 3 quantifiers), and show that
every such set is in the Boolean algebra that is generated by $AE$ sets 
([Se6],41). Hence, we obtained a quantifier elimination for definable sets in a free group.

In the 7th paper in the sequence we generalized the techniques and the results from free
groups to torsion-free hyperbolic groups, and in the 8th paper we used the techniques, that were
developed for quantifier elimination, to prove that the elementary theories of free and torsion-free
hyperbolic groups are stable.

In the 9th paper in the sequence we study definable equivalence relations (imaginaries) in free and hyperbolic groups.
Elimination of imaginaries is a fundamental notion in model theory that formalizes the idea that definable objects should admit concrete parameters. 
An $imaginary$ is, informally, an equivalence class of tuples under a definable equivalence relation. A theory $T$ has $elimination$ $of$ $imaginaries$ (EI) 
if every imaginary is interdefinable with a tuple from the home sort; equivalently, every definable quotient admits a canonical parameter 
in the original structure. Weaker forms relax this requirement. In $weak$ $elimination$ $of$ $imaginaries$ (WEI), every imaginary is only required to 
be interalgebraic with a tuple from the home sort. More generally, one may allow suitable auxiliary sorts parametrizing natural geometric objects; 
this leads to the notion of $geometric$ $elimination$ $of$ $imaginaries$.

The study of imaginaries is closely connected with the construction of $T^{eq}$, in which every definable quotient is added as a new sort. 
Thus, the question of elimination of imaginaries can be viewed as asking when these additional sorts can be coded, either in the home sort 
or in a natural collection of geometric sorts. This question has been particularly successful in algebraic and differential settings. 
The theory $ACF$ of algebraically closed fields has EI: every definable equivalence relation on a definable set admits a definable quotient 
in a finite Cartesian power of the field. Similarly, the theory $DCF_0$ of differentially closed fields of characteristic zero has EI, 
a result due to Poizat [Po], extended to several commuting derivations by McGrail [McG]. In this setting, EI is closely related to the algebraic-geometric 
description of differential-definable sets. 

Valued fields provide a particularly important example where full EI in the home sort fails, but a natural geometric form can be recovered. 
The theory $ACVF$ of algebraically closed valued fields does not eliminate imaginaries in the valued-field sort alone. However, Haskell, Hrushovski, and 
Macpherson [HHM] proved that $ACVF$ eliminates imaginaries after adding suitable geometric sorts, including sorts parametrizing lattices and their residue 
spaces. These sorts reflect the geometry naturally arising from the valuation and are, in this sense, more appropriate than attempting to code all 
imaginaries directly by field elements. Related results have been obtained for real closed valued fields and various other valued-field theories. 

Elimination of imaginaries is important because it provides a bridge between definability and geometry. It allows definable quotients to be represented 
by concrete parameters and makes their interaction with algebraic closure, types, automorphisms, and canonical bases more transparent. 
Even weaker forms can be sufficient for important applications, including the study of canonical bases and Galois groups (see [CA-FA]). In valued fields and 
other geometric theories, the use of suitable auxiliary sorts is often essential for describing interpretable objects and for transferring model-theoretic 
arguments into geometric or arithmetic settings.

Thus, the different forms of elimination of imaginaries, from full EI, through WEI, to geometric EI, measure different degrees to which definable objects 
can be represented by canonical parameters. Understanding which form holds in a given theory, and identifying the appropriate sorts in which imaginaries 
can be coded, is therefore an important structural problem with applications throughout modern model theory. 

Our goal in this paper is to prove weak elimination of imaginaries in the theories of free and torsion-free hyperbolic groups.
In an arbitrary group, there are 3 basic (not necessarily definable) families of equivalence relations: conjugation,
left (and right) cosets of subgroups, and double cosets of subgroups. 
As in general a subgroup may not be definable,
not all these equivalence relations are definable equivalence relations.

By results of M. Bestvina and M. Feighn [Be-Fe2] (an unpublished paper on negligible sets), 
of Kharlampovich-Myasnikov [KM],  and of Pillay, Perin, Sklinos and Tent [PPST], 
the only definable subgroups of a free or a torsion-free hyperbolic
group are (infinite) cyclic. Hence, the only basic equivalence relations over these groups that are definable,
are conjugation, and left, right and double cosets of cyclic 
groups${}^{3}$.

\footnote""{${}^3$ Note that we do not use the result on definable subgroups. We just indicate it for basic intuition in the introduction.} 

V. Guirardel pointed to us [Gu] that in the case of a cyclic group $C$, for every fixed non-zero pair, $\ell_1,\ell_2 \in Z$, the relation: 
$c^{\ell_1}uc^{\ell_2} \sim u$ 
where $c \in C$,
is an equivalence relation
(over a general group). Hence, in a free (or a torsion-free hyperbolic) group this is a definable equivalence relation that strictly refines
the definable equivalence relation corresponding to double cosets of two cyclic subgroups that are contained in the same maximal cyclic subgroup.

Guirardel also noticed  that in a free or a  torsion-free hyperbolic group $G$ commuting up to conjugation is a definable equivalence relation as well. i.e.,
the relation $z_1 \sim z_2$, when $z_1,z_2 \neq 1$, and there exists $g \in G$ such that $[gz_1g^{-1},z_2]=1$. 
Hence, in free (or torsion-free hyperbolic) groups there are additional basic families of definable equivalence relations, 
that generalize cosets and double cosets of cyclic groups.

Our first goal in this paper is to define the basic families of definable equivalence relations in a free group,
and show that these  basic families of definable equivalence relations
are imaginaries (i.e., not reals). In model theory, a (definable) equivalence relation is considered trivial
(and called real) if it is obtained from a definable function, i.e., if there exists some
definable function so that every equivalence class is the preimage of a point. In the second section of this paper
we define the basic families of equivalence relations in a free (or a torsion-free hyperbolic) group and prove that in free 
(and torsion-free hyperbolic) groups 
all these basic families of equivalence relations are not real.
i.e., there exist no definable functions 
so that classes in any of these equivalence relations are preimages of points (an easier proof that these basic equivalence relations are not real
appears in [Sk]).

The main goal of this paper is to classify (or represent) all the definable equivalence relations in free
and torsion-free hyperbolic groups. In particular, we aim at classifying all the imaginaries (non-reals) over these groups.
Our concluding theorems (theorems 4.4 and 4.5) show that the basic definable equivalence relations are the only
"essential" imaginaries. In particular, we show that if sorts are added for the  basic families of imaginaries
(that are defined in section 2),
then (definable) equivalence relations can
be weakly eliminated. Weak elimination means that if $G$ is a free or a torsion-free hyperbolic
group, $p$ and $q$ are $m$-tuples, and $E(p,q)$, is a definable equivalence relation,
then there exist some integers $s$ and $t$, and a
definable multi-function:
$$f: \, G^m  \ \to \ G^s \times R_1 \times \ldots \times R_t$$
where:
\roster
\item"{(1)}"  each of the $R_i$'s is a new sort for one of the  basic families of imaginaries (that are defined in section 2).

\item"{(2)}" the image of an element is uniformly bounded, and can be assumed to be of equal size if we are allowed to use constants from the group $G$.

\item"{(3)}" the
 multi-function
 is a class function, i.e., two elements in an equivalence class of $E(p,q)$ have the same
 image under the definable multi-function $f$.

\item"{(4)}" the multi-function
$f$ separates between classes. i.e., the images of elements from distinct equivalence classes
are distinct. 

\item"{(5)}" if 
 $E(p,q)$ is coefficient-free, then we can choose the definable multi-function $f$ to be
 coefficient-free.
\endroster

Note that if the definable equivalence relation is a partition of the set of tuples into finitely many definable sets, then the existence of a definable function
that satisfies (1)-(5) is obvious, as the function simply maps each of the finitely many definable sets into a different tuple. This may be either different 
elements from the group if one uses the coefficient group, or tuples of the form $e, (e,e),(e,e,e),\ldots$ where $e$ is the identity of the group, in case
one uses a coefficient free formula. Hence, the interesting case is when the equivalence relation has infinitely many equivalence classes.

In fact, we prove more than weak elimination of imaginaries, as we do get a representation of 
(generic points in) the equivalence classes
of a given definable equivalence relation as fibers of $completions$ over some (definable) parameters, where the definable parameters separate
between classes, and for each equivalence class the definable parameters admit only boundedly many values up
to the basic (definable) equivalence relations.

The main tool that we use for analyzing definable equivalence relations is the $Diophantine$ 
$envelope$ of a definable set. Given
a definable set, $L(p)$, we associate with it a Diophantine set, $D(p)$, its Diophantine envelope. 
$L(p) \subset D(p)$, and for a certain
notion of (combinatorial) genericity (which differs from the stability theoretic one), a generic point in 
the envelope, $D(p)$, is
contained in the original definable set, $L(p)$. 

If the free variables in the definable set are divided into two tuples,
$L(p,q)$, that can be interpreted as a formula that defines a parametric family of definable sets (where the free tuple $q$ is the parameter),
then we associate with $L(p,q)$ a $Duo$ (Diophantine) $Envelope$, $Duo(p,q)$, with a dual Duo limit group
(Duo limit groups were introduced in section 3 of [Se9], and served as the main tool in proving stability of free and
torsion-free hyperbolic groups). Again, $L(p,q) \subset Duo(p,q)$, and (combinatorial) generic points in $Duo(p,q)$
are contained in $L(p,q)$. 

The Duo envelope encodes all the product sets in the definable set $L(p,q)$, and such a description is crucial
in analyzing definable equivalence relations, and was also crucial in proving stability of the theory of free and hyperbolic groups in [Se9]. 

The paper is organized as follows. The first section of the paper constructs the Diophantine and Duo envelopes of definable sets. 
The construction of the envelopes
is based on the sieve procedure [Se6], that was originally used for quantifier elimination. The sieve procedure
 is the main technical
tool in proving stability of the theory [Se9], equationality of Diophantine sets [Se9], and is also the main technical
tool in analyzing equivalence relations in this paper. 

The sieve procedure itself is technical and requires constructions, techniques and notions that were introduced in our work
on quantifier elimination in the theory of free groups. Hence, we quote the description of definable sets that follows from the termination
of the sieve procedure in [Se6] (theorem 1.5), and use some of the tools and the structures that were used in the sieve procedure in this paper,
providing the exact references to the work on quantifier elimination.

In the second section we present the basic equivalence relations, that are added as new sorts in the next sections.
We use the existence of the Diophantine envelope and its properties that were presented and proved in the first section, to prove that
the  basic families of equivalence relations over free and torsion-free hyperbolic groups
are imaginaries (not reals). We believe that
envelopes can serve as an applicable tool to prove non-definability (and at times definability) in many other cases.

In the third section we start analyzing general definable equivalence relations. The Diophantine and Duo envelopes that
are constructed in the first section of the paper, depend on the defining formula, and in particular are not 
canonical. Our first goal is to find a multi-function of the type that is described above, or alternatively definable parameters, that satisfy
properties (1)-(3) from the list of properties that the eventual definable multi-function satisfies, but not necessarily part (4).

This means that the parameters that we associate with tuples in the free group in this section are indeed class functions, uniformly 
boundedly many with each class up to the basic families of equivalence relations that were presented in the second section, but these initial 
definable parameters may not separate between classes. 

To get such parameters we start with a given  (definable) equivalence relation, $E(p,q)$, we associate with it its Duo envelope, and gradually modify it.
By construction, into the Duo envelope of $E(p,q)$ 
there exists a map from a group that its values (or specializations)  provide valid proofs that testify that the specializations of the tuples, $(p,q)$, are
indeed in the given definable equivalence relation, $E(p,q)$. We use this map, i.e., the syntax of a valid proof, which is by itself the image of 
a homomorphism of a subgroup of the Duo envelope, to associate parameters with the
equivalence classes of the given equivalence relation, $E(p,q)$.  

The image of this map in the Duo envelope inherits a graphs of groups 
decomposition from the structure of the Duo envelope. In this graph of groups decomposition $<p>$ is contained in one vertex groups, and $<q>$ is contained
in another vertex groups. The parameters that we want to associate with the different equivalence classes are the values of the edge groups that connect 
between the vertex that is stabilized by $<p>$ and the vertex that is stabilized by $<q>$.

In general, it may be that the number of values that these edge groups obtain for elements in the same equivalence class is infinite (even up to
the basic equivalence relations that were defined in the second section), and we are interested
in parameters that obtain only boundedly many values for each equivalence class. Hence, we gradually modify the Duo envelope, constructing what we call
$uniformization$ limit groups. The iterative procedure that construct the uniformization limit groups terminates after finite time. 

When it terminates
in the graphs of groups that is inherited from the graph of groups of the terminal uniformization limit groups by the
image of the group that provides a valid proof that the pairs $p,q$ are indeed in the definable equivalence relation $L(p,q)$, the edge groups
admits values that depend only on the equivalence class of $p$ and $q$, and these values satisfy properties (1)-(3) of the desired multi-function. 
In particular, there is a uniform bound on the number of values of these edge groups for each equivalence class up to the basic equivalence relations that
were defined in the second section.

Hence, it is possible to construct a  definable (class) multi-function using these parameters.
However, these parameters are not guaranteed to separate between classes. Hence, the parameters and graphs of groups that we construct
in the third section
are not sufficient for obtaining weak elimination of imaginaries.

Still, the graphs of groups that are inherited from the constructed uniformization limit groups,
and their associated parameters, enable us to separate variables. i.e., to separate
the subgroup $<p>$ from the subgroup $<q>$. These subgroups 
are contained in two distinct vertex groups in the graphs of groups,
and the edge groups in between these two vertex groups are generated by finitely many elements that admit only boundedly many values (up to the basic 
imaginaries) for each equivalence class.

The separation of variables that is obtained by the graphs of groups that are inherited from the uniformization
limit groups, is the key for obtaining weak elimination of imaginaries in the fourth section of the paper.
In this section we present another iterative procedure, that combines the sieve procedure [Se6], with the procedure
for separation of variables that is presented in section 3. The combined
procedure,   iteratively constructs smaller and smaller (Duo)
Diophantine sets, that converge after finitely many steps to a Diophantine (Duo) envelope of the equivalence
relation, $E(p,q)$. Unlike what we managed to obtain in the 3rd section from the  (Duo) envelope that is constructed in the first section,
the parameters that are associated with 
the terminal envelope that is constructed by this combined iterative procedure satisfy all the properties (1)-(5) that we want our
definable multi-function to satisfy. This means that the values of the edge groups in the graph of groups that is inherited from
the graph of groups that is associated with the terminal (Duo) envelope by the subgroup
that provides valid proofs that $p$ and $q$ are in $L(p,q)$, admit only boundedly many values for each class of $L(p,q)$ 
(up to the basic equivalence relations from section 2), and that these bounded set of values distinguishes between different classes.

Therefore, the parameters that are
associated with the envelope that is constructed by the combined procedure, can be used to define the desired
multi-function, that finally proves weak elimination of imaginaries (theorems 4.4 and 4.5).

We believe that some of the techniques, notions and constructions that appear in this paper can be used to
study other model theoretic properties of free, torsion-free hyperbolic, and other groups. The arguments
that we use also demonstrate the power and the applicability of the sieve procedure [Se6] for tackling
model theoretic problems and properties. 

Finally, we would like to thank Gregory Cherlin for suggesting the study of equivalence relations to us, 
Ehud Hrushovski and Eliyahu Rips for their constant advice, Vincent Guirardel for finding mistakes in a previous version
of this paper, and the referee for many useful comments.

\vglue 1.5pc
\centerline{\bf{\S1. Diophantine and Duo Envelopes}}
\medskip
Before we analyze some of the basic imaginaries in free and hyperbolic groups, we present  
two of the main tools that are needed in order to classify the entire collection
of definable equivalence relations in these  groups. The objects that we construct may also serve as a tool in
proving that certain sets are not definable.  First, we recall the definitions of
a Duo limit group, and its associated Duo families, that were used in proving stability of free and hyperbolic groups,
and some of their properties (section 3 in [Se9]).
Then given a definable set, we associate with it a  finite collection of graded limit groups
that together form a $Diophantine$ $envelope$ of the definable set, and with them we associate a
canonical collection of duo limit groups, that in certain cases can be viewed as a $Duo$ $envelope$. 

We note that the constructions of the Diophantine and the Duo envelopes that we present in this paper
do not guarantee that these envelopes of a definable set are canonical. However, using the techniques 
that are used in the construction of the higher rank JSJ decomposition [Se10], it is possible to modify
the constructions of the envelopes and obtain canonical ones.
 
All the homomorphisms that we construct from a group or a limit group into a free or a hyperbolic group are assumed to
be $restricted$ $homomorphisms$. This means that the domain group of the homomorphism has a distinguished subgroup that is isomorphic 
to the target free or hyperbolic group, and the homomorphism maps this distinguished subgroup of the domain isomorphically onto the target. Note that given
a non-restricted homomorphism, $h:G \to \Gamma$, it is always possible to replace it by the restricted homomorphism, $\hat h:G*\Gamma \to \Gamma$,
that maps the factor $\Gamma$ in the domain isomorphically to the target, and extends the homomorphism $h$.  

Given a map $\eta: L \to Comp$, where $L$ is a limit group and $Comp$ is a completion of some resolution (not necessarily of $L$, we say that a 
homomorphism $h:L \to F_k$ $factors$ through the completion $Comp$ or through the map $\eta$, if there exists a homomorphism $\tilde h:Comp \to F_k$,
such that: $h=\tilde h \circ \eta$.

\vglue 1pc
\proclaim{Definition 1.1 ([Se9],3.1)} Let $F_k$ be a non-abelian free group, and let
$Rgd(x,p,q,a)$ ($Sld(x,p,q,a)$) be a rigid (solid) limit group
 with respect to the parameter subgroup $<p,q,a>$ (note that the parameter subgroup can be agreed to be $<p,q>$ or $<p,q,a>$ as the distinguished subgroup $<a>$ is the group
of constants and its image is fixed). Let $s$ be a (fixed) positive integer, and let
$Conf(x_1,\ldots,x_s,p,q,a)$ be a configuration limit group that is associated with the
limit group $Rgd(x,p,q,a)$ ($Sld(x,p,q,a)$) (see definition 4.1 in [Se3] for configuration limit groups). 

Recall that a $configuration$ $limit$ $group$ is obtained as a limit
of
a convergent sequence of specializations $(x_1(n),\ldots,x_s(n),p(n),q(n),a)$, that are called $configuration$
$homomorphisms$ ([Se3],4.1), in which each of the
specializations, $(x_i(n),p(n),q(n),a)$, is rigid (strictly solid), and $x_i(n) \neq x_j(n)$ for $i \neq j$
(belong to distinct strictly solid families for $i \neq j$). See section 4 of [Se3] for a detailed discussion
of these groups.

A $duo$ $limit$ $group$, $Duo(d_1,p,d_2,q,d_0,a)$, is a limit group with the following properties:
\roster
\item"{(1)}" with $Duo$ there exists an associated map: 
$$\eta: Conf(x_1,\ldots,x_s,p,q,a) \to Duo.$$
For brevity, we denote $\eta(p), \eta(q), \eta(a)$ by $p,q,a$ in correspondence.

\item"{(2)}" $Duo=<d_1>*_{<d_0,e_1>} \, <d_0,e_1,e_2>*_{<d_0,e_2>} \, <d_2>$, 

\noindent
$\eta(F_k)=\eta(<a>) \, < \, <d_0>$, $\eta(<p>) < <d_1>$,
and $\eta(<q>) < <d_2>$. 

\item"{(3)}" $Duo=Comp(d_1,d_0)*_{<d_0,e_1>} \, <d_0,e_1,e_2>*_{<d_0,e_2>} \, Comp(d_2,d_0)$, where 
$Comp(d_1,d_0)$ 
and $Comp(d_2,d_0)$ are 
(graded) completions 
with respect to the parameter subgroup $<d_0>$ (note that  $<a> \, < \, <d_0>$),
that terminate in the subgroup $<d_0>$.  $<e_1>$ and  $<e_2>$ are abelian subgroups with pegs in $<d_0>$ (i.e., abelian
subgroups that commute with non-trivial elements in the terminal limit group $<d_0>$). 

Note that we denote by $<d_0>$  the terminal limit group of the two completions,
$Comp(d_1,d_0)$ and $Comp(d_2,d_0)$, and not just a (parameter) subgroup that is contained in a distinguished vertex group in abelian decompositions of the base subgroups of
the two completions. 

\item"{(4)}" there exists a specialization $(x_1^0,\ldots,x_s^0,p_0,q_0,a)$ of the configuration
limit group $Conf$, for which the corresponding  elements $(x_i^0,p_0,q_0,a)$ are distinct
and rigid specializations of the rigid limit group, $Rgd(x,p,q,a)$ (strictly solid and belong
to distinct strictly solid families), that can be extended
to a homomorphism that factors through the duo limit group $Duo$ (i.e., there exists a configuration 
homomorphism that can be extended to a specialization of $Duo$). 
\endroster
\endproclaim

With Duo limit groups we naturally associate their collection of $duo$ $families$.

\vglue 1pc
\proclaim{Definition 1.2 } 
Given a duo limit group, $Duo(d_1,d_2,d_0)$, and a specialization of the variables
$d_0$, we call the set of specializations that
factor through $Duo$ for which the specialization of the variables $d_0$ is identical
to the given one, a $duo$-$family$. 

We say that a duo family that is associated with a duo
limit group $Duo$ is $covered$ by the  duo limit groups: $Duo_1,\ldots,Duo_t$, if there exists 
a finite collection of duo families that are associated with the duo limit groups, $Duo_1,\ldots,Duo_t$,
and a covering closure
of the duo family, so that each configuration homomorphism that can 
be extended to a specialization of a closure in the covering closure of the given duo family (of the duo limit group $Duo$), 
can also be extended to
a specialization that factors through one of the members of the  finite collection of duo families of
the duo limit groups $Duo_1,\ldots,Duo_t$ (see definition 1.16 in [Se2] for a covering closure).
\endproclaim

In [Se9] we used the sieve procedure [Se6] to prove the 
existence of a finite collection of duo limit groups, that cover
all the duo families that are associated with all the duo limit groups that are associated with a given rigid
limit group.

To state the theorem from [Se9] one needs to replace a duo limit group by a $covering$ $closure$ of itself. A (graded) $closure$ and a covering
closure are defined in section 1 in [Se2]. An  (ungraded) closure of a completion has the same structure as a completion, i.e., the same sequence of abelian
decompositions and retractions, except that abelian vertex groups along the completion are replaced by finite (free abelian) supergroups of finite 
index.

A graded closure is a closure of a graded completion. It has the same structure as the graded completion, some of the abelian vertex groups along the 
graded completion are replaced by finite index supergroups of finite index in the graded closure (like in an ungraded closure), and the bottom level 
of the graded closure is replaced by a possibly different rigid or solid limit group w.r.t\ the same subgroup of parameters. Furthermore, there is a map
from the terminal limit group of the graded completion to the terminal limit group of its closure, where this map maps rigid vertex groups and edge groups 
in the terminal limit group of the graded completion into 
rigid vertex groups and edge groups in the graded closure.

For the purpose of constructing closures, a duo limit group can be viewed as a graded completion. i.e., its two completions can be viewed as a single
completion, their union. Hence, a closure of a duo limit group has the same
structure as the duo limit group, except abelian vertex groups that are replaced by finite index supergroups, and the terminal limit group $<d_0>$ that
is replaced by a rigid or a solid limit group w.r.t.\ $<d_0>$.  

\vglue 1pc
\proclaim{Theorem 1.3 ([Se9],3.2)} Let $F_k$ be a non-abelian free group, let $s$ be a positive
integer, and let
$Rgd(x,p,q,a)$  be a rigid  limit group (with respect to the parameter subgroup $<p,q,a>$)
over $F_k$.  

There exists a finite collection of duo limit groups that are associated
with configuration homomorphisms of $s$ distinct rigid homomorphisms of $Rgd$, $Duo_1,\ldots,Duo_t$, 
and some global bound $b$, 
so that every duo family that is associated with a duo limit group $Duo$, that is associated
with configuration homomorphisms of $s$ distinct rigid homomorphisms of $Rgd$, 
is covered by the given finite collection $Duo_1,\ldots,Duo_t$. 

Furthermore, every duo family that is associated with an arbitrary duo limit group $Duo$, that is associated with $Rgd$, is covered by
at most $b$ duo families that are associated with the given finite collection,
$Duo_1,\ldots,Duo_t$.
\endproclaim

In this section we look for a partial generalization of theorem 1.3 to a general definable set. Given a definable
set $L(p,q)$ we associate with it  a finite collection of graded completions 
(with respect to the parameter subgroup
$<q,a>$). A "generic" point in each of these graded completions is contained in the definable set $L(p,q)$,
and given a value $q_0$ of the variables $q$ the (boundedly many) 
fibers that are associated with $q_0$ and the finite collection of graded completions, contain the projection 
$L(p,q_0)$ of the definable set $L(p,q)$. 

We call a collection of graded completions that satisfy these properties w.r.t.\ a given definable set $L(p,q)$, a  graded $envelope$ of the definable set.
Later on we associate with this  collection of graded
completions, a  collection of duo limit groups, that we call $duo$ $envelope$, and the obtained duo limit groups is the key tool
in our classification of imaginaries. 

\vglue 1pc
\proclaim{Definition 1.4}   
A completion is defined in section 1.12 in [Se2], and a graded completion in section 3 in [Se2]. A graded completion $GComp(z,q,a)$ 
 terminates with a rigid or a solid limit groups w.r.t.\ the parameter subgroup $<q,a>$. A fiber of values
of the graded completion $GComp$ is associated with an assignment of values $q_0$ to the defining parameter $q$, and with a rigid or a strictly solid
family of values of the (terminal) base  rigid or solid limit group of $GComp$ (note that there are boundedly many rigid or strictly solid families 
that can be associated
with $q_0$ by [Se3]). The fiber that is associated with such a rigid or a strictly solid family
is the set of values, $(z,q_0,a)$,  that factor through $GComp$ and restrict to the values of the (fixed) rigid or strictly solid family of values of 
the base (terminal) limit group of
$GComp$.
\endproclaim

\vglue 1pc
\proclaim{Theorem 1.5} Let $F_k$ be a non-abelian free group, and let $L(p,q)$ be a definable set  
in $F_k$.
There exists a finite collection of graded completions, $G_1(z,p,q,a),\ldots,G_t(z,p,q,a)$, 
that are associated with
$L(p,q)$ and graded w.r.t.\ $<q,a>$, that we call a $Diophantine$ $Envelope$ of $L(p,q)$, for which:
\roster
\item"{(1)}" For each $j$, $1 \leq j \leq t$, there exists a  test sequence $\{(z_n,p_n,q_0,a)\}$ 
of the completion $G_j(z,p,q,a)$, for which all the specializations $(p_n,q_0) \in L(p,q)$.

Furthermore, for each index $j$, $1 \leq j \leq t$, there exists a test sequence 
$\{(z_n,p_n,q_n,a)\}$ 
of the completion $G_j(z,p,q,a)$, for which all the specializations, $(p_n,q_n) \in L(p,q)$,
and the Zariski closure of the test sequence is the graded completion, $G_j(z,p,q,a)$. 

\item"{(2)}" Given a specialization  $(p_0,q_0) \in L(p,q)$,   
there exists an index $j$, $1 \leq j \leq t$,  and a test sequence 
$\{(z_n,p_n,q_0,a)\}$ 
of the completion $G_j(z,p,q,a)$ that are all in the same fiber of the graded completion, $G_j(z,p,q,a)$, and
for which all the specializations $(p_n,q_0) \in L(p,q)$.

Furthermore,  $(p_0,q_0)$
can be extended to a
specialization that factors through the same fiber of the completion $G_j(z,p,q,a)$ that contains
the test sequence, 
$\{(z_n,p_n,q_0,a)\}$.
\endroster
\endproclaim

\nfp Before we construct the $envelope$ of a definable set, i.e., a collection of graded closures that satisfy the properties in the theorem,
we recall some of the properties of the output of the sieve procedure for quantifier
elimination that is presented in [Se6].

Let $L(p)$ be a definable set. 
With a definable set $L(p)$ the sieve procedure in [Se6] associates a finite collection of graded (PS) resolutions
(graded w.r.t.\  the subgroup $<p,a>$), $PSRes_1,\ldots,PSRes_s$.
Each of these graded resolutions terminates in a rigid
or a solid limit group (with respect to the parameter subgroup $<p,a>$).
With each such graded resolution we associate its (graded) completion, and  the sieve procedure further associates with the graded
completion a finite collection of graded
closures 
(see definitions  1.25-1.30 in [Se5] for the exact definitions and construction of these graded closures). 

By the sieve procedure [Se6], that eventually leads to
quantifier elimination in the theory of a free group, the definable set $L(p)$ is equivalent to those 
rigid and strictly solid specializations
of the terminal rigid and solid limit groups of the PS resolutions, $PSRes_1,\ldots,PSRes_s$,
for which the  fibers of the graded completions of these PS resolutions, that are associated with the rigid and strictly solid specializations of
their terminal groups, are
not covered by the corresponding fibers of the collection of graded closures that are associated with the PS resolutions.

Therefore, using the output of the sieve procedure, i.e.,  the PS resolutions and the graded closures that are associated with them, it is possible to
associate definable sets with each terminal limit group of the PS resolutions, $PSRes_1,\ldots,PSRes_s$, and from these definable sets it is possible 
to redefine the original definable set  $L(p)$ using a formula that is in the Boolean algebra of AE sets, that implies quantifier elimination.
In the sequel we will need to use these definable sets and their properties.

With the terminal rigid or solid limit group of a PS resolution, $PSRes_j$, that we denote $Term_j$, one can associate  a sets of values
of the defining parameters $p$, $B(Term_j)$, 
that are all in the Boolean algebra of AE sets, and such that $B(Term_j) \subset L(p)$. Furthermore, the union of the subsets $B(Term_j)$, $1 \leq j \leq t$,
is the definable set $Lp)$.  To define and construct these sets we first need to recall the definition of strictly solid
families of a solid limit group w.r.t.\ a given set of closures.

In section 1 of [Se3] we defined rigid and strictly solid families of a graded limit group. In section 3 of [Se3] we further generalized the notion
of a strictly solid family of a solid limit group to a strictly solid family w.r.t.\ a covering closure. (Definition 2.12 in [Se3])
  
Given a graded limit group $G(x,p,a)$ we associated with it its graded Makanin-Razborov diagram (section 10 in [Se1]). The graded MR diagram
contains finitely many (graded) resolutions that terminate in either rigid or solid limit groups. Let $GRes_1,\ldots,GRes_t$ be the graded resolutions
in the graded MR diagram of the graded limit group $G(x,p,a)$. Suppose that with with each of these graded resolutions we associate a finite collection
of graded closures (graded closures are defined in section 3 of [Se2]). 

Let $Term_j$ be the terminal solid limit group of a graded resolution $GRes_j$ in the MR diagram of $G(x,p,a)$. Let $p_0$ be a value of the parameters $p$.
Suppose that for every graded resolution, $GRes_{\hat j}$, $1 \leq \hat j \leq t$, that starts with the path that is associated with the graded resolution
$GRes_j$ in the MR diagram of $G(x,p,a)$ and continues further in the MR diagram, the set of closures that is associated with $GRes_{\hat j}$ forms
a covering closure for all the fibers of $GRes_{\hat j}$ that are associated with $p_0$.

For such a $p_0$ a value $(p_0,x_0)$ belongs to a strictly solid family of $Term_j$ w.r.t.\ the given set of closures,
if there is a value in the same family of $(p_0,x_0)$ (a family w.r.t.\ the modular group of the solid limit group $Term_j$) that does 
not factor through any of the closures that are associated with the graded resolutions $GRes_{\hat j}$. In section 3 of [Se3] it is proved that 
the number of such strictly solid families, that are associated with a solid limit group $Term_j$ and a value $p_0$ of the defining parameters
$p$, is uniformly bounded for all the values of the defining parameter $p$ (Theorem 2.13 in [Se3]). 

The notion of rigid and strictly solid families w.r.t.\ a given collection of closures enables one to describe the main conclusion of the
sieve procedure in [Se6] that we use repeatedly in this paper, and in particular in the proof of theorem 1.5.

\vglue 1pc
\proclaim{Theorem 1.6 (Theorems 39 and 40  in [Se6])} Let $F_k$ be a non-abelian free group, and let $L(p)$ be a definable set. Then there exist finitely
many graded resolutions, $PSRes_1,\ldots,PSRes_s$, such that $PSRes_j$, $1 \leq j \leq s$, terminates in a rigid or solid limit group $Term_j$.
Furthermore, with each graded resolution $PSRes_j$ there are finitely many associated graded closures.

For each $j$, $1 \leq j \leq s$ there is a definable set $B(Term_j)$ that is associated with the terminal rigid or solid limit group $Term_j$.
of $PSRes_j$. The set $B(Term_j)$ have the following properties:

\roster
\item"{(1)}" $B(Term_j) \subset L(p)$. 

\item"{(2)}" $\Cup_{j=1}^t \, B(Term_j)=L(p)$.   

\item"{(3)}" The sets $B(Term_j)$ are in the Boolean algebra of AE sets.

\item"{(4)}" The set $B(Term_j)$ is the collection of values $p_0$  of the defining parameters $p$, for which $p_0$ extends to a rigid or a strictly solid
value of $Term_j$ w.r.t.\ the given set of graded closures of the graded resolutions, $PSRes_1,\ldots,PSRes_s$,  and such that the fiber that is 
associated with this rigid or strictly solid value is not covered by the collection of graded closures that are associated with $Term_j$.
\endroster
\endproclaim

\medskip
Let $L(p,q)$ be a definable set.
We start the proof of theorem 1.5, i.e., the construction of the finite collection of graded limit groups that are associated with the definable
set $L(p,q)$ and satisfy the conclusions of theorem 1.5,
by associating a finite collection of graded limit groups (that are graded with respect to $<q,a>$) with the sets $B(Term_j)$, where $B(Term_j)$ are the
definable sets that are associated with the definable set $L(p,q)$ according to theorem 1.6, i.e., according to the output of the
sieve procedure in [Se6].

\noindent
With the definable set $L(p,q)$ the sieve procedure associates finitely many PS resolutions, $PSRes_1,\ldots,PSRes_s$, that terminate in the
rigid or solid limit groups, $Term_1,\ldots,Term_t$. Note that these terminal limit groups are rigid or solid w.r.t.\ the parameter subgroup $<p,q,a>$.
Let $Term_j$ be one of these rigid or solid limit groups. With $Term_j$ and its associated graded resolution, $PSRes_j$, the sieve procedure associates 
finitely many graded closures, that are all graded closures of the (graded) completion of $PSRes_j$. 

Recall that by theorems 2.5, 2.9
and 2.13 in [Se3], there exists a global bound on the number of rigid specializations of a rigid limit group,
and a global bound on the number of strictly solid families of specializations of a solid limit group and of strictly
solid families of specializations w.r.t.\ a given set of closures, for
all the possible values of the parameter subgroup. Hence, with each value of the parameter
subgroup $<p,q,a>$, there are boundedly many rigid (strictly solid families of) specializations (w.r.t.\ a given collection of closures)  of the
(rigid or solid) terminal limit group $Term_j$, and of the the terminal (rigid or solid) limit groups of the graded closures that are associated with
$PSRes_j$. (note that $Term_j$ is mapped to each of the (rigid or solid) terminal limit groups of the graded closures that are associated with $PSRes_j$).

$Term_j$ is rigid or solid w.r.t.\ the parameter subgroup $<p,q,a>$. To associate the graded envelope of the definable set $L(p,q)$, i.e., to associate
with $L(p,q)$ a finite collection of graded completions and closures, we need to change the grading look at the terminal limit groups $Term_j$ as graded
limit groups w.r.t.\ the parameter subgroups $<q,a>$. Once we change the grading, we will proceed with a procedure that is similar to the sieve procedure
for quantifier elimination in [Se6] to obtain the graded completions and closures that are listed in the theorem 1.5. 

Given $Term_j$, we look at the collection of specializations of the form:
$$(x,y_1,\ldots,y_w,u_1,\ldots,u_m,v_1,\ldots,v_n,r,r_1,\ldots,r_m,p,q,a)$$ where:
\roster
\item"{(1)}" the integers $t$,$m$,$n$ are non-negative integers that are bounded by the sum of the global bounds on the number of rigid and 
strictly solid 
families w.r.t.\ the given set of closures of  the terminal rigid and solid limit groups of the graded resolutions $PSRes_1,\ldots,PSres_s$,
and the terminal limit groups of the graded closures that are associated with them.

Note that the triples: $w$, $m$, $n$ obtain all the possible values for the number of these rigid and strictly solid families of values of the corresponding
rigid and solid terminal groups.

\item"{(2)}" the specialization $(x,p,q,a)$ is a rigid or a strictly solid specialization of the terminal
rigid or solid limit group $Term_j$ with respect to the given collection of closures. The specializations $(y_i,p,q,a)$, $i=1,\ldots,w$, are rigid and 
strictly solid specializations of the terminal rigid and solid limit groups of the closures that are associated
with the graded resolution $PSRes_j$ and its terminal limit group $Term_j$, that extend the rigid or strictly solid family that is associated with
$(x,p,q,a)$.

The rigid
specializations are distinct and the strictly solid specializations belong to distinct strictly solid families,
and the finite collection of specializations $(y_i,p,q,a)$, $i=1,\ldots,w$, represent all the rigid and
strictly solid families of specializations that are associated with (i.e., that extend)
the rigid or strictly solid specialization
$(x,p,q,a)$.

\item"{(3)}" the specializations $(u_i,p,q,a)$, $i=1,\ldots,m$,  are all the distinct rigid and strictly solid 
specializations of the terminal limit groups $Term_{\hat j}$ that are terminal groups of graded resolutions $PSRes_{\hat j}$, that extend the path
that is associated with $PSRes_j$ in the  graded Makanin-Razborov diagram from which both are taken. If $Term_{\hat j}$ is solid then 
$(u_i,p,q,a)$ is from a strictly solid family w.r.t.\ the given set of closures that are associated with $PSRes_1,\ldots,PSRes_s$.

\item"{(4)}" the specializations $(v_i,p,q,a)$, $i=1,\ldots,n$,  are all the distinct rigid and strictly solid 
specializations
of the terminal (rigid and solid) of the graded closures that are associated with the various terminal limit groups $Term_{\hat j}$,
that extend the rigid or strictly solid values $u_1,\ldots,u_m$ of the rigid or solid terminal limit groups $Term_{hat j}$.

\item"{(5)}" the specializations $r$,  include primitive roots of the edge groups in the
graded abelian decomposition  
of the rigid or solid limit group $Term_j$ that are associated with the specialization $(x,p,q,a)$, 
and they indicate what powers of the primitive roots are
covered by the associated fibers of the graded closures that are associated with $Term_j$ and with the rigid and strictly solid families:
$y_1,\ldots,y_w$.

\item"{(6)}" the specializations $r_i$, $i=1,\ldots,m$,  include primitive roots of the values of the edge groups in the
graded abelian decompositions of the solid limit groups $Term_{\hat j}$ that are associated with the strictly solid families that are
associated with the values: $u_1,\ldots,u_m$. 
\endroster

The specializations $(x,y,u,v,r,p,q,a)$ that satisfy properties (1)-(6) form "proof statements" for validation
that $(p,q) \in L(p,q)$. i.e., a value $(p_0,q_0) \in L(p,q)$ if and only if there is an extended value $(x_0,y_0,u_0,v_0,r_0,p_0,q_0)$
that testifies that $(x_0,p_0,q_0)$ is a rigid or a strictly solid value w.r.t.\ the given closures of the graded resolutions,
$PSREs_1,\ldots,PSRes_s$ of one of their terminal limit groups $Term_j$, and the graded closures that are associated with $PSRes_j$ do not
cover the fiber that is associated with $(x_0,p_0,q_0,a)$.

By our standard methods (presented in section 5 of [Se1]), with this collection of specializations, for all the terminal limit groups
$Term_j$ of the graded resolutions, $PSRes_1,\ldots,PSRes_s$,  we can
canonically associate a finite collection of graded limit groups (which is the Zariski closure of
the collection). We view each of these (finitely many) limit groups, as graded with respect to the
parameter subgroup $<q,a>$. We associate with each such limit group its graded taut Makanin-Razborov diagram
 (see proposition 2.5 in [Se4] for
the construction of the taut Makanin-Razborov diagram), 
 and
with each (taut) graded resolution in the diagram we associate its (graded) completion. 

The graded resolutions that appear in the graded MR diagram contains all the "proof statements" $(x_0,y_0,u_0,v_0,r_0,p_0,q_0,a)$ 
for pairs $(p_0,q_0) \in L(p,q)$. However, it is not true that generic points (i.e., test sequences) in these resolutions, or generic points in any
fiber in them (a fiber that is associated with a value $q_0$ of the variables $q$), restrict to pairs $(p_0,q_0) \in L(p,q)$.
In order to gradually transfer these initial graded resolutions into graded resolutions or their completions, that satisfy the conclusions of 
theorem 1.5, we need to apply (a version of) the sieve procedure that was introduced in [Se6].

The sieve procedure is rather involved technically, and to learn its full details the interested reader is referred to [Se6]. The rest of the proof
of theorem 1.5 adapts it to our context.
We continue with each of the obtained graded completions in parallel. With each graded completion we associate canonically a finite set of closures
precisely as was done definitions 1.25-1.30 in [Se5]. 

We look at all the test sequences of the graded completions for which the "proof statement" fails. That is  rigid or  strictly solid values
from the proof statements along a test sequence, $\{(x_{\ell},y_{\ell},u_{\ell},v_{\ell},r_{\ell}p_{\ell},q_{\ell}.a)\}$, are in fact not rigid or not
strictly solid (w.r.t.\ given closures of the initial resolutions $PSRes_1,\ldots,PSRes_s$), 
or there exist extra rigid or strictly solid values that do not appear in the "proof statement",
or the fibers that are associated with $\{(x_{\ell},p_{\ell},q_{\ell},a)\}$ are covered by the graded closures that are associated with the graded resolution
$PSRes_j$ that is associated with these rigid or strictly solid values, or the elements $\{r_{\ell}\}$ that are supposed to be primitive roots are actually not 
primitive roots.

This collection of test sequences that (potentially) fail to be valid proof statements enable us to associate with each of the graded completions that we
started with a finite collection of graded closures. We construct these graded closures using the constructions of graded formal limit groups that appear 
in sections 2-3 in [Se2], precisely as is done in the sieve procedure (definitions 1.25-1.30 in [Se5]). 

With those graded closures that demonstrate that there 
are extra rigid or strictly solid families that do not appear in the proof statement, that we (still) denote $ExtraPS$ (as their notation in the sieve 
procedure in [Se6]), we further associate $generic$ $collapse$ (graded) resolutions. We construct these generic collapse resolutions, that are graded closures 
of the $ExtraPS$ resolutions, by 
looking at all the test sequences of the ExtraPS resolutions for which the (potential) extra rigid or strictly solid family that s
supposed not to  appear in proof statements, 
is either not rigid or not strictly solid or it belongs to a rigid or a strictly solid family that actually appears in the proof statement.

These last test sequences enable us to construct a finite collection of graded closures of the original completion, for which potential 
obstruction to the validity of "generic" proof statements actually collapse (these closures are called $generic$ $collapse$ in [Se6]).

\medskip
As in the sieve procedure in [Se6], we continue to the next step only with the graded closures of the graded resolutions that we started the first step with,
that we called
ExtraPS resolutions. i.e., those graded closures that collect test sequences of the original completions for which it is possible to find extra rigid
or strictly solid families of values (w.r.t.\ the original set of graded closures of $PSRes_1,\ldots,PSRes_s$) 
that do not appear in the proof statements along the test sequences. 
Note that such a closure, as well as all the graded resolutions and closures that we construct along the procedure, is graded with respect to the parameter subgroup $<q,a>$.
 
Given an extra PS resolution $EXtraPS$ we look at all the values $(x_0,y_0,u_0,v_0,r_0,,b_0,p_0,q_0,a)$ that factor through it, where $b_0$ is the value
of the (potential) extra rigid or strictly solid value that is not part of the proof statement (but $b_0$ are elements in the graded closure
$ExtraPS$), for which the potential extra value $b_0$ is either not
rigid or not strictly solid (w.r.t.\ the given original collection of closures of the original resolutions, $PSRes_1,\ldots,PSRes_s$), or it belongs
to the same family as a rigid or a strictly solid family that are part of the proof statement. i.e., $(b_0,p_0,q_0,a)$ belongs to the same rigid or
strictly solid family as one of the elements, $y_0,u_0,v_0$.

Demonstrating that an element is not rigid nor strictly solid (w.r.t.\ a given collection of closures), 
or that it is in the same rigid or strictly solid of one of finitely many given families,
is equivalent to imposing one out of finitely many  Diophantine conditions on the elements $(x_0,y_0,u_0,v_0,r_0,b_0,p_0,q_0,a)$. Hence, we impose
each of these (finitely many $collapse$) Diophantine conditions on the elements that factor through the graded closures, $ExtraPS$, and look at the collection
of elements that factor through these extra PS resolutions (or closures) $ExtraPS$ and satisfy one of the finitely many collapse Diophantine conditions.  

With the collection of these specializations, that are composed from values that factor through the graded closures ExtraPS and values that satisfy one 
of the finitely many collapse Diophantine conditions, we naturally 
associate (canonically) a collection of (maximal) graded (with respect to $<q,a>$)
limit groups (which is the Zariski closure of the collection). With these limit groups we associate a collection 
of graded taut resolutions (with respect to the subgroup $<q,a>$), 
that are not taken from the graded Makanin-Razborov diagrams that are associated with the graded limit groups, but rather 
constructed according to the first step
of the sieve procedure (that is presented in [Se6]). Some of these resolutions are graded closures of the
graded completion that we have started with, (i.e., of the graded closures $ExtraPS$) and others have strictly smaller complexity (in the sense of the complexities
that are used along the sieve procedure in [Se6]). We continue to the next step 
only with the completions of each of these graded resolutions (with respect to the parameter subgroup $<q,a>$) that are
not graded closures of the graded completion that we have started the first step with. i.e., only with those graded resolutions that have strictly 
smaller complexities than the complexity of the original ExtraPS resolution we started this step with. 

For a detailed description of the construction of the obtained graded resolutions see the first step of the sieve procedure in [Se6]. With the constructed
graded resolutions, we associate finitely many graded closures as we did with the original graded resolutions: $PSRes_1,\ldots,PSRes_s$.  To construct these graded 
closures we look at all the test sequences of the constructed graded resolutions for which elements that were supposed to be rigid or strictly solid are in fact
not rigid or not strictly solid,
or for which there exist additional rigid or strictly solid families of values that do not appear in the proof systems in the test sequences, or
for which elements that are supposed to be primitive roots in the proof system are not, or for which the fibers that are associated with the
elements $(x_{\ell},p_{\ell},q_{\ell},a)$ from the test sequences are covered by the graded closures that are associated with the original graded resolutions,
$PSRes_1,\ldots,PSRes_s$.  

We continue iteratively to the next steps as was done in the sieve procedure in [Se6].
At each step we start with the ExtraPS resolutions that were constructed in the previous step and impose on them one of finitely many 
Diophantine conditions that force the additional rigid or strictly solid values (w.r.t.\ the given set of graded closures of
the original graded resolutions $PSRes_1,\ldots,PSRes_s$) that do not appear in the associated proof statements to either be non-rigid or
not strictly solid or to coincide with one of the families that are part of the corresponding proof statements.

We collect the specializations that factor through the ExtraPS resolutions that were constructed in the previous step and satisfy one of
the collapse Diophantine conditions, and analyze them using the procedure that was used in the general step of the sieve procedure in [Se6].
Since the analysis and the procedure that we use is identical to the one that is used in the sieve procedure for quantifier elimination in [Se6],
the iterative procedure that we use terminates after finitely many steps, according
to the proof of the termination of the sieve procedure that is given in section 4 in [Se6] (theorem 22 in [Se6]).

Let $\hat G_1(z,p,q,a),\ldots,\hat G_{t'}(z,p,q,a)$ be the graded completions that are constructed along the various steps of the iterative 
procedure presented above and associated with $L(p,q)$.
$\hat G_i(z,p,q,a)$, $1 \leq i \leq t'$ are graded completions w.r.t.\ the parameter subgroup $<q,a>$.
We further replace each of the graded completions $\hat G_i(z,p,q,a)$, by the graded completions that are associated with the
Zariski closure of all the test sequences, $(z_{\ell},p_{\ell},q_{\ell},a)$,
that factor through $\hat G_i(z,p,q,a)$, and for which the proof statements that the elements in the test sequences restrict to are valid proof
statements. 

Note that the graded completions that are associated with the Zariski closure of these test sequences that factor through the graded completion 
$\hat G_i(z,p,q,a)$ is either empty 
(if there is no such test sequence), or it is the graded
completion $\hat G_i(z,p,q,a)$ itself, or there are finitely many graded closures that are associated with the Zariski closures, and they all have the same
structure as $\hat G_i(z,p,q,a)$ and they are proper quotients of $\hat G_i(z,p,q,a)$.
  
We denote the obtained graded completions that are associated with the Zariski closures of the above test sequences:  $G_1(z,p,q,a),\ldots,G_t(z,p,q,a)$.
With each of the obtained completions we further associate finitely many graded closures, as we did along the procedure. Recall that these graded closures are
constructed from test sequences of the graded completions, $G_i(z,p,q,a)$, $1 \leq i \leq t$, for which elements that are supposed to be rigid or strictly solid
in the corresponding proof statements are not, or for which there are extra rigid or strictly solid elements that are not in the families that
appear in the proof statements, or for which elements that are supposed to be primitive roots are not, or for which the fibers that are associated with
the elements $(x_{\ell},p{\ell},q_{\ell},a)$ in the proof statements are covered by the graded closures that are associated with the corresponding PS
resolution $PSRes_j$.  

The collection of graded completions,
$G_1(z,p,q,a),\ldots,G_{t}(z,p,q,a)$, clearly satisfies part (1)  of theorem 1.5. The iterative construction that was used to construct these 
graded completions 
guarantees that they satisfy part (2) as well.

\line{\hss$\qed$}

Given a definable set, $L(p,q)$, theorem 1.5 associates with it a Diophantine envelope. Diophantine envelopes were implicitly used in the proof of stability
of free and hyperbolic groups, and found other applications and different constructions that enabled one to construct canonical ones.
Starting with the Diophantine
envelope of a definable set,
we further associate with it a $Duo$ $Envelope$, that were also used implicitly in the proof of stability, and are the main tool that we use
in classifying
imaginaries over free and hyperbolic groups.

\vglue 1pc
\proclaim{Theorem 1.7} Let $F_k$ be a non-abelian free group, let $L(p,q)$ be a definable set  
over $F_k$, and let $G_1(z,p,q,a),\ldots,G_t(z,p,q,a)$ be a finite collection of graded completions, that form a 
Diophantine envelope of $L(p,q)$ (see theorem 1.5).
There exists a finite collection of duo limit groups: 
$$Duo_1(d_1,p,d_2,q,d_0,a),\ldots,Duo_r(d_1,p,d_2,q,d_0,a)$$ 
which is  associated with
$L(p,q)$, that we call the $Duo$ $Envelope$ of $L(p,q)$, for which:
\roster
\item"{(1)}" For each index $i$, $1 \leq i \leq r$, there exists a (duo) test sequence of the duo limit group
$Duo_i$, that restricts to a sequence of couples, $\{(p_n,q_n)\}$, so that for every index $n$,
$(p_n,q_n) \in L(p,q)$, and the test sequence converges into the Duo limit group itself (and not to a proper quotient of it). 

A (duo) test sequence of a duo limit group is a sequence of values that factor through the Duo
limit group that restricts to test sequences of the two completions from which it is composed. Furthermore, the sequence of values of the duo test sequence
converges into a limit group that has the same structure as the Duo limit group, only the base $<d_0>$ may be replaced by a  quotient.



\item"{(2)}" Given a specialization  $(p_0,q_0) \in L(p,q)$,   
there exists an index $i$, $1 \leq i \leq r$,  and a duo family of the duo limit group $Duo_i$, so that
$(p_0,q_0)$ extends to a specialization that factors through the duo family, and there exists a test
sequence of the duo family that restricts to specializations $\{(p_n,q_n)\}$, so that for every index $n$,
$(p_n,q_n) \in L(p,q)$.

\item"{(3)}" Let $Duo$ be some duo limit group of $L(p,q)$.  
Every duo family of $Duo$ that admits a duo test sequence 
that restricts to a
sequence of specializations, $\{(p_n,q_n)\}$, for which $(p_n,q_n) \in L(p,q)$ for every index $n$, 
is boundedly covered
by the Duo envelope, $Duo_1,\ldots,Duo_r$. 

i.e., there exists some constant $b>0$ (that depends  on $Duo$ and the definable set $L(p,q)$),
for which  given a duo family of the duo limit group
$Duo$, and a duo test sequence in the duo family that restrict to elements $\{(p_n,q_n)\}$ that are all in  $L(p,q)$,
there exist at most $b$ duo families of the Duo envelope,
$Duo_1,\ldots,Duo_r$, 
such that all but finitely many of the pairs $\{(p_n,q_n)\}$ can be extended to specializations that belong to one of the $b$ duo families of
the Duo limit group $Duo_1,\ldots,Duo_r$.

Furthermore, if we fix a duo family of the duo limit group $Duo$, then at most $b$ duo families of the duo limit groups $Duo_1,\ldots,Duo_r$
satisfy the above condition for all the test sequences in the (fixed) duo family of $Duo$ that restrict to pairs $\{(p_n,q_n)\}$ that are in $L(p,q)$. 
\endroster
\endproclaim

\nfp In proving stability of the theory of a free group in [Se9], in studying the equationality of  Diophantine sets (section 2 in [Se9]), 
or in proving that the existence of rigid or strictly solid values of a rigid or solid limit group is in the Boolean algebra of equational sets
(section 3 in [Se9]), 
it was immediate to obtain a canonical collection of duo limit groups
that cover any duo family that is associated with a Diophantine set or with a rigid or a solid
limit group (see section 3 in [Se9]). When we study a general definable set, $L(p,q)$, 
the construction of such a universal
family of duo limit groups as stated in theorem 1.7  is more involved. 

We will construct the  duo envelope by starting with the same construction that was applied in constructing the duo limit groups that are
associated with a rigid or a solid limit group 
in section 3 of [Se9],
and then iteratively refine this construction using the sieve method [Se6],
in a similar way to the construction of the Diophantine envelope  (theorem 1.5).

Let $G_1,\ldots,G_t$ be a collection of graded completions that forms a Diophantine envelope of the definable set $L(p,q)$ (i.e., that satisfies the
conclusions of theorem 1.5). $G_1,\ldots,G_t$ are graded w.r.t.\ the parameter subgroup $<q,a>$.
We start with the graded completions, $G_1,\ldots,G_t$,
in parallel.
With each graded completion $G_j$, $1 \leq j \leq t$,  we first associate a finite collection  of duo
limit groups (that are not yet necessarily part of the duo envelope).

Each of the duo limit groups that we construct is going to have a closure of one of the graded completions, $G_j$, as one of the completions 
in the duo limit group, the one that contains the subgroup $<p,a>$, and we gradually construct the other completions in these duo limit groups, the ones 
that contain the subgroup $<q,a>$. The construction of these duo limit groups uses the construction of formal graded limit groups in section 3 of [Se2] and
(as the construction of the Diophantine envelope) the sieve procedure that is presented in [Se6].

To construct these (initial) duo limit groups, we look at the entire collection
of graded test sequences that factor through the given graded completion, $G_j$, for which:
\roster
\item"{(i)}"  the 
(restricted) sequence of specializations $\{p_n\}$ can be extended to specializations of one of the 
graded completions (the Diophantine envelope) $G_1(z,p,q,a)\ldots,G_t(z,p,q,a)$. 

\item"{(ii)}"  the values of the variables
$q$ in the extension to specializations of one of the graded completions, $G_1,\ldots,G_t$, are fixed for the entire test sequence.
\endroster

With this entire collection of graded test sequences, and their extensions to specializations of
$G_1,\ldots,G_t$,
we associate finitely many  graded Makanin-Razborov diagrams, precisely as we did in constructing
the formal graded Makanin-Razborov diagram in section 3 of [Se2]. As in the
formal Makanin-Razborov diagram, each (graded formal) resolution in the diagrams that we construct terminates with
a (graded) closure of the given graded completion, $G_j$, that we have started with, 
amalgamated with another group (that contains the subgroup $<q,a>$) along a group that is generated by (image of) the terminal graded limit group
of the graded closure of $G_j$
(the terminal limit group of the graded closure of $G_j$ is either rigid or solid), and the abelian vertex groups
that commute with non-trivial elements in the terminal limit group (i.e., abelian vertex groups that have their pegs in the base terminal group of 
the graded closure of $G_j$).
  
We continue by associating a graded Makanin-Razborov diagram with the terminal group of each of the obtained groups. i.e., we look at a graded 
Makanin-Razborov diagram of the group that contain the image of the terminal rigid or solid  group of the graded completion $G_j$ that we
started with, relative to this subgroup. When we join the completion of such a graded resolution,  the other higher levels of the graded closure
of the graded completion $G_j$, i.e. when we amalgamate these two completions along their bases, we get a duo limit group.

This is the preliminary collection of duo limit groups. Note that  duo test sequences in one of these duo limit groups do not restrict to pairs $(p_{\ell},q_{\ell})$
that are in the definable set $L(p,q)$. To construct duo limit groups that do have duo test sequences that restricts to pairs $(p_{\ell},q_{\ell})$ that
are in the definable set $L(p,q)$,  
we proceed by constructing an iterative procedure that is similar to the one that is used in the proof of theorem 1.5 to
construct the Diophantine envelope, 
and is based on the sieve procedure [Se6]. 

At each step we first
collect all the duo test sequences of the current duo limit group
(that is associated with $G_j$),
for which for every specialization from the sequence, that restricts to a proof statement from one of the graded completions, $G_1,\ldots,G_t$,  
there exists an extra rigid or strictly 
solid specialization. 

We analyze
the obtained collection of sequences of specializations using the analysis of (graded) formal limit 
groups that appears in section 3 in [Se2], in the same way that we constructed the collection of preliminary duo limit groups that are associated with 
the graded completion $G_j$. Therefore,
with these test sequences  and their Zariski closure we associate finitely many closures of the duo limit group
that we started with. 

As in the sieve procedure, and in the procedure that was used to prove theorem 1.5,  we look at all the sequences of specializations  that factor through 
the constructed closures of the duo limit groups that we started with, such that:
\roster
\item"{(1)}" The sequences restrict to a test sequence of the completion that is associated with the original graded completion $G_j$, i.e., the 
completion that contains the subgroup $<p,a>$ in the closure of the duo limit group.

\item"{(2)}" 
The extra 
rigid or strictly solid specializations that were added when we constructed the closures of the duo limit group collapses. i.e., they satisfy one of finitely
many Diophantine conditions that demonstrate that either they are nor rigid or not strictly solid, or that they belong to the same family as one of the
rigid or strictly solid family in the corresponding proof statement. 
\endroster

We analyze
the obtained collection of specializations that satisfy conditions (1) and (2) using the analysis that is used in the general step of the sieve procedure,
and finally continue to the next step only with those duo limit groups for which the graded completion in them that in the duo limit group that contains the
subgroup $<q,a>$, i.e., the graded completion that is paired with the graded completion that contains the subgroup $<p,a>$ and has the same structure
as the graded completion $G_j$ we have started with, has strictly smaller complexity (in the sense of the sieve procedure in [Se6]) 
than the corresponding graded completion in the duo limit group that we started this step with, and from which the current duo limit group was
constructed from.
 
We continue iteratively, and by the termination of the sieve procedure
(theorem 22 in [Se6]), the obtained procedure terminates after finitely many steps.

Using this procedure, we obtain finitely many duo limit groups. Given one of the obtained duo limit groups we look at all the duo test sequences
that factor through it and restrict to values $(p_{\ell},q_{\ell}) \in L(p,q)$. With the collection of these duo test sequences we naturally associate 
a finite collection of duo limit groups, where each of these duo limit groups is a quotient of the one we started with.

We denote the obtained finite collection of duo limit groups:   $Duo_1,\ldots,Duo_{r_1}$. This collection of duo limit groups is associated with the graded
resolutions $G_1,\ldots,G_t$ that are the Diophantine envelope of the definable set $L(p,q)$. They satisfy property (1) in theorem 1.7 but in general
they do not property (2).

To get a family that satisfies properties (2) and (3) as well,
we need to extend the collection of duo limit groups. To get this extension we combine the constructions of the diagrams that
are used in sections 2 and 3 in [Se9] (for studying equational properties of Diophantine and rigid and solid limit groups), 
with the sieve procedure that was used in constructing the duo envelope and the initial collection of duo limit groups $Duo_1,\ldots,Duo_{r_1}$.

\medskip
We start with an iterative procedure that is similar to the one that we used in sections 2 and 3 in [Se9], that enables one to associate more graded
completions, that are not part of the Diophantine envelope (theorem 1.5), and are associated with finite intersections of definable sets:
$\Cap_{i=1}^m \, L(p,q_i)$ for finite subsets of values $q_1,\ldots,q_m$ of the variables $q$, where $m$ is an arbitrary positive integer.

To construct the finite collection of graded closure we construct iteratively a finite diagram, where at each vertex in the finite diagram we place a 
graded closure w.r.t.\ finite values of the variables $q$.
We start with each of the graded resolutions, $G_1,\ldots,G_t$, that form the Diophantine envelope, in parallel. 
With each graded resolution $G_j$, $1,\ldots,t$, we associate a finite diagram. With each vertex in this diagram, on which we place
a graded completion we later
associate a finite collection of duo limit groups, in a similar way to the construction of the duo limit groups,
$Duo_1,\ldots,Duo_{r_1}$, that were associated with the initial graded completions (the Diophantine envelope)
$G_1,\ldots,G_t$. The union of these finite collections of finite duo limit groups are going to be the duo envelope of the definable set
$L(p,q)$. 

We start with a graded completion  $G_j$ from the Diophantine envelope. $G_j$ is a graded completion with respect to the parameter subgroup $<q,a>$, that
we denote $<q_1,a>$. We look at all the values of the pair $(p,q_2)$ in L(p,q) for which:
\roster
\item"{(i)}" $q_1 \neq q_2$ 

\item"{(ii)}" $(p_0,q_1)$ extends to a specialization of $G_j$ such that the the extended specialization restricts to a valid
proof statement that $(p_0,q_1) \in L(p,q)$. 

\item"{(iii)}"  $(p_0,q_2)$ extends to a specialization of one of the graded completions,
$G_1,\ldots,G_t$, and the extended specialization restricts to a valid proof statement that $(p_0,q_2) \in L(p,q)$.
\endroster

We analyze these collections of (combined) specializations of $G_j$ and each of the graded resolutions $G_1,\ldots,G_t$, using the analysis of 
quotient limit groups and their (graded) resolutions  that is used
in the first step of the sieve procedure in [Se6]. Each of the graded resolutions that are obtained in this analysis has complexity 
(the one that is used in the sieve procedure in [Se6])
that is bounded above by the complexity of the original completion $G_j$.

We continue iteratively as what we did in constructing the Diophantine envelope in theorem 1.5. At each step we construct finitely many graded closures
that are associated with the graded resolutions and completions that were constructed in the previous step. As in the proof of theorem
1.5 we first look at all the graded test sequences  in which in one of the associated proof statements (the ones that are associated with the pairs $(p,q_1)$
or $(p,q_2)$) an element that is supposed to be rigid or strictly solid is not, or an element that is supposed to be a primitive root is not. With this
collections of test sequences we associate finitely many graded closures. We also look at test sequences for which 
there exists an extra rigid or strictly solid value that does not belong to one of the families that are part of the two proof statements.  
With these test sequences associate finitely many graded closures that we call (as in the proof of theorem 1.5 and in the sieve procedure, $ExtraPS$). 

We continue iteratively to the next step only with the ExtraPS resolutions, on which we impose one of finitely many Diophantine conditions
that demonstrate that the potential extra rigid or solid value is in fact nor rigid nor strictly solid, or that it belongs to a family that
is part of the proof statement.

We collect the specializations of the ExtraPS resolutions with the values of the elements that demonstrate the collapse of the extra rigid or 
strictly solid value in finitely many
limit groups and analyze these limit groups using the analysis that appears in the general step of the sieve procedure  in [Se6]
The procedure terminates after finitely many steps by theorem 22 in [Se6] (the termination of the sieve procedure).

Once the first procedure terminates we are left with finitely many graded completions that were constructed along the procedure, and with each 
graded completion there are finitely many associated graded closures. The graded completions that are not graded closures of the original graded completion $G_j$,
for the graded completions that appear in the second level of the diagram that we construct. i.e., the graded completions that are placed on vertices with
a directed edge from the vertex in which $G_j$ is placed. 

The graded completions that we constructed are associated with pairs of values $q_1,q_2$ of the variables $q$. We continue iteratively, and add a value
$q_3$ of the variables $q$. Given one of the graded completions $G^1_i$ that we constructed, and appears in the second level of the diagram, we look at
values that factor through it, that restrict to valid proof statements that both $(p_0,q_1)$ and $p_0.q_2)$ are in $L(p,q)$, 
and such that $(p_0,q_3)$  further extends to a specializations of one of the original graded
completions $G_1,\ldots,G_t$, and that specialization restricts to valid proof statement that demonstrates that $(p_0,q_3) \in L(p,q)$.

With the collection of specializations that are associated with the triple $q_1,q_2,q_3$ we associate a finite collection of graded completions
using the analysis in the first step of the sieve procedure. Then we use an iterative procedure similar to the one that we applies for the pair
$(q_1,q_2)$, and associate with the triple $(q_1,q_2,q_3)$ a finite collection of graded completions, where with each  completion we associate a 
finite collection of graded closures. The graded completions are placed on vertices in the third level of the constructed diagram.

We continue iteratively. In each step we add a new value $q_n$ of the variables $q$ and use the sieve procedure to associate with the $n$ values of
the variables $q$, a finite collection of graded completions, and with each completion a finite collection of graded closures. These are the
graded completions that are placed on vertices in the $n$-th level of the constructed directed diagram.

The iterative process that iteratively adds a value of the variables $q$, $q_1,q_2,\ldots$ and constructs new graded completions terminates after finitely
many steps by the termination of the sieve procedure (theorem 22 in [Se6]). This termination, and the whole construction of the diagram is similar to the 
termination of the constructions of the diagrams  in section 2 and 3 of [Se9], diagrams that were used  to prove equationality properties of 
Diophantine sets and of rigid and strictly
solid homomorphisms that we generalize here to general Diophantine sets.

\smallskip
Once we have constructed all the graded completions and form the directed diagram, we start with each of the graded completions in parallel,
and associate with it a finite collection of duo limit groups, precisely as we associated the duo limit groups $Duo_1,\ldots,Duo_{r_1}$ with the
graded completions, $G_1,\ldots,G_t$, that were assumed to be a Duo envelope of the set $L(p,q)$, and are placed in the first level of the diagram
that we constructed. We denote the union of the collections of duo limit groups that are associated with all the graded completions in the diagram,
$Duo_1,\ldots,Duo_r$.

By the iterative procedures that were used to construct them, the collection of Duo limit groups, $Duo_1,\ldots,Duo_r$, satisfies properties (1) and (2) 
in theorem 1.7. To prove part (3) let $Duo$ be a duo limit group that has duo test sequences that restrict to pairs $\{(p_n,q_n)\}$ that are in the
definable set $L(p,q)$. i

We look at the collection of duo test sequences of $Duo$ that restrict to pairs $\{(p_n,q_n)\}$ that are in the definable
set $L(p,q)$. By property (2) of the Duo envelope, $Duo_1,\ldots,Duo_r$, each pair $(p_n,q_n)$ in these test sequences extends to a specialization
of one of these duo limit groups. Using the techniques that construct formal limit groups in section 3 of [Se2], from the collection of test sequences
of $Duo$ it is possible to construct finitely many duo closures of the duo limit group $Duo$, that we denote $Cl_i(Duo)$, with (duo) maps:
$\eta_i: Duo_{j(i)} \to Cl_i(Duo)$, i.e., the maps $\eta_i$ send each of the two completions in $Duo_{j(i)}$ into the corresponding completions in 
$Cl_i(Duo)$, and the elements $p$ and $q$ in $Duo_{j(i)}$ to the elements $p$ and $q$ in correspondence in $Cl_i(Duo)$  
(recall that a duo closure is composed from graded closures of the two completions from which the duo limit group is
built, amalgamated along their common (terminal) base group, which is rigid or solid w.r.t.\ the base group of the (image of) the base of original
duo limit group). 

The existence of the closures $Cl_i(Duo)$ and the maps $\eta_i$ imply that there exists a positive integer $b$, such that for a given  duo family
of $Duo$, the collection of duo test sequences in this duo family that restrict to pairs $(p_n,q_n)$ in $L(p,q)$, is covered by at most $b$ duo families
of $Duo_1,\ldots,Duo_r$. i.e., from every duo test sequence in the given duo   family of $Duo$  that restrict to pairs $\{p_n,q_n)\}$ in $L(p,q)$, it is possible
to omit finitely many elements, and the rest of the sequence extend to specializations of the $b$ duo families of $Duo_1,\ldots,Duo_r$. This proves part
(3) of theorem 1.7.   

\line{\hss$\qed$}

\vglue 1.5pc
\centerline{\bf{\S2. Few Basic Imaginaries}}
\medskip
Our goal in this paper is to study definable equivalence relations over free and hyperbolic groups.
Before we analyze the structure of a general definable equivalence relation over such groups,
we introduce some basic well-known and less known definable equivalence relations over a free (or a hyperbolic) group,
and prove, using Diophantine envelopes that were presented in the first section, 
that these equivalence relations are  imaginaries (non-reals).  
Later on we show that if one adds these basic imaginaries as new sorts to the model of a free or a hyperbolic
group, then for any definable equivalence relation there exists a definable multi-function that separates classes
and maps every equivalence class into a uniformly bounded set, i.e., that adding the basic imaginaries
as sorts geometrically eliminates imaginaries over free and hyperbolic groups.
We start by proving that conjugation is not real.

\vglue 1pc
\proclaim{Remark} It was pointed out to us that in [Sk] a stronger result than theorem 2.1 below is proved using
an easier argument. Still we preferred to include our proof since it presents our perspective on how to prove non-definability
using the geometric description of definable sets from [Se6] (theorem  1.6) and the existence of a Diophantine envelope (theorem 1.5).
\endproclaim

\vglue 1pc
\proclaim{Theorem 2.1} Let $F_k=<a_1,\ldots,a_k>$ be a non-abelian free group, and let:
$$Conj(x_1,x_2) \ = \ \{ \, (x_1,x_2) \, | \, \exists u \ ux_1u^{-1}=x_2 \, \}$$
be the definable equivalence relation that is associated with conjugation over $F_k$. Then $Conj(x_1,x_2)$ is not real.
\endproclaim

\nfp To prove that conjugation is not real, we need to show that there is no positive integer $m$
and a (definable) function:  
$f:F_k \to F_k^m$ that maps each conjugacy class in $F_k$ to a tuple, and distinct conjugacy classes to
distinct tuples in $F_k^m$. To prove that there is no such function,
we use the precise (geometric) description of a definable set from theorem 1.6, 
and the Diophantine envelope of
a definable set that was constructed in theorem 1.5.


Suppose that conjugation is real. Then there exists a definable function: $f:F_k \to F_k^m$,
that  maps each conjugacy class to a tuple, and different conjugacy classes to distinct tuples. Let
$L(v,x_1,\ldots,x_m)$ be the graph of the function $f$, that by our assumption on the definability of $f$
has to be a definable set, for which:

\roster
\item"{(i)}" for every possible value of the variable $v$ there exists a unique value
of the tuple,  $(x_1,\ldots,x_m)$.  

\item"{(ii)}" the value of the tuple $(x_1,\ldots,x_m)$ is the same for elements $v$ in the same conjugacy class, 
and distinct for different conjugacy classes.   
\endroster

Using theorem 1.5, we associate with the definable set $L(v,x_1,\ldots,x_m)$, its Diophantine envelope
with respect to the parameter $v$ (i.e., in the statement of theorem 1.5, we set $v$ to be $q$, and 
$(x_1,\ldots,x_m)$ to be $p$). Let $G_1,\ldots,G_t$ be the graded completions that form this Diophantine envelope (note that by theorem 1.5
each of the graded completions $G_1,\ldots,G_t$ contains the coefficient subgroup $F_k$ in a distinguished vertex group (coefficient subgroup that is a subgroup
of the terminal vertex group in the graded completions).

Recall that given a test sequence of a graded completion, all the elements that are not in the terminal vertex group of the graded completion,
and if the terminal vertex group is solid, all the elements that are not in the distinguished vertex group in the graded JSJ decomposition of the
terminal vertex group, obtain infinitely many distinct values along the test sequence.

Since for every possible value of $v$ there exists a unique value of
$x_1,\ldots,x_m$ in $L(v,x_1,\ldots,x_m)$, each of the graded completions, $G_j(z,x_1,\ldots,x_m,v,a)$, is either
a rigid limit group (with respect to the parameters $<v,a>$), or it is a solid limit group where the subgroup 
$<x_1,\ldots,x_m,v,a>$ is contained in the distinguished vertex group (that is stabilized by
$<v,a>$).

With the definable set $L(v,x_1,\ldots,x_m)$ we
associate three sequences of specializations:
$\{u_n\}$, $\{(v_n,x_1^n,\ldots,x_m^n)\}_{n=1}^{\infty}$ and 
$\{(\hat v_n,\hat x_1^n,\ldots,\hat x_m^n)\}_{n=1}^{\infty}$,  
so that:
\roster
\item"{(1)}" for every index $n$, the tuples $(v_n,x^n_1,\ldots,x^n_m)$ and
$(\hat v_n,\hat x_1^n,\ldots,\hat x_m^n)$,  
are in the definable set $L(v,x_1,\ldots,x_m)$. 

\item"{(2)}" $v_n$ and $u_n$ are taken from a test sequence in the free group $F_k$ (see theorem 1.1 and 
lemma 1.21 in [Se2] for a test sequence), and $\hat v_n=u_n  v_n u_n^{-1}$.

\item"{(3)}" $2 \cdot |u_n| \, > \, |v_n| \, > \, \frac {1}{2}  \cdot |u_n|$, 
where $|w|$ is the length of the word 
$w \in F_k$ with respect to a fixed generating set of $F_k$. 
\endroster

By theorem 1.5, each of the specializations 
$(v_n,x^n_1,\ldots,x^n_m)$ extends to a specialization that factors through one of the graded completions,
$G_1,\ldots,G_t$, that form the Diophantine envelope of the definable set $L(v,x_1,\ldots,x_m$). 
By passing to
a subsequence, and changing the order of the graded completions, we can assume that they all factor through
the graded completion $G_1$. 

$\{v_n\}$ was chosen to be a test sequence in $F_k$, This implies that the sequence, $\{(v_n,a_1,\ldots,a_k)\}$, converges into an action of
the free product $<v>*F_k$ on a limit tree, $T_v$, which is a simplicial tree, with an associated graph of groups that has a single vertex
stabilized by $F_k$, and a single loop with trivial stabilizer and a Bass-Serre element which is the element $v$.

Since $G_1$ is rigid or solid with respect to the parameter subgroup $<v,a>$, and the homomorphisms, $h_n: <v,a> \to F_k$, $h_n(v)=v_n$, extend to
homomorphisms of the envelope $G_1$, $\hat h_n: G_1 \to F_k$, that we can choose shortest possible if $G_1$ is solid w.r.t.\ $<v,a>$ (and just a rigid
extension if $G_1$ is rigid),  the homomorphisms $\hat h_n$ subconverge into an action of $G_1$ on the limit tree $T_v$. Hence, the homomorphisms 
$\hat h_n$ subconverge into a limit action of the limit group $<v>*F_k$ on $T_v$ with a homomorphisms $\eta:G_1 \to <v>*F_k$. Furthermore: 
\roster
\item"{(1)}" By the properties of the envelope (theorem 1.5), and from the existence of the map $\eta$ that is derived from a subsequence of a 
sequence of values, $\{(v_n,x_1^n,\ldots,x_m^n)\}$,
by passing to the subsequence, and still use the index $n$, for each index $n$, there is a retraction, $\nu_n:H=<v,a> \to <a>=F_k$, given by: 
$\nu_n(v)=v_n$, and for every index $i$, $1 \leq i \leq m$, $x_i^n=\nu_n \circ \eta(x_i)$.

\item"{(2)}" By the construction of the Diophantine envelope in the proof of theorem 1.5, the graded completion $G_1$  contains a proof statement.
i.e.,  elements
that together are supposed to 
validate that the elements $(v,x_1,\ldots,x_m)$ are indeed elements of the definable set
$L(v,x_1,\ldots,x_m)$ (see the proof of theorem 1.5). 
For each index $n$, the restriction of the composition $\nu_n \circ \eta : G_1 \to F_k$ 
to these elements validates that $(v_n,x_1^n,\ldots,x_m^n) \in L(v,x_1,\ldots,x_m)$. 
\endroster

\vglue 1pc
\proclaim{Lemma 2.2}  Recall that $\{u_n\}$ and $\{v_n\}$ are taken from two test sequences, and:
$\hat v_n=u_nv_nu_n^{-1}$. If we set
$\hat \nu_n:H \to F_k$ to be the retraction given by: $\hat \nu_n(v)=\hat v_n=u_nv_nu_n^{-1}$, 
then after possibly throwing finitely many elements from the two test sequences, for every index $n$ and every index $i$, 
$1 \leq i \leq m$: $\hat x_i^n=\hat \nu_n \circ \eta (x_i)$. 
\endproclaim

\nfp The sequence $\{(u_n,v_n)\}$ is composed from two test sequences. Hence, the sequence of homomorphisms $s_n:<u,v,a> \to F_k$, 
$s_n(u)=u_n$, $s_n(v)=v_n$, converges into a free action of the free  limit group, $L_{uv}=<u>*<v>*F_k$ on a Bass-Serre tree $T_{uv}$, 
with a graph of groups decomposition that has a single vertex stabilized by $F_k$ and two loops with trivial stabilizers, and Bass-Serre generators 
$u$ and $v$.

We define an embedding: $\mu:H \to L_{uv}$, $\mu(a)=a$ and $\mu(v)=uvu^{-1}$. Hence, we get a map $\tau_{uv}: G_1 \to L_{uv}$ given by
$\tau_{uv}=\mu \circ \eta$. Hence, we get the map $s_n \circ \tau_{uv}:G_1 \to F_k$. Since $G_1$ contains elements for proof statements
that values of $(v,x_1,\ldots,x_m)$ are in the definable set $L(v,x_1,\ldots,x_m)$, the images $s_n \circ \tau_{uv}(G_1)$ provide proof statements that
the elements $s_n \circ \tau_{uv}(v,x_1,\dots,x_m)$ are in $L(v,x_1,\ldots,x_m)$. 

We argue that for large $n$ these are valid proof statements, using the fact that for large $n$ the maps $\nu_n \circ \eta$ restrict to valid
proof statements. A proof statement contains several rigid and strictly solid values of certain rigid
and solid limit groups w.r.t.\ the parameter subgroup $<v,x_1,\ldots,x_m,a>$ (that we denoted $Term_j$ in theorem 1.6). Since $\mu$ is an embedding 
for large $n$ the maps of every rigid limit group in the proof statement has to be rigid, since it satisfies the finite set of equalities and inequalities that are
required from a rigid value.  

Suppose that the images of elements that are suppose to be strictly solid are in fact not strictly solid for some subsequence of indices. In that case
we can add elements that will demonstrate that the images of that fake strictly solid element are in fact not strictly solid. By passing to a limit
of a subsequence, and using the
construction of formal limit groups section 3 of [Se2], the elements that demonstrate the failure of the images to be strictly solid can be chosen to be in
$L_{uv}$.

Since $uvu^{-1}$,
$u$ and free basis of $F_k$ is a free basis of $L_{uv}$, and $\tau_{uv}(G_1)$ is in $<uvu^{-1}>*F_k$ we can choose the elements that demonstrate the
failure of the images to be strictly solid to be in $<uvu^{-1}*F_k$. But this implies that for all $n$, $\nu_n \circ \eta(G_1)$ is not a valid proof statement
either, since an element that is supposed to be strictly solid is not, a contradiction. 

We further need to argue that for large $n$ there are no extra rigid or strictly solid values for the proof statements that the maps 
$s_n \circ \tau_{uv}$ restrict to. First, suppose that for some subsequence of indices there exist an extra rigid value, By possibly passing to a subsequence
that converges to a limit, the extra rigid values converge to an element in $L_{uv}$, and since $uvu^{-1}*F_k$ is a free factor in $L_{uv}$, the extra rigid value is in
$uvu^{-1}*F_k$. Hence, the extra rigid value is in $\tau_{uv}(G_1)$, so it can be pulled back by the inverse of the embedding, $\mu:H \to L_{uv}$ to an extra rigid 
value in $H$. Hence, for large enough $n$, $s_n(H)$ contains extra rigid values, a contradiction since $s_n(H)$ were assumed to restrict to valid proof statements.

Suppose that for some subsequence of indices there are extra solid values for the proof statements that the maps $s_n \circ \tau_{uv}$ restrict to. As in the rigid case,
we further pass to a convergent subsequence, such that the extra solid values converge into elements in $\mu(H) < L_{uv}$, so they can be pulled back to elements in $H$.
The values $s_n(H)$ restrict to valid proof statements, so the values of the elements that are associated with the extra solid value in $H$ are mapped by $s_n$ to values
that are either not strictly solid or they are in the same strictly solid family with elements that are part of the proof statement that $s_n(H)$ restrict to.

Hence the extra solid elements in $s_n(H)$ collapse, and it is possible to demonstrate that collapse by extra elements that demonstrate that they
satisfy some Diophantine condition. By passing to a further subsequence and to a limit, the elements that demonstrate the collapse Diophantine condition are in fact in $H$.
Since $\mu$ embeds $H$ in $L_{uv}$, the collapse elements are mapped to $<uvu^{-1}>*F_k$. For large $n$, the elements that are suppose to extra strictly solid
(from the subsequence that converged to a limit) are given by a restriction of the map $s_n \circ \tau_{uv}$. But these elements collapse since they satisfy the collapse 
Diophantine condition in $\mu(H)$. 

Therefore, for large $n$ the proof statements that the images $s_n \circ \tau_{uv}(G_1)$ restrict to are valid proof statements and the lemma follows.

\line{\hss$\qed$}

Since the
elements $x_1,\ldots,x_m$ are distinct for distinct conjugacy classes of specializations of $v$, and the
elements $v_n$ are not conjugate for different indices, since they are taken from a test sequence, the 
tuples $(x_1^n,\ldots,x_m^n)$ are distinct for distinct indices $n$. Hence, there exists an index $i$, 
$1 \leq i \leq m$, for which $\eta(x_i) \notin F_k$. Therefore, for large enough $n$:
$$x_i^n=\nu_n \circ \eta (x_i) \,  \neq \, \hat \nu_n \circ \eta (x_i)= \hat x_i^n.$$ 
But, $v_n$ is conjugate to $\hat v_n$, and both tuples $(v_n,x_1^n,\ldots,x_m^n)$ and
$(\hat v_n,\hat x_1^n,\ldots,\hat x_m^n)$ are contained in the definable set $L(v,x_1,\ldots,x_m)$. Hence,
for every index $n$, and every $i$, $1 \leq i \leq m$, $x_i^n=\hat x_i^n$, and we get a contradiction.

\line{\hss$\qed$}

\vglue 1pc
\proclaim{Remark} Note that the argument that was used to prove theorem 2.2 implies that there is no
no definable multi-function, $f:F_k \to F_k^m$ that sends each conjugacy class to a finite set (that has to be uniformly bounded
by the existence of the envelope in theorem 1.5), and separates between conjugacy classes. 
\endproclaim

Having proved that conjugation is not real, we argue that left (and right) cosets of cyclic groups are real.

\vglue 1pc
\proclaim{Lemma 2.3} Let $F_k=<a_1,\ldots,a_k>$ be a non-abelian free group, let $m$ be a 
positive integer and let:
$$Left(x_1,y_1,x_2,y_2) \ = \ \{ \, (x_1,y_1,x_2,y_2) \ | \ y_1,y_2 \neq 1 \wedge [y_1,y_2]=1 \wedge
\exists y \, [y,y_1]=1 \wedge x_1^{-1}x_2=y^m \, \}$$
be the definable equivalence relation associated with left 
cosets of cyclic 
subgroups over $F_k$. 

Note that $Left$ is not defined on all of $F_k^2$, but it can be naturally (definably) extended to $F_k^2$ by defining all the pairs,
$\{(x,1) \ x\in F_k\}$, to be in the same equivalence class in $F_k^2$. Then $Left$ 
is not  a real definable equivalence relation for any fixed positive integer $m$.
\endproclaim

\nfp The proof that we give is similar to the one that we used in theorem 2.1.
Suppose that Left is real. Then there exists a definable function: $f:F_k^2 \to F_k^r$,
that maps each left coset (of the corresponding cyclic subgroup) to a tuple, and different left cosets 
to distinct tuples. Let
$L(t,s,x_1,\ldots,x_r)$ be the definable set that is associated with the definable function $f(t,s)$. Note that
for every possible specialization  $(t_0,s_0)$ of $(t,s)$, there exists a unique specialization
$(t_0,s_0,x_1^0,\ldots,x_r^0)$ that belongs to the set $L(t,s,x_1,\ldots,x_r)$, and if 
$(t_0,s_0,x_1^0,\ldots,x_r^0) \in L(t,s,x_1,\ldots,x_r)$ then
$(t_0s_0^{\ell  m},s_0^v,x_1^0,\ldots,x_r^0) \in L(t,s,x_1,\ldots,x_r)$ for every non-zero integers $\ell,v$.
   
We proceed as in the proof of theorem 2.1. Using theorem 1.5, we associate with the definable set 
$L(t,s,x_1,\ldots,x_m)$, its Diophantine envelope
with respect to the parameter subgroup $<t,s,a>$ (i.e., in the statement of theorem 1.5, we set $t,s$ to be $q$, and 
$x_1,\ldots,x_r$ to be $p$). Let $G_1,\ldots,G_t$ be the graded completions that form this Diophantine envelope (note that each of
the graded completions $G_j$ contains a copy of the coefficient group $F_k$ in its terminal rigid or solid vertex group).
Since for every possible value of the pair $(t,s)$ there exists a unique value of
$x_1,\ldots,x_r$ in $L(t,s,x_1,\ldots,x_r)$, each of the graded completions, $G_j(z,x_1,\ldots,x_r,t,s,a)$, 
is either
a rigid limit group (with respect to the parameter subgroup $<t,s,a>$), or it is a solid limit group where 
the subgroup 
$<x_1,\ldots,x_r,t,s,a>$ is contained in the distinguished vertex group (that is stabilized by
$<t,s,a>$).

With the definable set $L(t,s,x_1,\ldots,x_r)$ we
associate two sequences of specializations:
$\{(t_n,s_n,x_1^n,\ldots,x_r^n)\}_{n=1}^{\infty}$ and 
$\{(\hat t_n,\hat s_n,\hat x_1^n,\ldots,\hat x_r^n)\}_{n=1}^{\infty}$,  
so that:
\roster
\item"{(1)}" for every index $n$, the tuples $(t_n,s_n,x^n_1,\ldots,x^n_r)$ and
$(\hat t_n,\hat s_n,\hat x_1^n,\ldots,\hat x_r^n)$,  
are in the definable set $L(t,s,x_1,\ldots,x_r)$. 

\item"{(2)}"  Let $u_n$ and $v_n$  be taken from   test sequences in the free group $F_k$ (see theorem 1.1 and 
lemma 1.21 in [Se2] for a test sequence). 

\item"{(3)}" $s_n=v_n$ and $t_n=u_n$,
 $\hat s_n=v_n$ and $\hat t_n=u_n{(v_n)}^{m \cdot n }$.

\item"{(4)}" $2 \cdot |u_n| \, > \, |v_n| \, > \, \frac {1}{2}  \cdot |u_n|$, 
where $|w|$ is the length of the word 
$w \in F_k$, with respect to a fixed generating set of $F_k$. 
\endroster

By theorem 1.5, each of the specializations 
$(t_n,s_n,x^n_1,\ldots,x^n_r)$ extends to a specialization that factors through one of the graded completions,
$G_1,\ldots,G_t$, that form the Diophantine envelope of the definable set $L(t,s,x_1,\ldots,x_r)$. By passing to
a subsequence, and changing the order of the graded completions, we can assume that they all factor through
the graded completion $G_1$. 

Since $G_1$ is rigid or solid with respect to the parameter subgroup $<t,s,a>$, and the homomorphisms, $h_n: <t,s,a> \to F_k$, $h_n(t)=t_n$, 
$h_n(s)=s_n$, extend to
rigid or strictly solid homomorphisms of the envelope $G_1$, $\nu_n: G_1 \to F_k$, that we can choose shortest possible if $G_1$ is solid w.r.t.\ $<t,s,a>$ (and just a rigid
extension if $G_1$ is rigid). 

Since the elements $u_n,v_n$ were taken from a test sequence, the sequence of homomorphisms $\{h_n\}$ converges into a faithful action of the group
$H=<t>*<s>*F_k$ on a simplicial tree $T_{ts}$, with a Bass-Serre tree that consists of a single vertex stabilized by $F_k$ and two edge with
trivial stabilizers and Bass-Serre elements $u$ and $v$. Since the homomorphisms $\nu_n$ were chosen to be shortest possible, after possibly passing
to a subsequence, the homomorphisms $\{\nu_n\}$ subconverge to a map $\eta: G_1 \to H$.

We further define $\hat h_n: <u,t,s,a> \to F_k$ by setting: $\hat h_n(u)=u_n$, $\hat h_n(t)= \hat t_n$, $\hat h_n(s)=\hat s_n=s_n$.
The sequence of homomorphisms $\{\hat h_n\}$ converges into a limit group $\hat H=(<u>*<v>*F_k)*_{<v>}<v,f>$, where $<v,f>$ is isomorphic to $Z^2$,
$s$ is mapped to $v$ and $t$ is mapped to $uf$. By Noetherianity, for large $n$ the homomorphisms $\hat h_n$ can be identified with homomorphisms:
$\tilde h_n: \hat H \to F_k$. 

There is a natural embedding $\tau:H \to \hat H$ that maps the elements $s$ and $t$ in $H$ to their images in $\hat H$. This embedding enables
one to define a map $\hat \eta:G_1 \to \hat H$ as $\hat \eta=\tau \circ \eta$.  
%

\vglue 1pc
\proclaim{Lemma 2.4}  
For large $n$ the homomorphisms $\tilde h_n \circ \hat \eta$
restrict to valid proof statements that $\tilde h_n \circ \hat \eta (s,t,x_1,\ldots,x_r)$ are in the definable set: $L(s,t,x_1,\ldots,x_r)$. 
\endproclaim

\nfp The proof is similar to the proof of lemma 2.2.
We argue that for large $n$ these are valid proof statements, using the fact that for large $n$ the maps $h_n \circ \eta$ restrict to valid
proof statements. 
Since $\tau$ is an embedding 
for large $n$ the maps of every rigid limit group in the proof statement has to be rigid.

Suppose that the images of elements that are suppose to be strictly solid are in fact not strictly solid for some subsequence of indices. In that case
we  add elements that will demonstrate that the images of that fake strictly solid element are in fact not strictly solid. By passing to a limit
of a subsequence, and using the
construction of formal limit groups section 3 of [Se2], the elements that demonstrate the failure of the images to be strictly solid can be chosen to be in
$\hat H$.  $\tau(H)=<v>*<uf>*F_k \ < \ \hat H$ so $\hat H=\tau(H)*_{<v>} \ <v,f>$, where $<v,f>$ is isomorphic to $Z^2$, so there is a retraction from
$\hat H$ onto $\tau(H)$ that maps $<v,f>$ to $<v>$. 

Therefore, if strictly solid elements in the proof statements that are restriction of the maps: $\tilde h_n \circ \hat \eta$ are not strictly solid for
some subsequence of indices, then for large indices in this subsequence,  strictly solid  elements in the proof statements that are restrictions of the maps:
$h_n \circ \eta$ are not strictly solid, a contradiction to our initial assumptions.

We further need to argue that for large $n$ there are no extra rigid or strictly solid values for the proof statements that the maps 
$s_n \circ \tau_{uv}$ restrict to. First, suppose that for some subsequence of indices there exist an extra rigid value, By possibly passing to a subsequence
that converges to a limit, the extra rigid values converge to an element in $\tau(H)$.
Hence, the extra rigid value can be pulled back by the inverse of the embedding, $\tau:H \to \hat H$ to an extra rigid 
value in $H$, so the original proof statements that are restrictions of $h_n \circ \eta$ are not valid proof statements, a contradiction.

Suppose that for some subsequence of indices there are extra strictly solid values for the proof statements that the maps $tilde h_n \circ \hat \eta$ restrict to. 
Because of the retraction from $\hat H$ to $\tau(H)$ we can modify the elements that represent the extra strictly solid values to be in $\tau(H)$. Since
$\tau$ is an embedding,
we can pull these elements back to $H$. 

For large $n$ the proof statements that the maps: $h_n \circ \eta$ restrict to are valid proof statements. Hence, the elements that represent extra strictly solid 
element are either not strictly solid or they are in the same family of a strictly solid value that is in the proof statement. Therefore, these collapse 
forms of the
elements correspond to an extra solid value can be expressed by Diophantine conditions, and these Diophantine conditions can be expressed by elements in
$H$. By mapping $H$ to $\tau(H)$ it follows that the elements that are supposed to be associated with extra strictly solid values in
$\hat H$ are in fact not extra strictly solid values for large $n$, a contradiction, and the lemma follows. 

\line{\hss$\qed$}

Since the
elements $x_1,\ldots,x_m$ are distinct for distinct conjugacy classes of specializations of $v$, and the
elements $v_n$ are not conjugate for different indices, since they are taken from a test sequence, the 
tuples $(x_1^n,\ldots,x_m^n)$ are distinct for distinct indices $n$. Hence, there exists an index $i$, 
$1 \leq i \leq m$, for which $\eta(x_i) \notin F_k$. Therefore, for large enough $n$:
$$x_i^n=\nu_n \circ \eta (x_i) \,  \neq \, \hat \nu_n \circ \eta (x_i)= \hat x_i^n.$$ 
But, $v_n$ is conjugate to $\hat v_n$, and both tuples $(v_n,x_1^n,\ldots,x_m^n)$ and
$(\hat v_n,\hat x_1^n,\ldots,\hat x_m^n)$ are contained in the definable set $L(v,x_1,\ldots,x_m)$. Hence,
for every index $n$, and every $i$, $1 \leq i \leq m$, $x_i^n=\hat x_i^n$, and we get a contradiction.

\line{\hss$\qed$}

Since the
elements $x_1,\ldots,x_r$ are distinct for distinct left cosets,
and the pairs: $\{(t_n.s_n)\}$
belong to distinct left cosets, 
there exists an index $i$, 
$1 \leq i \leq r$, for which $\eta(x_i) \notin <s,a>$. Therefore, for large enough $n$:
$$x_i^n=h_n \circ \eta (x_i) \,  \neq \, \tilde h_n \circ tilde \eta (x_i)= \hat x_i^n.$$ 
But, $(t_n,s_n)$ is in the same left coset as $(\hat t_n,\hat s_n)$,
and both tuples $(t_n,s_n,x_1^n,\ldots,x_r^n)$ and
$(\hat t_n,\hat s_n,\hat x_1^n,\ldots,\hat x_r^n)$ are contained in the definable set $L(t,s,x_1,\ldots,x_r)$. Hence,
for every index $n$, and every $i$, $1 \leq i \leq r$, $x_i^n=\hat x_i^n$, and we get a contradiction.

\line{\hss$\qed$}

Since left cosets of cyclic subgroups are not real, right cosets of cyclic subgroups are not real as well  (though as basic sorts for the elimination
of imaginaries we do not need them).
Having Proved that cosets of cyclic groups are not reals, we further show that double cosets of cyclic groups
are also not real.

\vglue 1pc
\proclaim{Theorem 2.5} Let $F_k=<a_1,\ldots,a_k>$ be a non-abelian free group, let $m_1,m_2$ be  
positive integers and let:
$$Dcoset(y_1,x_1,z_1,y_2,x_2,z_2) \ = \ \{ \, (y_1,x_1,z_1,y_2,x_2,z_2) \ | \ y_1,z_1,y_2,z_2 \neq 1 \wedge 
[y_1,y_2]=[z_1,z_2]=1 \wedge $$
$$ \wedge \exists y,z \, [y,y_1]=[z,z_1]=1 \wedge y^{m_1}x_1z^{m_2}=x_2 \, \}$$
be the definable equivalence relation associated with double cosets of cyclic 
subgroups over $F_k$. Then $Dcoset$ is not real.
\endproclaim

\nfp A straightforward modification of the proof for left  cosets (theorem 2.3), proves that
double cosets are not reals.

\line{\hss$\qed$}

So far we proved that conjugation, left (and right) cosets of cyclic groups and double cosets of cyclic groups are not reals. V. Guirardel
pointed out some additional definable equivalence relations over free (and hyperbolic) groups, that can not be reduced to the 3 basic families.
We start with a refinement of the double coset equivalence relation, in case the left and right cyclic groups that define the double coset belong
to the same maximal cyclic subgroup.

\vglue 1pc
\proclaim{Theorem 2.6} Let $F_k=<a_1,\ldots,a_k>$ be a non-abelian free group, let $m_1,m_2$ be  
positive integers and let:
$$RDcoset(y_1,x_1,y_2,x_2) \ = \ \{ \, (y_1,x_1,y_2,x_2) \ | \ y_1,y_2 \neq 1 \wedge \ [y_1,y_2]=1 \ wedge \
\exists y \, [y,y_1]=1 \wedge y^{m_1}x_1y^{m_2}=x_2 \, \}$$
Then $RDcoset$ is a (definable) equivalence relation, and it is not real.
\endproclaim

\nfp $RDcoset$ is an equivalence relation by definition. It is not real, by the same argument that was used in the proof of theorem 2.3.

\line{\hss$\qed$}

Cosets, double cosets, and the refined double cosets (theorem 2.6) can be further generalized to what we call $generalized$ $double$ $cosets$.

\vglue 1pc
\proclaim{Definition 2.7} Let $F_k=<a_1,\ldots,a_k>$ be a non-abelian free group. Let $q,r,s$ be positive integers. For each $t$, $1 \leq t \leq s$,
we associate a positive integer (length) $\ell_t$, and a tuple of integers: $1 \leq b(1,t),\ldots,b(\ell_t,t) \leq q$.
For each $i$, $1 \leq i \leq r$, we associate a positive integer (rank), $d_i$.

Given all these tuples of integers, we set the generalized double coset (definable) equivalence relation to be: 
$$GDcoset(u_1,\ldots,u_q,y_1,\ldots,y_r,\hat u_1,\ldots,\hat u_q,\hat y_1,\ldots,\hat y_r) \ = \ 
 y_1,\ldots,y_r,\hat y_1,\ldots,\hat y_r \neq 1 \ \wedge $$
$$ \wedge \ [y_i, \hat y_i]=1, \  1 \leq i \leq r \ \wedge$$  
$$\wedge \ \forall \ 
z_{1,1},\ldots,z_{1,d_1},\ldots,z_{r,1},\ldots,z_{r,d_r} \ [z_{i,j},y_i]=1 \ 1\leq i \leq r, \ 1 \leq j \leq d_i $$ 
$$x_t=w_{1,t}(z_{1,1},\ldots,z_{r,d_r})u_{b(t,1)}
w_{2,t}(z_{1,1},\ldots,z_{r,d_r})u_{b(t,2)} \ldots
w_{\ell_t,t}(z_{1,1},\ldots,z_{r,d_r})$$ 
$$u_{b(t,\ell_t)}
w_{\ell_t+1,t}(z_{1,1},\ldots,z_{r,d_r}),\ 1 \leq t \leq s$$
$$\exists \ 
\hat z_{1,1},\ldots, \hat z_{r,d_r} \ [\hat z_{i,j},y_i]=1, \ 1\leq i \leq r, \ 1 \leq j \leq d_i $$
$$x_t=w_{1,t}(\hat z_{1,1},\ldots,\hat z_{r,d_r}) \hat u_{b(t,1)}
w_{2,t}(\hat z_{1,1},\ldots,\hat z_{r,d_r}) \hat u_{b(t,2)} \ldots
w_{\ell_t,t}(\hat z_{1,1},\ldots,\hat z_{r,d_r})$$ 
$$\hat  u_{b(t,\ell_t)}
w_{\ell_t+1,t}(\hat z_{1,1},\ldots,\hat z_{r,d_r}),\ 1 \leq t \leq s$$
\endproclaim 
 
\vglue 1pc
\proclaim{Theorem 2.8} Generalized double cosets are equivalence relations, and they are not real.  
\endproclaim

\nfp A generalized double coset relation $GDcoset$ is reflexive  by definition. If two tuples $u_1,\ldots,u_q,y_1,\ldots,y_r$ and
 $\hat u_1,\ldots,\hat u_q,\hat y_1,\ldots,\hat y_r$ are in the relation, then for elements $z_{i,j}$, that satisfy
$[z_{i,j},y_i]=1$, $1 \leq i \leq r, 1 \leq j \leq d_i$,
that are all much longer than the word length of the $u$'s, $\hat u$'s and the $y$'s,  the corresponding values of $x_1,\ldots,x_s$ satisfy the
two systems of equalities that are specified in definition 2.6 for some values of the elements $\hat z_{i,j}$. 
In particular, these values of the  elements $\hat z_{i,j}$ are also much longer than the lengths of the 
$u$'s, $\hat u$'s and the $y$'s.  

The possibility to present the corresponding values of the $x$'s as words in long $z_{i,j}$'s and long $\hat z_{i,j}$ imply that for any values of 
the elements $\hat z_{i,j}$, that satisfy $[\hat z_{i,j},y_i]=1$, the corresponding values of the elements $x_1,\ldots,x_s$ satisfy the
two systems of equalities that are specified in definition 2.6 for some values of the elements $z_{i,j}$. This proves the symmetry of the relation
$GDcoset$. 

Exactly the same argument proves the transitivity of the relation $GDcoset$, hence, $GDcoset$ is an equivalence relation. $GDcoset$ is not real
by the same argument that was used to prove theorems 2.1 and 2.3.

\line{\hss$\qed$}

Generalized double cosets can be generalized further, by combining them with conjugations. We call this generalized equivalence relation,
$generalized$ $conjugated$ $double$ $cosets$.

\vglue 1pc
\proclaim{Definition 2.9} Let $F_k=<a_1,\ldots,a_k>$ be a non-abelian free group. Let $p,r$ be positive integers. 
For each $i$, $1 \leq i \leq r$, we associate a positive integer (rank), $d_i$.

For each $m$, $1 \leq m \leq p$, there are two positive integers, $q_m,s_m$. For each $m$ and $t$, $1 \leq t \leq s_m$,
we associate a positive integer (length) $\ell_{m,t}$, and a tuple of integers: $1 \leq b(m,1,t),\ldots,b(m,\ell_{m,t},t) \leq q_m$.

Given all these tuples of integers, we set the generalized conjugated double coset (definable) equivalence relation to be:
$$GCDcoset(u^1_1,\ldots,u^1_{q_1},y^1_1,\ldots,y^1_r,\ldots, u^p_1,\ldots, u^p_{q_p}, y^p_1,\ldots, y^p_r,$$
$$\hat u^1_1,\ldots,\hat u^1_{q_1},\hat y^0_1,\ldots,\hat y^0_r,\ldots, \hat u^p_1,\ldots, \hat u^p_{q_p}, \hat y^p_1,\ldots, \hat y^p_r) \ = $$
$$= \  y^m_1,\ldots,y^m_r,\hat y^m_1,\ldots,\hat y^m_r \neq 1, \ 1 \leq m \leq p \ \wedge \
\exists \ g_1,\ldots,g_p, \hat g_1,\ldots,\hat g_m $$
$$g_1=1 \ [{g_m}^{-1}y^m_ig_m,y^1_i]=[{\hat g_m}^{-1}\hat y^m_i{\hat g_m},y^1_i]=1 \ 1 \leq m \leq p, \ 1 \leq i \leq r \ \wedge$$
$$ \wedge \ \forall \
z_{1,1},\ldots,z_{1,d_1},\ldots,z_{r,1},\ldots,z_{r,d_r} \ [z_{i,j},y^1_i]=1 \ 1\leq i \leq r, \ 1 \leq j \leq d_i $$

$$x_{m,t}=w_{m,1,t}(z_{1,1},\ldots,z_{r,d_r})g_m^{-1}u_{b(m,t,1)}g_m
w_{m,2,t}(z_{1,1},\ldots,z_{r,d_r})g_m^{-1}u_{b(m,t,2)}g_m \ldots$$
$$w_{m,\ell_{m,t}}(z_{1,1},\ldots,z_{r,d_r})g_m^{-1}u_{b(m,t,\ell_t)}g_m
w_{m,\ell_{m,t}+1,t}(z_{1,1},\ldots,z_{r,d_r}), \ 1 \leq m \leq p, \ 1 \leq t \leq s_m$$
$$\exists \
\hat z_{1,1},\ldots,\hat z_{1,d_1},\ldots,\hat z_{r,1},\ldots,\hat z_{r,d_r} \ [\hat z_{i,j},y_i]=1, \ 1\leq i \leq r, \ 1 \leq j \leq d_i $$ 
$$x_{m,t}=w_{m,1,t}(\hat z_{1,1},\ldots,\hat z_{r,d_r}){\hat g_m}^{-1} \hat u_{b(m,t,1)}\hat g_m
w_{m,2,t}(\hat z_{1,1},\ldots,\hat z_{r,d_r}){\hat g_m}^{-1} \hat u_{b(m,t,2)} \hat g_m \ldots$$
$$w_{m,\ell_{m,t}}(\hat z_{1,1},\ldots,\hat z_{r,d_r}) {\hat g_m}^{-1} \hat u_{b(m,t,\ell_t)} \hat g_m
w_{m,\ell_{m,t}+1,t}(\hat z_{1,1},\ldots,\hat z_{r,d_r}), \ 0 \leq m \leq p, \ 1 \leq t \leq s_m$$

And a similar formula in which $\hat g_1=1$ (and not just $g_1$).
\endproclaim 

\vglue 1pc
\proclaim{Theorem 2.10} Generalized conjugated double cosets are equivalence relations, and they are not real.  
\endproclaim

\nfp The argument is a straightforward generalization of the argument for generalized double cosets (theorem 2.8).

\line{\hss$\qed$}

,

\vglue 1.5pc
\centerline{\bf{\S3. Separation of Variables}}
\medskip
In section 3 of [Se9] we have introduced Duo limit groups, and associated a finite
collection of Duo limit groups with a given rigid or solid limit group. In the first section
of this paper
we have constructed the Diophantine envelope of a definable set (theorem 1.5), and then 
used it to construct the Duo envelope of a definable set (theorem 1.7). 

Recall that by its definition (see definition 1.1), a
Duo limit group $Duo$ admits an amalgamated product: 
$Duo=<d_1,p>*_{<d_0,e_1>} \, <d_0,e_1,e_2>*_{<d_0,e_2>} \, <d_2,q>$
where $<e_1>$ and $<e_2>$ are free abelian groups with pegs in $<d_0>$, i.e., free abelian groups that 
commute with non-trivial elements in $<d_0>$.
A specialization of the parameters $<d_0>$ of a Duo limit group gives us
a Duo family of it.

To analyze definable 
equivalence relations over a free (or a hyperbolic) group, we will need to further study
the parameters ($<d_0>$) that are associated with the Duo families that are associated with the Duo
limit groups that form the Duo envelope of a  definable equivalence relation. To do that we will need to get a
better "control" or understanding of the parameters
$<d_0>$ that are associated with a duo family and then with an equivalence class of a given
equivalence relation.

In this section we modify and further analyze   
the construction of the Duo envelopes that were presented in theorem 1.7, 
in the special case of a definable equivalence relation.
We carefully study the set of values of the parameters that are associated with the duo families that are
associated with each equivalence class.
This careful study, that uses what we call $uniformization$ limit groups that we associate with the
Duo envelope, enables one to associate  
a "bounded" set of values of certain subgroups of the parameters that are associated 
with the Duo families of the Duo envelope, for each
equivalence class of a definable 
equivalence relation (the bounded set of values of the  subgroups of parameters 
is modulo the basic imaginaries that were presented
in the previous section). 

The bound that we achieve on the number of specializations of the subgroups that we
look at, allows us to obtain what we view as 
"separation of variables". This means that with the original subgroups of parameters, $<p>$ and $<q>$, we associate
a bigger subgroup, for which there exists a graph of groups decomposition, where $<p>$ is contained in one vertex group,
$<q>$ is contained in a second vertex group, and the number of specializations of the edge groups (up to the
imaginaries that were presented in the previous section) is (uniformly) bounded for
each equivalence class of $E(p,q)$. 

In particular, the bounded set of specializations of the edge groups are a class function for the given 
equivalence relation. However, this class function is not guaranteed to separate between classes.    

The   
approach to the separation of variables with the uniform bound on the number of specializations of the edge groups for each
equivalence class,
combines the techniques that were used
in constructing the Duo envelope in theorem 1.7 (mainly the sieve procedure for quantifier elimination
that was presented in [Se6]), 
together with the techniques that were used to
construct formal (graded) Makanin-Razborov diagrams in sections 2-3 in [Se2], and the proof of the existence
of a global bound on the number of
rigid and strictly solid families of specializations of rigid and solid limit groups, that was presented in 
sections 1-2 in [Se3]. 

In the next section we show how to use the separation of variables that is obtained in
this section to finally analyze definable equivalence relations. This means how to use the separation of variables to improve the class function
that we get in this section to separate between classes.  

\smallskip
Let $F_k=<a_1,\ldots,a_k>$ be a non-abelian free group, 
and let $E(p,q)$ be a definable equivalence relation over $F_k$. 
Recall that with the definable equivalence relation, $E(p,q)$, being a definable set, one associates (using
the sieve procedure) finitely many (terminal)
rigid and solid  limit groups, $Term_1,\ldots,Term_s$ (see theorem 1.6). With each of the terminal limit group
 $Term_i$ there is an  associated set,
$B(Term_i)$, and that the definable set $E(p,q)$ is the set:
$$ E(p,q) \, = \, \cup_{i=1}^s \, B(Term_i)$$



By theorems 1.5 and 1.7, with the given definable  equivalence relation $E(p,q)$, being a
definable set, we can associate a Diophantine and a Duo envelopes. 
Let $G_1,\ldots,G_t$ be the graded Diophantine envelope of the given definable equivalence relation $E(p,q)$ (with respect to the parameters $<q,a>$),
and let $Duo_1,\ldots,Duo_r$, be those duo limit groups in the Duo envelope that are associated with the graded resolutions
$G_1,\ldots,G_t$. i.e., the duo limit groups that are associated with the the graded completions in the graded Diophantine envelope, and
not with the graded completions of smaller complexities that are constructed from them by increasing the set of parameters and applying the sieve procedure, 
in the  
iterative procedure that 
constructs the Duo envelope in the proof of theorem 1.7.

Let $Duo$ be one of the Duo limit groups,
$Duo_1,\ldots,Duo_r$. By definition (see definition 1.1), $Duo$ can be represented as an amalgamated product:  
$$Duo=<d_1,p>*_{<d_0,e_1>} \, <d_0,e_1,e_2>*_{<d_0,e_2>} \, <d_2,q>.$$ By construction, in $Duo$ there exists a subgroup that demonstrates that 
generic elements $<p,q>$ in $Duo$ are indeed in the equivalence relation $E(p,q)$. We denote this subgroup $<u,p,q,a>$. The subgroup contains elements
that are rigid or strictly solid of the limit groups $Term_i$ and terminal groups of some of their associated graded closures, together with elements
that are assumed to be primitive roots of edge groups in these terminal limit groups (see the construction of the duo envelope in theorem 1.5 and of the duo
envelope in theorem 1.7). 

The subgroup $<u,p,q,a>$,
being a subgroup of $Duo$, 
inherits a graph of groups decomposition from the presentation of $Duo$ as an amalgamated
product. We denote the subgroup $<u,p,q,a>$ of $Duo$ by $Ipr$.

$Ipr$, being a subgroup of $Duo$, inherits a graph of groups decomposition from its action on the
Bass-Serre tree that is associated with the amalgamated product:
$$Duo=<d_1,p>*_{<d_0,e_1>} \, <d_0,e_1,e_2>*_{<d_0,e_2>} \, <d_2,q>.$$
Suppose that in this graph of groups decomposition of the subgroup $Ipr$ there exists an edge with a
trivial edge group. In that case $Ipr$ admits a non-trivial free decomposition, $Ipr=A*B$. 
Since $<p> < <d_1,p>$ and $<q> < <d_2,q>$, so $<p>$ and $<q>$ are either trivial or can be embedded into
vertex groups in the graph of groups decomposition that is inherited by $Ipr$ and so is the coefficient group $F_k$. Therefore,  either:  
\roster
\item"{(i)}" $<p,q,a>$ is a subgroup of $A$.  

\item"{(ii)}" $<p>$ is a subgroup of $A$ and $<q>$ is a subgroup of $B$ and $F_k$ is a subgroup of either $A$ or $B$. 

\item"{(iii)}" $<q>$ is a subgroup of $A$ and $<p>$ can be conjugated into $A$ and $F_k$ is a subgroup of either $A$ or the conjugate of $A$ that
contains $<q>$.

\endroster
If (i) holds, then the restrictions of the specializations of $Duo$ to specializations  of $Ipr$ do not restrict to
rigid and strictly solid specializations of the terminal limit groups of $Term_i$ and the terminal groups of their associated graded closures,
and that contradicts the way the duo limit group $Duo$ was constructed in theorem 1.7.

To deal with case (ii), 
we present the following theorem, that
associates finitely many rigid and solid limit groups (with respect to the parameter subgroup $<p,q,a>$)
with the equivalence relation $E(p,q)$, so that each
pair, $(p_0,q_0) \in E(p,q)$, can be proved to be in $E(p,q)$ by a rigid or a strictly solid
homomorphism from one of these limit
groups into $F_k$. Furthermore,
there exist at most finitely many equivalence 
classes of the equivalence relation $E(p,q)$, such that pairs $(p_0,q_0) \in E(p,q)$ that do not belong to these 
exceptional classes, 
can be proved to be in $E(p,q)$  using rigid or strictly solid homomorphisms of one of the finitely many associated
rigid and solid limit groups
(with respect to $<p,q,a>$), and these rigid and strictly solid homomorphisms
do not factor through a free product (of limit groups) as in case (ii). 

To prove the theorem we apply once again the
sieve procedure [Se6], which was originally used for quantifier elimination, and was also the main tool
in the construction of the Diophantine and Duo envelopes in theorems 1.5 and 1.7.

\vglue 1pc
\proclaim{Theorem 3.1} Let $F_k=<a_1,\ldots,a_k>$ be a non-abelian free group, 
and let $E(p,q)$ be a (coefficient-free) definable equivalence relation over $F_k$. There exist finitely many rigid and solid limit groups
(with respect to $<p,q,a>$) that we denote: $Ipr_1,\ldots,Ipr_w$, so that:
\roster
\item"{(1)}" for every pair $(p_0,q_0) \in E(p,q)$ there exists a rigid or a strictly solid homomorphism
(with respect to $<p,q,a>$) $h:Ipr_i \to F_k$, for some
index $i$, that contains a valid proof statement that $(p_0,q_0) \in E(p,q)$. Note that valid proof statements were used to construct the duo envelope,
and it is essentially the restriction of  the homomorphism $h$ to the image of the group $Ipr=<u,p,q,a>$ in $IPr_i$.

\item"{(2)}" there exist at most finitely many equivalence classes of $E(p,q)$, so that for every
pair $(p_0,q_0) \in E(p,q)$ that does not belong to one of these finitely many classes, there exists a 
rigid or a strictly solid homomorphism, $h:Ipr_i \to F_k$, for some index $i$, that restricts to a valid proof statement that $(p_0,q_0) \in E(p,q)$, and
so that $h$, and every strictly solid homomorphism in the strictly solid family of $h$, 
does not factor through a homomorphism $\nu$ onto a non-trivial free product (of limit groups)  $A*B$, in which $\nu(p) \in A$,
$\nu(q) \in B$ and $\nu(F_k)$ is in either $A$ or $B$.  
\endroster
\endproclaim

\nfp 
%
%
We look at the collection of all the homomorphisms from a subgroup $Ipr=<u,p,q,a> \to F_k$,  
that verify that  pairs $(p_0,q_0) \in E(p,q)$. With this collection we can naturally associate
(by section 5 in [Se1]) a finite collection of limit groups that we denote: $V_1,\ldots,V_f$. With each of these limit groups
we can associate its graded Makanin-Razborov diagram with respect to the parameter subgroup, $<p,q,a>$. We further look at the
collection of rigid and strictly solid homomorphisms of rigid and solid limit groups in these diagrams, that verify that
pairs $<p_0,q_0> \in E(p,q)$ (note that if a pair $(p_0,q_0) \in E(p,q)$, then there exists such a rigid or a strictly solid
homomorphism). With this collection of rigid and strictly solid homomorphisms (with respect to $<p,q,a>$), we can associate finitely
many limit groups, that we denote, $L_1,\ldots,L_g$.

At this point we 
collect the subcollection of this collection of homomorphisms (i.e., rigid and strictly solid homomorphisms with respect
to $<p,q,a>$ that verify that pairs $(p_0,q_0) \in E(p,q)$),
that
factor through a non-trivial free product of the form $A*B$, where $A$ and $B$ are limit groups,
so that $<p> \, < A$, $<q> \, < B$ and $F_k$ is in either $A$ or $B$. By the standard methods
of section 5 in [Se1], with the subcollection of such rigid and strictly solid 
homomorphisms we can naturally associate a finite
collection of limit groups (graded with respect to $<p,q,a>$),  that we denote $M_1,\ldots,M_e$.

By successively
applying the 
shortening argument to the subcollection of homomorphisms that factor through a free product of
limit groups (by considering
the actions of the graded limit groups $M_1,\ldots,M_e$ on the Bass-Serre trees corresponding to the free
products $A*B$ through which they factor), we can replace this subcollection of homomorphisms with a new
subcollection, and the finite collection of limit groups, $M_1,\ldots,M_e$, with a new
finite subcollection, $GFD_1,\ldots,GFD_d$, 
for which each of the (graded) limit groups, $GFD_1,\ldots,GFD_d$, admits a non-trivial free decomposition $A_j*B_j$, where
$<p> < A_j$ and $<q> < B_j$ and $F_k$ is in either $A_j$ or $B_j$.

Let $GFD_j=A_j*B_j$, so that $<p> \, < \, A_j$ and $<q> \, < \, B_j$. With $A_j$ and $B_j$ viewed as (ungraded) limit
groups, we can naturally associate their taut Makanin-Razborov diagrams (see section 2 in [Se4] for the
construction and properties of the taut diagram). With a taut resolution of $A_j$ and a taut resolution of $B_j$,
we naturally associate their free product which is a resolution of $GFD_j=A_j*B_j$. Let $Res$ be such a resolution
of $GFD_j$. $Res$ is composed from two ungraded taut resolutions, $Res_1$ of $A_j$ and $Res_2$ of $B_j$.

Before we continue to the next step of the construction of the limit groups and resolutions that we'll use
in order to prove theorem 3.1, we prove the following fairly straightforward lemma on the finiteness of
equivalence
classes that contain generic points in one of the resolutions $Res_1$ or $Res_2$ from which the resolution $Res$ is composed.

\vglue 1pc
\proclaim{Lemma 3.2} Let $GFD$ be one of the graded limit groups, $GFD_1,\ldots,GFD_d$ constructed above,
and let $Res$ be one of its constructed taut resolutions. The resolution $Res$ is composed from two separate resolutions, $Res_1$ that contains the subgroup $<p>$,
$Res_2$ that contains the subgroups $<q>$ and $F_k$ is contained in either $Res_1$ or $Res_2$. 
Then there exist at most finitely many $exceptional$ equivalence 
classes of the definable equivalence relation $E(p,q)$, for which:
\roster
\item"{(1)}" for each of the finitely many exceptional equivalence classes there exist a test sequence 
of $Res_1$ and some value $q_0$ (that depends on the class)  of the variables $q$, that restrict to pairs  
$\{(p_n,q_0)\}_{n=1}^{\infty}$ that are in the exceptional equivalence class, and so that the
specializations in the test sequence
restrict to valid proof
statements that the pairs 
$\{(p_n,q_0)\}_{n=1}^{\infty}$ 
are in the (definable) equivalence relation $E(p,q)$.

\item"{(2)}" for each test sequence of $Res_1$, and for an arbitrary value $q_0$ of the variables $q$,
for which the restricted pairs
$\{(p_n,q_0)\}_{n=1}^{\infty}$ are in $E(p,q)$, and so that the
specializations in the sequence restrict to valid proof statements that the pairs
are in the equivalence relation $E(p,q)$, all the pairs in the sequence $\{(p_n,q_0)\}$ belong to one of the finitely
exceptional equivalence classes. 
\endroster
\endproclaim

\nfp With each of the finitely many limit groups, $GFD_1,\ldots,GFD_d$, we have associated finitely many
(taut Makanin-Razborov) resolutions. Let $Res$ be a resolution of one of these limit groups $GFD_j$. $GFD_j=A_j*B_j$, $<p> \, < \, A_j$,
$<q> \, < \, B_j$ and $F_k$ is either in $A_j$ or $B_j$.  $Res$
is built from two (ungraded) resolutions, $Res_1$ of $A_j$ and $Res_2$ of $B_j$.

The limit group $GFD_j$ contains an image of the group $Ipr=<u,p,q>$. Hence, every specialization of $GFD_j$ restricts to a specialization
of $Ipr$, that is a proof statement, not necessarily a valid proof statement, that the restriction of the specialization  of $Ipr$ to
specializations of of $p$ and $q$ are in the equivalence relation $E(p,q)$.

We look at all the sequences of specializations that factor through the resolution $Res$ that restrict to valid proof statements (i.e., the restriction of the
sequence to specializations of $Ipr$ are valid proof statements), for which the specializations of $Res_2$ are fixed,
and hence the specialization $q_0$ of $q$ is fixed, and the specializations of $Res_1$ form a test sequence. Note that
test sequences of $Res_1$ do not depend on 
the value of $q$. 

By the construction of formal limit groups in section 3 of [Se2], with the collection of these sequences of specializations that factor through $Res$
it is possible to associate a finite collection of formal graded limit groups, and each such formal graded limit group is a graded closure of the
completion of $Res_2$ (which is an ungraded resolution), where the terminal  limit group in the base of each of these graded closures is a rigid 
or a solid limit group w.r.t.\ the parameter subgroup $<q>$. 

Given such a graded closure of $Res_1$, that we denote, $GCl_i(Res_1)$, we look at all the test sequences of the graded closure for which the restrictions
of these test sequences to the image of $Ipr$ fail to be valid proof statements, precisely as we did in the construction of the Diophantine and Duo
envelopes in section 1. Using the construction of formal limit groups (section 3 in [Se2]), with the collection of these graded test sequences we associate
finitely many graded closures of (the graded closure) $GCl_i(Res_1)$. 



The graded closures that were associated with $GCl_i(Res_1)$ indicate what test sequences fail to restrict to valid proof statements. There
are only finitely many such graded closures, and each graded closures is associated with a collection of cosets of finitely generated abelian groups
that are dual to the  finitely generated abelian groups  that appear as vertex group along $GCl_i(Res_1)$, For every value $q_0$ of the variables $q$ for which
the boundedly many   fibers of $GCl(Res_1)$ that are associated with $q_0$ are not covered by the fibers of the graded closures of $GCl(Res_1)$ that are
associated with $q_0$, there is a test sequence of $GCl(Res_1)$ in one of the fibers that are associated with $q_0$ that restrict to valid
proof statements. Hence, a test sequence that restrict to values $\{(p_n,q_0)\}$ that are in $E(p,q)$.

Since such test sequences have to be belong to finitely many possible cosets, and test sequences in the same coset that restrict to valid
proof statements have to be in the same equivalence class of $E(p,q)$, there are finitely many exceptional classes of $E(p,q)$ that contain
such test sequences. By construction, every such coset contains a test sequence that restrict to valid proof statements, so (1) and (2) hold.

\line{\hss$\qed$}



Given lemma 3.2, In order to prove theorem 3.1, we still need to study  pairs $(p_0,q_0) \in E(p,q)$ that:
\roster
\item"{(1)}" do not belong to the finitely many equivalence classes that were excluded in lemma 3.1 and contain generic points in
one of the resolutions $Res_1$ that contain the subgroup $<p>$ in a resolution $Res$ of one of the limit groups, $GFD_j$.

\item"{(2)}" $(p_0,q_0)$ can be extended to specializations of $GFD_1,\ldots,GFD_d$, and these extended specializations restrict
to  valid proof statements that demonstrate that the  pair $(p_0,q_0) \in E(p,q)$. 
\endroster

To analyze these pairs   we need
to construct new limit groups with new (valid) proof statements that admit homomorphisms 
that do not factor through a free product $A*B$,
in which $A$ and $B$ are limit groups, $<p> < A$,  $<q> < B$ and $F_k$ is in either $A$ or $B$.  
To construct these new limit groups and new proof statements, we apply (once again) the sieve procedure [Se6], that was originally
presented as part of the quantifier elimination procedure.

Let $GFD_j$ be one of the limit groups $GFD_1,\ldots,GFD_d$, and let $Res$ be one of the resolutions
in its taut Makanin-Razborov diagram. $GFD_j=A_j*B_j$ and $Res$ is composed from two resolutions, $Res_1$ of $A_j$
and $Res_2$ of $B_j$. 


We view the completion of the resolution $Res$, $Comp(Res)$, as a graded completion of the resolution $Res_1$. i.e., we place $Comp(Res_2$ as the base of
the completion, $Comp(Res)$, where the structure of the completion is the structure of $Comp(Res_1)$. Note that $Comp(Res)=Comp(Res_1)|*Comp(Res_2)$.
Viewing $Comp(Res)$ as a graded completion of $Comp(Res_1)$, we associate with it finitely many graded closures that collects all the graded test sequences
of $Comp(Res)$ (viewed as a graded completion of $Comp(Res_1)$) that restrict to proof statements that fail to be valid proof statements (the construction of
these graded closures appear in 
definitions 1.25-1.30 in [Se5]). We continue with the graded closures that we call $Extra$ it collect those test sequences that restrict to 
proof statements for which there exists an extra rigid or strict solid family that does not appear in the proof statement. Note that like the completion,
$Comp(Res)$, that decomposes as a free product $Comp(Res_1)*Comp(Res_2)$, the Extra graded closures decompose as a similar free product as well.

We continue with all the (finitely many) Extra resolutions and Extra graded closures that were constructed from resolutions $Res$
of the limit groups $GFD_j$, $1 \leq j \leq d$. 
We collect all the specializations
that factor and are taut with respect to the Extra graded closures, 
for which the elements that are
supposed to be extra rigid or strictly solid specializations and are
specified by these specializations collapse. i.e., satisfy one of finitely many Diophantine conditions that force the extra
rigid and strictly solid specializations to be non-rigid or non strictly solid, or to be in the families that are specified in the 
associated specialization of the proof statement.

By our standard methods (section
5 in [Se1]), with the entire collection of specializations that factor through an Extra graded closure
and restrict to elements that are taut with respect to the (taut) resolution, $Res$,
 and for which the elements that are supposed to be
extra rigid or strictly solid specialization satisfy one of the finitely many
possible (collapse) Diophantine conditions, together with specializations of elements that demonstrate the
fulfillment of these Diophantine conditions, we can associate finitely many limit groups. We denote these 
limit groups, $Collape_1,\ldots,Collapse_f$, and call them $Collapse$ limit groups. 

Let $Collapse$ be one of the Collapse limit groups, $Collapse_1,\ldots,Collapse_f$. 
With $Collapse$ we associate its graded Makanin-Razborov diagram with respect to the
parameter subgroup $<p,q,a>$. We continue with all the rigid and strictly solid homomorphisms of 
rigid and solid limit groups in this Makanin-Razborov diagram. We look at all the rigid and 
strictly solid specializations of rigid and solid limit groups in this diagram, so that their restrictions
to  
specializations of the corresponding limit group $GFD$ are valid proofs that $(p,q) \in E(p,q)$, and
for which the specializations  factor through a non-trivial free product $A*B$, so that $A$ and $B$ are limit groups, 
$<p> < A$, $<q> < B$ and $F_k$ is in one of them. 
With this collection of (rigid and strictly solid)
homomorphisms we can associate a finite collection of limit groups (by the standard techniques that are
presented in section 5 of [Se1]), that we denote, $R_1,\ldots,R_m$.

Given a limit group $R_j$, we can associate with it a graded Makanin-Razborov diagram in which every
 graded resolution (with respect to the parameter
subgroup $<p,q>$) terminates in a rigid or in a solid limit group, and this terminal rigid or solid limit group
admits a non-trivial free decomposition $A*B$, in which $A$ and $B$ are limit groups, $<p> < A$, $<q> < B$ and $F_k$ is in either $A$ or $B$. 

At this point we combine the  graded Makanin-Razborov diagram of $Collapse$ with the graded Makanin-Razborov
diagrams of each of the limit groups $R_1,\ldots,R_m$. Each of the resolutions in the graded Makanin-Razborov
diagram of $Collapse$ terminates in a rigid or a solid limit group. We replace this graded resolution with
finitely many resolutions. First we replace its terminating rigid or solid limit group by each of
the quotients that
are associated with it from the set $R_1,\ldots,R_m$. We continue each of the obtained resolutions (after performing
the replacements) with the graded Makanin-Razborov diagram of the corresponding limit group $R_j$. By construction,
each of the
constructed  resolutions starts with $Collapse$ and terminates with a limit group that admits a free
product in which $<p>$ is a subgroup of one factor and $<q>$ is a subgroup of the second factor.

Given the obtained (graded) diagram of the limit group $Collapse$, we replace it with a strict (graded) diagram,
according to the finite iterative procedure that is presented in proposition 1.10 in [Se2]. Note that each resolution
in the strict diagram starts with a quotient of $Collapse$ and terminates with a limit group that admits a free 
product in which $<p>$ is contained in one factor and $<q>$ is contained in the second factor.

Let $CRes_1$ be a (graded) resolution in the strict diagram that is associated with $Collapse$.
Note that given a  homomorphism
$h:Collapse \to F_k$ that factors through $CRes_1$,
and $h$ restricts to a specialization of $GFD$ which
is a valid proof that $(p,q) \in E(p,q)$, the rigid vertex groups in the graded abelian decompositions that are
associated with  the rigid and strictly solid specializations in $GFD$, remain elliptic through the entire
combined resolution $CRes_1$ (i.e., remain elliptic in both resolutions from which $CRes_1$ is constructed).

Suppose that a rigid vertex group or an edge group
in the abelian decomposition that is associated with the specialization of the extra graded closure,
from which $Collapse$ was constructed, does not remain elliptic through $CRes_1$. In that case a non-trivial relation
that does not hold in the Extra graded closure holds in the terminal limit group of $CRes_1$.

Suppose that all the  rigid vertex groups and edge groups 
in the abelian decomposition that is associated with the specialization of the extra rigid or strictly
solid specialization in $Collapse$ remain elliptic through $CRes_1$.
 Then all these rigid vertex groups can be conjugated into the factors of the terminal limit group 
of $CRes_1$, and hence, 
for the relevant specializations that
factor through the Extra graded closure
the Diophantine condition that we imposed on the extra rigid or strictly
solid specialization in the Extra graded closure  is in fact rephrased as a Diophantine condition on the terminal limit group
of $CRes_1$, and indeed it is part of the terminal limit group of $CRes_1$ which is a free product $A_1*B_1$ in which 
$<p> < A_1$, $<q> < B_1$ and $F_k$ is in either $A_1$ or $B_1$.


Since both the Extra graded closure, 
and the Diophantine conditions that are imposed on them, admit
 a free product in which $<p>$ is contained in one factor, and $<q>$ is contained in the second factor, to
analyze the set of specializations that factor through  an extra graded closure 
are taut with respect to an original resolution, $Res$, of one of the limit groups, $GFD_1,\ldots,GFD_d$,
and so that these (extra rigid and strictly solid) 
specializations extend to specializations that satisfy the Diophantine conditions that are imposed on
them, we can use the analysis of such resolutions that was presented in the sieve procedure [Se6]. 

\smallskip 
We continue iteratively as in the sieve procedure. At each step we start with the collection of
Extra graded closures (w.r.t. the subgroup $<q>$) that are constructed from the resolutions 
that were constructed in the
previous step. 
We look at the collection of specializations that factor through
and are taut with respect to these resolutions, that satisfy one of finitely many (collapse) Diophantine conditions, 
and with this collection we associate (using section 5 in [Se1]) finitely many limit groups that we denote $Collapse$.

With each of the obtained limit groups $Collapse$ we associate its graded Makanin-Razborov diagram
(with respect to the subgroup $<p,q,a>$). Given each of the rigid and solid limit groups in this diagram,
we collect all the rigid or strictly solid homomorphisms of it that factor through a free product of limit
groups in which $<p>$
is contained in one factor, and $<q>$ is contained in the second factor and $F_k$ is either in the first or the second factor. 

We collect all these rigid and strictly
solid homomorphisms in finitely many limit groups, and with each such limit group we associate a graded
Makanin-Razborov diagram (with respect to $<p,q,a>$) that terminates in rigid and solid limit groups that admit
 free products in which $<p>$ is contained in one factor, $<q>$ is contained in the second factor and $F_k$ is either in the first or the second factor.

We further combine the graded Makanin-Razborov diagrams of each of the limit groups $Collapse$ with the graded
Makanin-Razborov diagrams of the limit groups that are associated with the rigid and solid limit groups in these
diagrams. Given a combined diagram we replace it by a strict diagram using the iterative procedure that appears in
proposition 1.10 in [Se2]. Each resolution in the strict diagram that is associated with a limit group
$Collapse$, starts with a quotient of $Collapse$ and terminates in a rigid or a strictly solid homomorphism that admits a free product in which $<p>$ 
is contained in one factor,  $<q>$ is contained in the second factor, and $F_k$ is either in the first or the second factor.
 

We continue by associating (taut) resolutions with these limit groups according to the construction that
is used in the sieve procedure [Se6]. 
By lemma 3.2 there are at most finitely many equivalence classes of the equivalence relation $E(p,q)$ for which a test
sequence of one of the constructed (taut) resolutions   restricts to valid proofs that the corresponding
couples $\{(p_n,q_0\}$ are in the set $E(p,q)$. 

By the termination of the sieve procedure (theorem 22 in [Se6]) the procedure terminates after finitely many steps.
By applying 
lemma 3.2 in the various steps of the iterative procedure, 
there exist at most finitely many equivalence classes of $E(p,q)$ so that if $(p,q) \in E(p,q)$,
and $(p,q)$ does not belong to one of these classes, then there exists a rigid or a strictly solid 
homomorphism from one of the limit groups that were constructed along the procedure,
$Ipr_1,\ldots,Ipr_w$, that restricts to a valid proof that the couple $(p,q)$ is in the set $E(p,q)$, and furthermore,
this rigid homomorphism and every strictly solid homomorphism which is in the same strictly solid family of the
strictly solid homomorphism does not factor through a free product of limit groups
in which $<p>$ is contained in one factor, 
$<q>$ is contained in the second factor, and $F_k$ is contained in either the first or the second factor. Hence, theorem 3.1
follows.

\line{\hss$\qed$}

\smallskip
Theorem 3.1 associates with the given definable equivalence relation, $E(p,q)$, finitely many rigid and solid limit
 groups, $Ipr_1,\ldots,Ipr_w$, so that apart from finitely many equivalence classes, for each couple,
$(p,q) \in E(p,q)$, there exists a rigid or a strictly solid family of 
homomorphisms from at least one of the limit groups, $Ipr_1,\ldots,Ipr_w$,
to the coefficient group $F_k$, so that  the rigid homomorphisms or the strictly solid homomorphisms  
from the given strictly solid family
do not factor through a non-trivial free product of limit groups, $A*B$, in which $<p> < A$, $<q> < B$ and $F_k$ is in either $A$ or $B$, 
and each of these homomorphisms
restricts to
a valid proof that $(p,q) \in E(p,q)$.

The rigid and solid limit groups $Ipr_1,\ldots,Ipr_w$ and their rigid and strictly solid families of
homomorphisms that 
do not factor through  graded free products and restrict to valid proofs, 
are the starting point to our approach to 
associating (definable) parameters with the equivalence classes of the definable equivalence relation $E(p,q)$.

Recall that by theorems 1.5 and 1.7, with the given definable  equivalence relation $E(p,q)$, being a
definable set, we can associate a Diophantine and a Duo envelopes. 
We denoted by $G_1,\ldots,G_t$ the Diophantine envelope of the given definable equivalence relation $E(p,q)$,
and by $Duo_1,\ldots,Duo_r$, its Duo envelope. 

We continue by modifying the construction of the Duo envelope that is presented in theorem 1.7, and use the
collection of homomorphisms from the limit groups, $Ipr_1,\ldots,Ipr_w$, that do not factor through a free
product in which (the image of) $<p>$ is contained in one factor, and (the image of) $<q>$ is contained in
the second factor, and $F_k$ is contained in one of the factors  (see the proof of theorem 1.7). 

\medskip
Let $G_1,\ldots,G_t$ be the Diophantine envelope of the definable equivalence relation,  $E(p,q)$ 
(see theorem 1.5).
Since we are analyzing definable equivalence relations and not general
definable sets, we are interested only in duo limit groups that are associated with the (graded) resolutions $G_1,\ldots,G_t$ that form
the Diophantine envelope of the definable equivalence relation $E(p,q)$. We do not need the procedure that gradually increases the set of 
parameters: $q_1,q_2,\ldots$, and constructs the associated resolutions,
that we used in theorem 1.7 as part of the construction of the Duo envelope of a general definable set.  

We start with the graded completions $G_1,\ldots,G_t$
in parallel.
With each graded completion $G_j$, $1 \leq j \leq t$,  we first associate a finite collection  of duo
limit groups.

To construct these duo limit groups, we look at the entire collection
of graded test sequences that factor through the given graded completion, $G_j$, for which the 
(restricted) sequence of specializations $\{(p_n,q_n)\}$ can be extended to specializations of one of the 
limit groups, $Ipr_1(f,p,q,a),\ldots,Ipr_w(f,p,q,a)$, so that these specializations of the subgroups $Ipr_s$ are:
\roster
\item"{(1)}" rigid or 
almost shortest in their strictly solid family (see definition 2.8 in [Se3] for an almost shortest 
specialization).
\item"{(2)}" the images of the subgroups $Ipr_s$ do not factor through a free product in which the subgroup
$<p>$ is in one factor, and the subgroup $<q>$ s in the  second factor and $F_k$ is in one of the factors.

\item"{(3)}" in each such test sequence of $G_j$ the specializations of the subgroup $<q>$ are fixed.
\endroster

With this entire collection of graded test sequences, and their extensions to specializations of
the limit groups $Ipr_1,\ldots,Ipr_w$,
we associate finitely many  graded Makanin-Razborov diagrams, precisely as we did in constructing
the formal graded Makanin-Razborov diagram in section 3 of [Se2]. As in the
formal Makanin-Razborov diagram, each resolution in the diagrams we construct terminates with
a (graded) closure of the given graded completion, $G_j$, we have started with, 
amalgamated with another group along its base
(which is the terminal rigid or solid limit group of the graded completion $G_j$), and the abelian vertex groups
that commute with non-trivial elements in the base).
  
We continue as in the proof of theorem 1.7. 
By construction, a completion of a resolution in one of the constructed graded diagrams is
a duo limit group.
We take the completions of the resolutions that appear in the 
finitely many diagrams that are associated with the graded completion $G_j$,
to be the  preliminary (finite) collection of duo limit groups that are associated with $G_j$.
We proceed by applying the sieve procedure to the constructed duo limit groups, precisely as we
did in the (second part of the) construction of the duo envelope in proving theorem 1.7.

Finally, we set the Duo envelope of the definable equivalence relation $E(p,q)$, that we denote,
$TDuo_1,\ldots,TDuo_m$, to be those duo limit groups
that are associated with the Diophantine envelope, $G_1,\ldots,G_t$, for which there exists a duo family
having a test sequence, so that all the specializations in the test sequence restrict to
elements $(p,q)$
in $E(p,q)$, and for which the associated specializations of the subgroup, $Ipr_s$, testify that indeed
the elements $(p,q)$ are in $E(p,q)$
(i.e., in particular, a "generic point" in $TDuo_i$ restricts to elements in $E(p,q)$, and the corresponding restrictions
to (the image of) $Ipr_s$
are form a proof of that).

Note that by construction, the collection of duo limit groups that we constructed, $Tduo_1,\ldots,Tduo_m$,
satisfy   properties (1) and (2) in theorem 1.7. They satisfy the universality property (3) in theorem 1.7, only with respect to
duo limit groups that are associated with the graded resolutions that form the Duo envelope of $E(p,q)$, i.e., $G_1,\ldots,G_t$.

\vglue 1pc
\proclaim{Proposition 3.3} The Duo limit groups, 
$Tduo_1,\ldots,Tduo_m$, 
and the rigid and solid limit groups, $Ipr_1,\ldots,Ipr_w$,
have the following properties:
\roster
\item"{(1)}" With each of the Duo limit groups $Tduo_i$ there is an associated homomorphism from 
an associated graded completion, $G_{j(i)}$,  which is one of the graded completions
that form the Diophantine envelope of $E(p,q)$.
Furthermore, the graded completion $G_{j(i)}$ has the same structure as one of the
two graded completions that are associated with $Tduo_i$. In fact, the corresponding graded completion
in $Tduo_i$ is a graded closure of $G_{j(i)}$, and $G_{j(i)}$ is mapped into this closure preserving the
level structure.

\item"{(2)}" With each Duo limit group $Tduo_i$, there is an associated homomorphism from
one of the limit groups, $Ipr_1,\ldots,Ipr_w$, into $Tduo_i$ that does not factor through a non-trivial free
product in which $<p>$ is contained in one factor,  $<q>$ is contained in the second factor and $F_k$ is contained in one of the two factors.
We denote the image of this limit group in $Tduo_i$, $<f,p,q,a>$. 

\item"{(3)}"  $Tduo_i$, being a Duo limit group, admits the amalgamated product:
$Tduo_i=<d_1^i,p>*_{<d_0^i,e_1^i>} \, <d_0^i,e_1^i,e_2^i>*_{<d_0^i,e_2^i>} \, <d_2^i,q>$. 
If both subgroups $<p>$ and $<q>$ are non-trivial in $<f,p,q,a>$, 
then the subgroup $<f,p,q>$ intersects non-trivially some
conjugates of the distinguished vertex group in $Tduo_i$, $<d_0^i,e_1^i,e_2^i>$. 
\endroster
\endproclaim

\nfp Parts (1) and (2) follow from the construction of the duo limit group, $Tduo_i$, that starts with one
of the graded completions, $G_{j(i)}$, and continues by collecting all the test sequences of $G_{j(i)}$,
for which their restrictions to the subgroup  $<p>$  
can be extended to specializations of one of the limit groups, $Ipr_1,\ldots,Ipr_w$. 

Let $<f,p,q,a>$ be the image of one of the rigid or solid limit groups, $Ipr_1,\ldots,Ipr_w$, in $Tduo_i$. 
The subgroup, $<f,p,q,a>$, inherits a graph of groups decomposition from the amalgamation of the ambient group
$Tduo_i$:  
$$Tduo_i=<d_1^i,p>*_{<d_0^i,e_1^i>} \, <d_0^i,e_1^i,e_2^i>*_{<d_0^i,e_2^i>} \, <d_2^i,q>$$ 
If both subgroups $<p>$ and $<q>$ are non-trivial in $<f,p,q,a>$, and $<f,p,q>$ intersects trivially all the 
conjugates of the vertex group, $<d_0^i,e_1^i,e_2^i>$, the graph of groups decomposition that is inherited by
$<f,p,q>$ collapses into a non-trivial free product of $<f,p,q,a>$ in which $<p>$ is contained in one factor, 
$<q>$ is contained in a second factor, and $F_k$ is in one of the factors. However, the duo limit group $Tduo_i$ was constructed from 
specializations of one of the limit groups, $Ipr_1,\ldots,Ipr_w$, that do not factor through a free product
of limit groups in which $<p>$ is contained in one factor and $<q>$ is contained in a second factor, a 
contradiction. Hence, in case both subgroups $<p>$ and $<q>$ are non-trivial, $<f,p,q,a>$ intersects non-trivially
some conjugate of the vertex group,
$<d_0^i,e_1^i,e_2^i>$, and part (3) follows.

\line{\hss$\qed$}

Part (3) of proposition 3.3 uses the fact that the homomorphisms from the rigid and solid limit groups,
$Ipr_1,\ldots,Ipr_w$, that we use to verify that (generic) couples
$(p,q)$ in the Duo limit groups, $TDuo_1,\ldots,TDuo_m$, are in the given definable equivalence relation
  $E(p,q)$, do not factor
through a non-trivial free product in which $<p>$ is contained in one factor,  $<q>$ is contained in the second
factor and $F_k$ is contained in one of the factors, to deduce that  
$<f,p,q,a>$ intersects non-trivially some conjugates of $<d_0^i,e_1^i,e_2^i>$. The analysis of the specializations of these intersections is
a key in our approach to associating parameters with the families of equivalence classes of the equivalence relation,
$E(p,q)$. 

For presentation purposes, we first continue by assuming that the Duo limit groups, $Tduo_1,\ldots,Tduo_m$, terminate
in rigid limit groups, i.e., that the abelian decomposition that is associated with the limit group $<d_0^i>$ is
the trivial (graded) decomposition. We further assume that the graded closures that are associated with the duo limit
groups, $Tduo_1,\ldots,Tduo_m$, do not contain abelian vertex groups in any of their levels. Hence, in particular
the subgroup $<d_0^i,e_1^i,e_2^i>$ is simply $<d_0^i>$. 
In the sequel, we will further assume that the 
Duo (and uniformization) limit groups that are constructed from these Duo limit groups and are associated with
them, terminate in rigid limit groups, and the graded closures that are associated with them do not contain 
abelian vertex groups in any of their levels   as well. These assumptions will allow us to present our approach
to separation of variables, and to associating parameters with the equivalence classes of $E(p,q)$, 
while omitting some technicalities.
Later on we omit these assumptions, and generalize our approach to work in the presence of both rigid and solid 
terminal limit groups, and when abelian groups do appear as vertex groups in the abelian decompositions that are
associated with the graded closures that are associated with the constructed duo limit groups.

\vglue 1pc
\proclaim{Proposition 3.4} Let $Tduo_i$ be one of the Duo limit groups, 
$Tduo_1,\ldots,Tduo_m$, and suppose that $Tduo_i$ terminates in a rigid limit group, i.e.,
that the abelian decomposition that is associated with $<d_0^i>$ is trivial. Suppose further that the two graded
completions that are associated with $Tduo_i$ contain no abelian vertex groups in any of their levels.

$Tduo_i$ being a duo limit group, admits a presentation as an amalgamated product:
$Tduo_i=<d_1^i,p>*_{<d_0^i>} \, <d_2^i,q>$. 
Suppose that the subgroups $<p>$ and $<q>$ are both non-trivial in $Tduo_i$. 
By proposition 3.3, the subgroup  $<f,p,q>$, which is the image of one of the rigid or solid limit groups,
$Ipr_1,\ldots,Ipr_w$, in $Tduo_i$, intersects
non-trivially some conjugates  of the edge group $<d_0^i>$. 
Let $H_i^1,\ldots,H_i^e$ be these (conjugacy classes of) intersection subgroups.

 Let 
$G_{j(i)}$ be the graded completion from the Diophantine envelope, 
$G_1,\ldots,G_t$, that is mapped into $Tduo_i$.  
$G_{j(i)}$, being a graded completion, has a distinguished vertex group, which is a subgroup of its
base subgroup.

 Then there exists a global integer $b_i>0$, so that for any specialization
of the distinguished vertex group of $G_{j(i)}$, 
there are at most $b_i$ rigid  specializations of $<d_0^i>$ that extend to 
generic (i.e., restrictions of (duo) test sequences) rigid and strictly solid specializations
of $<f,p,q,a>$, that form valid proofs that the pairs  $(p,q)$ are in  $E(p,q)$. 
In particular, these specializations of $<d_0^i>$ restrict to
 at most $b_i$   conjugacy classes of specializations of the subgroups,    
$H_i^1,\ldots,H_i^e$.
\endproclaim

\nfp Since the subgroup, $<d_0^i>$, of the duo limit group, $Tduo_i$, is assumed to be rigid, the proposition 
follows from the existence of a uniform bound on the number of rigid specializations of a rigid 
limit group with a fixed value of the defining parameters (i.e., a bound that does not depend on the
specific value of the defining parameters) that was proved in theorem 2.5 in [Se3].

\line{\hss$\qed$}

Proposition 3.4 proves that for a given specialization of the distinguished vertex group in $G_{j(i)}$,
there are at most boundedly many conjugacy classes of
specializations
of the corresponding subgroups, $H_i^1,\ldots,H_i^e$,  that may be associated with it.
However, given an equivalence class of $E(p,q)$ we can't, in general, associate with it only finitely
many conjugacy classes of specializations of the subgroups $H_i^1,\ldots,H_i^e$. 
Hence, to obtain only boundedly many specializations or conjugacy classes of specializations
of some "preferred" groups of parameters 
that are associated with each equivalence class (and not only with
a specialization of the distinguished vertex group in $G_{j(i)}$),
we need to construct $uniformization$  $limit$ $groups$.

\smallskip
Let $Tduo$ be one of the Duo limit groups, $Tduo_1,\ldots,Tduo_m$. We assume as in proposition 3.4,
 that $Tduo$ terminates in a rigid limit group 
(i.e., the subgroup $<d_0>$ admits a trivial graded JSJ decomposition), and that the two graded completions
that are associated with $Tduo$ contain no abelian vertex group in any of their levels. 
$Tduo$, being a Duo limit group (with no abelian vertex groups that appear along the levels of its
two associated graded completions), admits the amalgamated product:
$Tduo=<d_1,p>*_{<d_0>} \, <d_2,q>$.  

By part (1) of proposition 3.3, with $Tduo$ there is an associated homomorphism from 
an associated graded completion, $G_j$,  which is one of the graded completions in the Diophantine envelope 
of $E(p,q)$.
By part (2) of proposition 3.3, with $Tduo$ there is also an associated homomorphism from one of the rigid
or solid limit groups, $Ipr_1,\ldots,Ipr_w$, into $Tduo$.

We denote the image of this homomorphism in $Tduo$, $<f,p,q,a>$. 
Note that  by proposition 3.3 if both subgroups $<p>$ and $<q>$ in $<f,p,q,a>$ are non-trivial, then
the intersection of $<f,p,q,a>$ with some conjugates of $<d_0>$ is non-trivial. 
We denote by 
$H^1,\ldots,H^e$
these intersection subgroups.

Suppose
that there exists an equivalence class of $E(p,q)$, for which there exists an infinite sequence of conjugacy
classes of specializations
of $H^1,\ldots,H^e$ that can be extended to couples of test sequences of the two graded completions 
that are associated with
$Tduo$, so that restrictions of generic  elements in these test sequences, $(f_n,p_n,q_n,a)$, prove that the couples
$(p_n,q_n) \in E(p,q)$, these test sequences restrict to valid proofs that the couples $(p_n,q_0(n))$ and
$(q_0(n),q_n)$ belong to $E(p,q)$ (recall that $q_0(n)$ is the restriction of the specializations $d_0(n)$ to
the elements $q_0$), and furthermore these test sequences restrict to sequences of 
distinct couples, $\{(p_n,q_n)\}$.

Note that it is  a corollary of the quantifier elimination procedure (and the uniform bounds on rigid and
strictly solid specializations of rigid and solid limit groups obtained in [Se3]), that there is a global
bound on the size of all the finite equivalence classes of a definable equivalence relation.
To analyze the sets of specializations of the subgroups $H^1,\ldots,H^e$ that 
are associated with the same infinite equivalence
classes of $E(p,q)$, we construct finitely many limit groups (that are all associated with $Tduo$),
that we call $uniformization$ $limit$ $groups$. To construct these limit groups 
we look at the collection of all the sequences:
 $$\{(d_1(n),p_n, d_0(n), d_2(n),q_n,\hat f, \hat d_0(n))\}$$ for which:
\roster
\item"{(1)}" $\{(d_1(n),p_n,d_0(n))\}$ is a test sequence of the first graded completion  that 
is associated with $TDuo$, 
and $\{(d_2(n),q_n,d_0(n))\}$ is a test sequence of the second graded completion  that 
is associated with $TDuo$. These sequences restrict to proofs that the couples $(p_n,q_0)$ and $(q_0,q_n)$ are in
$E(p,q)$, and the couples, $\{(p_n,q_n)\}$, in such a sequence are distinct.

\item"{(2)}" the sequence $\{(d_1(n),p_n,d_0(n),d_2(n),q_n)\}$ restricts to a sequence of specializations: 
$\{(f_n,p_n,q_n)\}$, that  are rigid or strictly solid specializations of the rigid or solid
limit group $Ipr$ that is associated with $Tduo$ (see part (2) of proposition 3.3).  
Furthermore, the sequence of specializations,
$\{(f_n,p_n,q_n,a)\}$, restricts to proofs that the sequence $\{(p_n,q_n)\}$ are in the same equivalence class of the equivalence  relation $E(p,q)$.
As in constructing the Duo limit groups, $Tduo$, we  assume that the couples $(p_n,q_n)$
do not belong to the finitely many equivalence classes that are specified in theorem  3.1.

\item"{(3)}" the  elements  $\{(\hat f, d_0(n),\hat d_0(n))\}$ restrict to rigid or  strictly solid
specializations of one of the rigid and solid limit groups $Ipr_1,\ldots,Ipr_w$, that prove that the specializations $<d_0(n)>$ 
and the specialization $\hat d_0(n)$
belong to the same equivalence class in $E(p,q)$.
\endroster

Using the techniques of sections 2 and 3 in [Se2] for the construction of graded limit groups, 
we can associate with the above collection of sequences (for the entire collection of
Duo limit groups  
$Tduo_1,\ldots,Tduo_m$) a finite collection of
Duo limit groups, $Dduo_1,\ldots,Dduo_u$. 
By construction, the subgroup $<\hat f,d_0,\hat d_0>$ form the distinguished vertex 
groups of the constructed Duo limit groups, and the two graded completions that are associated with each such Duo
limit group have the same structure as those of the Duo limit group $Tduo$ from which they were constructed.

We continue with each of the distinguished vertex groups $<\hat f,d_0,\hat d_0>$ of the Duo limit groups,
$Dduo_1,\ldots,Dduo_u$. We view each of the vertex groups, $<\hat f,d_0,\hat d_0>$, as graded limit groups
with respect to the parameter subgroups $<\hat d_0>$, and associate with them their graded taut Makanin-Razborov 
diagrams. With each resolution in these diagrams we naturally associate its graded completion (see
definition 1.12 in [Se2] for the completion of a well-structured resolution). 

Given each graded completion that is associated with a limit group $<\hat f,d_0,\hat d_0>$, 
we construct a new limit group that starts with the graded completion of a resolution of a subgroup
$<\hat f,d_0,\hat d_0>$ (so that the terminal limit group of this graded completion contains the
subgroup $<\hat d_0>$), and on top of this completion we amalgam the two 
graded completions that were associated with
the associated subgroup $Dduo$. i.e., we get a Duo limit group that has the same structure as the 
associated Duo limit group $Dduo$, but the distinguished vertex in $Dduo$ is replaced with a
graded completion that terminates in a rigid or a solid limit group that contains $<\hat d_0>$ (which is the
parameter subgroup of this terminal rigid or solid limit group). We denote the limit groups that
are constructed in this way from the Duo limit groups, $Dduo_1,\ldots,Dduo_u$, $Cduo_1,\ldots,Cduo_v$.

The distinguished vertex group in the Duo limit groups $Dduo_1,\ldots,Dduo_u$, that was the limit group 
$<d_0,\hat f, \hat d_0>$, was replaced by the completions of the  resolutions in the taut graded 
Makanin-Razborov diagrams of these groups (with respect to $<\hat d_0>$) to obtain the Duo limit groups
$Cduo_1,\ldots,Cduo_v$. Each of the obtained Duo limit groups, that we denote $Cduo$, has the structure of
a completion (of a resolution of $<d_0,\hat f,\hat d_0>$), and on top of this completion we amalgam two
additional completions, which are the two completions that are associated with the Duo limit group, $Tduo$,
from which it was constructed.

Since the Duo limit group, $Cduo$, is constructed from 3 completions, we can naturally associate generic
points with it, i.e., test sequences that are composed from  test sequence of the completion of the
resolution of the limit group, $<d_0,\hat f, \hat d_0>$, that is extended to be a test sequence of the
two completions that are amalgamated with that completion and these have the structure of the two completions
that are associated with the Duo limit group, $Tduo$, from which the limit group $Cduo$ was constructed.

By construction, there exist generic points of the Duo limit group $Tduo$, and its associated Duo 
limit group, $Dduo$, i.e., (double) test sequences of the two completions that are associated with 
each of them, that restrict to proofs that the couples $(p_n,q_0(n))$ and $(q_0(n),q(n))$ 
are in $E(p,q)$, and restrict to
 specializations $\{(f_n,p_n,q_n,a)\}$ of the associated limit group
$Ipr$, that restrict to proofs that the couples $\{(p_n,q_n)\}$ are in the given definable set
$E(p,q)$. 

However, there is no guarantee that there exists a generic point of the limit group
$Cduo$, that is constructed from $Dduo$, with these properties, i.e., that there exists a (triple) test sequence,
that is composed from test sequences of the 3 completions that form the Duo limit group $Cduo$, so
that this (triple) test sequence restricts to proofs that the couples $(p_n,q_0(n))$ and $(q_0(n),q_n)$
are in $E(p,q)$ and to
specializations, $\{(f_n,p_n,q_n)\}$ and $\{(\hat f_n,d_0(n),\hat d_0(n))\}$,
that restrict to proofs that  the couples $\{(p_n,q_n)\}$ are in $E(p,q)$, and that for each $n$ the couple
$(d_0(n),\hat d_0(n))$ are in the same equivalence class of $E(p,q)$.

Therefore, we start with the Duo limit group $Cduo$, and apply the sieve procedure to it, 
in the same way that we constructed the
Diophantine and Duo envelopes in theorems 1.5 and 1.7. First, we look at all the (triple) test sequences of $Cduo$
for which specializations of subgroups of $<f,p,q,a>$ and $<\hat f,q_0,\hat q_0,a>$ and specializations of the
the subgroup that is supposed to
demonstrate that the specializations of $(p,q_0)$ and $(q_0,q)$ 
are in $E(p,q)$ that were supposed to be rigid or strictly
solid do not have this property. With this collection of test sequences we associate finitely many graded closures (following 
definitions 1.25-1.30 in [Se5]). 

In particular, these graded closure contain graded closures that collect test sequences for which there exist
extra rigid or strictly family elements that are not part of the proof systems. Given such an Extra  graded closure we impose one of finitely many 
collapse Diophantine conditions that force the extra rigid or strictly solid element to be either not rigid or not strictly solid, or to be part of a family
that appears in the proof statement.
By adding elements that demonstrate these collapse Diophantine conditions we get finitely many $Collapse$ limit groups.

We analyze its (Duo) resolutions
using the analysis of quotient resolutions, that is used in the sieve procedure [Se6]. However, when we analyze the
resolutions of a Collapse limit group, we analyze only those resolutions that have a similar structure as that of
$Cduo$, i.e., that are built from a completion to which two other completions are amalgamated and these two
completions are closures of the two completions that are amalgamated to the completion of a resolution of
$<\hat f, d_0,\hat d_0>$ in the construction of $Cduo$. 

By continuing this construction iteratively, according to
the steps of the sieve procedure, we finally obtain a finite collection of Duo limit groups, that we denote
$Sduo_1,\ldots,Sduo_h$. Each of these Duo limit groups is constructed from a completion to which we amalgamate two
closures of the completions that are amalgamated in the construction of $Cduo$. The sieve procedure that was
used to construct the Duo limit groups, $Sduo_1,\ldots,Sduo_h$, guarantees that they have the following 
properties.

\vglue 1pc
\proclaim{Proposition 3.5} The Duo limit groups, 
$Sduo_1,\ldots,Sduo_h$, that are associated with the Duo limit groups, $Cduo_1,\ldots,Cduo_v$, 
have the following properties:
\roster
\item"{(1)}" Each of the Duo limit groups, $Sduo_i$, is constructed from a completion that we denote $B(Sduo_i)$,
that contains the subgroup $<\hat f,d_0,\hat d_0>$, to which we amalgamate two completions that are closures
of the completions that are amalgamated to the base completion in the 
associated Duo limit group $Cduo$. We denote these two
completions $Cl_p(Sduo_i)$ and $Cl_q(Sduo_i)$.

\item"{(2)}" there is a homomorphism from the associated Duo limit group $Cduo$ into $Sduo_i$ that maps
the base completion in $Cduo$ into $B(Sduo_i)$, and
the two completions that are amalgamated to the base completion in $Cduo$ into their closures, $Cl_p(Sduo_i)$
and $Cl_q(Sduo_i)$, so that the map preserves the level structure of the two completions.

\item"{(3)}" The base completion $B(Sduo_i)$ terminates in either a rigid or a solid limit group with
respect to the subgroup, $<\hat d_0>$.

\item"{(4)}" for each $Sduo_i$ there exists two maps from either one or two  of the limit groups,
$Ipr_1,\ldots,Ipr_w$, into $Sduo_i$. Their images are the subgroup $<f,p,q,a>$, and a subgroup of
$<\hat f,d_0,\hat d_0>$ that we denote $<\hat f, q_0,\hat q_0,a>$.
\endroster
\endproclaim

\nfp All the properties (1)-(4) follow in a straightforward way from the construction of the duo limit groups,
$Sduo_1,\ldots,Sduo_h$, from the duo limit groups, $Cduo_1,\ldots,Cduo_v$.

\line{\hss$\qed$}

As we did with the Duo limit groups $Tduo$, for presentation purposes we assume
that the terminal limit groups of each of the Duo limit groups, $Sduo_1,\ldots,Sduo_h$, are rigid,
and there are no abelian vertex groups that appear in any of the levels of the 3 graded completions from which 
each of the duo limit groups, $Sduo_1,\ldots,Sduo_h$, is constructed.

The construction of the Duo limit groups, $Sduo_1,\ldots,Sduo_h$, allows us to present the construction of
uniformization limit groups, that are the main tool that we use in order to obtain separation of variables,
with which we will eventually be able to construct the set of parameters
that are associated with the equivalence classes that are defined by the given equivalence relation $E(p,q)$.

\vglue 1pc
\proclaim{Definition 3.6}  
Let $Sduo$ be one of the obtained limit groups, $Sduo_1,\ldots,Sduo_h$.
$Sduo$ is by construction a Duo limit group, and it contains a base (graded) completion that contains the
subgroup
$<\hat f,d_0,\hat d_0>$, 
and on top of that completion, there are two amalgamated  graded closures, one containing the
subgroup $<d_1,p>$ that we denote $Cl_p(Sduo)$, and the second containing the subgroup $<d_2,q>$ that we
denote $Cl_q(Sduo)$. We further denote the completion obtained from the base completion in $Sduo$, $B(Sduo)$,
to which we amalgam the closure, $Cl_p(Sduo)$, $GC_p(Sduo)$, and the completion obtained from $B(Sduo)$ 
to which we amalgam the closure, $Cl_q(Sduo)$, $GC_q(Sduo)$.

Starting with a limit group $Sduo$, and its associated graded completion, $GC_p(Sduo)$,
we can apply the construction of the Duo envelope, that is presented in theorem 1.7, and associate with
$GC_p(Sduo)$ a finite collection of Duo limit groups, so that one of the graded completions that 
is associated with these Duo limit groups is a graded closure of $GC_p(Sduo)$, 
and the distinguished vertex groups of these Duo limit groups contain the distinguished vertex group of 
$GC_p(Sduo)$,
 and in particular contain the subgroup $<\hat d_0>$.

Since the graded completion that we started the construction with, $GC_p(Sduo)$, is contained in $Sduo$, 
and $Sduo$ is obtained from $GC_p(Sduo)$ by amalgamating to it the closure, $Cl_q(Sduo)$, 
each of the Duo limit groups
that is obtained from $GC_p(Sduo)$ using the construction of the Duo envelope (theorem 1.7) can be
extended to $Sduo$ itself. i.e., with $Sduo$ we associate finitely many limit groups, where each of these limit
groups is obtained from $Sduo$ by amalgamating it with the second graded completion (the one that is associated 
with the subgroup $<\tilde q>$) of one of the Duo limit groups that are constructed from it. We call the obtained
limit groups $uniformization$ $limit$ $groups$, and denote their entire collection (i.e., all the limit groups 
of this form that are obtained from the various Duo limit groups $Sduo$), $Unif_1,\ldots,Unif_d$.
\endproclaim

Uniformization limit groups is the main tool that will serve us to obtain separation of variables,
which will eventually enable us to find the class functions that we are aiming
for.
 i.e., class functions that separate classes and associate a bounded set of elements with
each equivalence class (up to the basic equivalence relations). 

As we did with the Duo limit groups $Tduo$, for presentation purposes we assume
that the terminal limit groups of each of the Duo limit groups that are associated with the various
completions, $GC_p(Sduo_i)$, are rigid (and not solid) with respect to the subgroup $<\hat d_0>$, and that
the graded completions that are associated with these duo limit groups contain no abelian vertex groups in any of their
levels. Since
uniformization limit groups were constructed from these Duo limit groups, this implies
that the terminal limit group of all the uniformization limit groups are rigid as well, and the graded completions
that are associated with the constructed uniformization limit groups contain no abelian vertex groups in any
of their levels. Later on we 
generalize our arguments and omit these assumptions.

Before we use the uniformization limit groups in further constructions we list some of their basic properties
that will assist us in the sequel. These have mainly to do with the various maps from the rigid and solid
limit groups, $Ipr_1,\ldots,Ipr_w$, into uniformization limit groups.

\vglue 1pc
\proclaim{Proposition 3.7} Let $Unif$ be one of the uniformization limit groups, $Unif_1,\ldots,Unif_d$. Then: 
\roster 
\item"{(1)}" With $Unif$ there are 3 associated maps from the various rigid and solid limit groups,
$Ipr_1,\ldots,Ipr_w$ into $Unif$. Two of these maps are associated with and are into the Duo limit
group $Sduo$ from which $Unif$ was constructed. 
The images of these two maps are the subgroup $<f,p,q,a>$, and a subgroup of
$<\hat f,d_0,\hat d_0>$ that we denote $<\hat f, q_0,\hat q_0,a>$. The third homomorphism is from one of the
rigid or solid limit group, $Ipr_1,\ldots,Ipr_w$, into the Duo limit group that was constructed from 
$GC_p(Sduo)$, and from which $Unif$ was constructed. We denote the image of this third map,
$<\tilde f, p,\tilde q,a>$.
Furthermore, none of these 3 maps factor
through a free product of limit groups in which $<p>$ is contained in one factor,  $<q>$ or $<\tilde q>$ is contained in the second
factor and $F_k$ is contained in one of the two factors.

\item"{(2)}" there exist generic points of $Unif$, i.e.,  sequences of specializations of $Unif$ that
are composed from test sequences of the 4 completions from which $Unif$ is built, ($Cl_p$, $Cl_q$, $B(Sduo)$), and
the completion that is associated with $<\tilde q>$ in the Duo limit group that is associated with 
$GC_p(Sduo)$ from which $Unif$ was constructed), for which the restrictions to specializations of the 3
subgroups, $<f,p,q,a>$, $<\hat f,q_0,\hat q_0,a>$, and 
$<\tilde f, p,\tilde q,a>$, restrict to proofs that the specializations of the couples $(p,q)$,$(q_0,\hat q_0)$,
and $(p,\tilde q)$ are in the same equivalence class of $E(p,q)$. Furthermore, these test sequences 
restrict to proofs that the specializations of the couples $(p,q_0)$ and $(q_0,q)$ are in the same equivalence 
class of $E(p,q)$, and these test sequences restrict to distinct sequences of specializations $\{(p_n,q_n\}$.

\item"{(3)}" If the subgroups $<p>$ and $<q>$ in $Unif$ are non-trivial, then the subgroup 
$<\tilde f,p,\tilde q,a>$ intersects non-trivially some conjugates of the
terminal rigid vertex group in $Unif$. 
Let $\tilde H^1,\ldots,\tilde H^c$ be the (conjugacy classes of) subgroups of intersection.

Then there exists a global bound $U$, so that
for every possible value of $\hat d_0$ for which there exists a test sequence with the properties that
are described in part (2), there are at most $U$ possible conjugacy classes of specializations of
the subgroups $\tilde H^1,\ldots,\tilde H^c$ that extend 
$\hat d_0$ and together are restrictions of conjugates of a rigid specialization of 
the terminal limit group of $Unif$. 
\endroster
\endproclaim

\nfp Parts (1) and (2) follow from the construction of the uniformization limit group, $Unif$, and from the
construction of the duo limit group, $Sduo$, from which it was constructed. Part (3) follows from the uniform bound
on the number of rigid specializations of a rigid limit group with the same value of the defining parameters, that
was proved in theorem 2.5 in [Se3].

\line{\hss$\qed$}

Uniformization limit groups are constructed as an amalgamation of a (Duo) limit
group $Sduo$, and a Duo limit group that was associated with a graded completion, $GC_p(Sduo)$, that is 
contained in $Sduo$. This structure of uniformization limit groups enable them to reflect properties of
generic points in fibers on one side (the $Sduo$ side) and universal properties (w.r.t.\ duo families)  on the other side.

\vglue 1pc
\proclaim{Theorem 3.8} Let $Unif$ be one of the uniformization  limit groups that are associated with the given
definable equivalence relation $E(p,q)$. Then:
\roster
\item"{(i)}" Let $DE$ be the distinguished vertex group in $Tduo$, the Duo limit
group from which $Sduo$, the Duo limit group that is associated with $Unif$,  was constructed. 
Let $H^1,\ldots,H^e$ be the (conjugacy classes) of intersections between the subgroup $<f,p,q,a>$ and 
conjugates of the distinguished vertex group of $Tduo$, $DE$.
Recall that by part (3) of proposition 3.3, if both subgroups $<p>$ and $<q>$ are non-trivial, then at least one of these 
subgroups of intersection, and all but at most one (conjugacy classes), are non-trivial.

Let $\hat d_0'$ be a specialization of $\hat d_0$ 
that does not belong to the finitely many equivalence
classes of $E(p,q)$ that are singled out in theorem 3.1. Suppose that with the equivalence class of
$\hat d_0'$ there are infinitely many 
specializations of $H^1,\ldots,H^e$ 
that are associated with the equivalence class of $\hat d_0'$, and  these
specializations of $H^1,\ldots,H^e$ can be extended to conjugates of rigid specializations of $DE$ 
(the terminal limit group in $Tduo$). Furthermore,  these rigid specializations of $DE$
can be extended to test sequences of $Tduo$ that restrict to valid proofs that
the corresponding couples: $(p_n,q_n)$, and $(p_n,q_0)$ (in $Tduo$), belong to the same equivalence
class of $E(p,q)$, and each of these test sequences restrict to an infinite distinct sequence of couples,
$\{(p_n,q_n)\}$. Then there exists a Duo limit group  $Sduo_1$ that is constructed from $Tduo$, 
in which 
$H^1,\ldots,H^e$ are not all contained in the distinguished vertex of $Sduo_1$.

If there are in addition infinitely many conjugacy classes of specializations of the subgroups $H^1,\ldots,H^e$
with the same properties, then at least one of the subgroups $H^1,\ldots,H^e$ is not contained in a conjugate of
the distinguished vertex in $Sduo_1$.

\item"{(ii)}" 
Let $\hat d_0'$ be a specialization of $\hat d_0$ that does not belong to the finitely many
equivalence classes that are singled out in theorem 3.1. Suppose that there are only finitely many
conjugacy classes of specializations of $H^1,\ldots,H^e$ that are associated with the equivalence class of 
$\hat d_0'$, so that these
specializations of $H^1,\ldots,H^e$ can be extended to conjugates of rigid specializations of $DE$ 
(the terminal limit group in $Tduo$)
that can be extended to test sequences of $Tduo$ that restrict to valid proofs that
the corresponding couples: $(p_n,q_n)$, and $(p_n,q_0)$, belong to the same equivalence
class of $E(p,q)$, and so that these test sequences restrict to sequences of distinct couples, $\{(p_n,q_n)\}$.
Then there exists a global bound (that does not depend on $\hat d_0'$ nor on its class) on the possible
values of the conjugacy classes of the subgroups $H^1,\ldots,H^e$  that can extend such a specialization 
of the elements $d_0$ which is in the 
equivalence class of $\hat d_0'$.

\item"{(iii)}" Let $\hat d_0'$ be a specialization of the elements $\hat d_0$, so that $\hat d_0'$ restricts to
specializations $\hat q_0'$, and $\hat q_0'$ does not belong to the finitely many equivalence classes 
of $E(p,q)$ that
are excluded in theorem 3.1. Suppose that $\hat d_0'$ extends to a test sequence of $Sduo$ that restricts to proofs
that the couples $(p_n,q_0(n))$, and $(q_0(n),\hat q_0')$ belong to the same equivalence class in $E(p,q)$. 

Let $\tilde q'$ be a specialization of the elements $\tilde q$ that belongs to  
the equivalence class of $\hat q_0'$.
Then  $\tilde q'$ extends to
a specialization of at least one of the uniformization limits groups $Unif$ that are constructed from $Sduo$, 
so that for this uniformization limit group $Unif$, $\tilde q'$ extends to a sequence of specializations,
$\{(\tilde f,p_n,\tilde q')\}$ that prove that the sequence of couples $\{(p_n,\tilde q')\}$ is in the given definable
set $E(p,q)$, and the sequence $\{p_n\}$ is a test sequence for the completion, $GC_p(Sduo)$, that is contained in $Sduo$.
Furthermore, the specializations of the corresponding test sequence of $GC_p(Sduo)$ restrict to proofs that 
the couples $\{(p_n,q_0(n))\}$, and $\{(q_0(n),\hat q_0')\}$ are in the same equivalence class of $E(p,q)$. 
\endroster
\endproclaim

\nfp Part (i) follows from the constructions of the duo limit group, $Sduo$, and the uniformization
limit groups that are associated with it. Part (ii) follows from the construction of a uniformization limit group,
and from the existence of a uniform bound on the number of rigid specializations of a rigid limit group, for
any given value of the defining parameters (theorem 2.5 in [Se3]). Part (iii) follows since the uniformization limit
groups that are associated with the duo limit group $Sduo$, collect all the test sequences of $GC_p(Sduo)$ and all the
specializations of the subgroup $<q>$, so that the restriction of the test sequences of $GC_p(Sduo)$ to the subgroup
$<p>$, and the values 
of the subgroup, $<q>$,  can be extended 
to a sequence of specializations of one of the rigid or solid limit groups, $Ipr_1,\ldots,Ipr_w$. If $\tilde q'$ belong
to the same equivalence class as $\hat q_0'$, such a test sequence clearly exists for a test sequence of 
$GC_p(Sduo)$ and $\tilde q'$, and part (iii) follows.

\line{\hss$\qed$}

By construction, Uniformization limit groups admit 3 different maps from the rigid and solid limit groups, 
$Ipr_1,\ldots,Ipr_w$, into them. These maps prove that (generic specializations of)
the couples, $(p,q)$, $(q_0,\hat q_0)$, and $(p,\tilde q)$,
are in the definable set $E(p,q)$. The structure of the uniformization limit groups, and in particular their
ability to use both generic points of their associated Duo limit group $Sduo$, and the universality 
property of the collection of specializations of the elements $\tilde q$, that is associated with
the (finite) collection of uniformization limit groups that is associated with a Duo limit group $Sduo$,
enable us to "compare" between two of these maps, those that verify that generic specializations of the 
couples $(p,q)$ and $(p, \tilde q)$ are in $E(p,q)$. 
This comparison is crucial in our approach
to constructing the desired class functions from the given definable equivalence relation $E(p,q)$.

Recall that by part (3) of proposition 3.7, if the subgroups $<p>$ and $<q>$ are non-trivial, then the subgroup
 $<\tilde f,p,\tilde q>$ of a uniformization limit group $Unif$, intersects non-trivially conjugates of
the distinguished vertex
group in $Unif$, in (conjugacy classes of) the subgroups:  
$\tilde H^1,\ldots,\tilde H^c$. 

The collection of duo limit groups, $Tduo_1,\ldots,Tduo_m$, collect
all the possible extensions of test sequences of the graded completions, $G_1,\ldots,G_t$,
that form the Diophantine envelope of $E(p,q)$, to rigid and almost shortest strictly solid specializations of the
rigid and solid limit groups, $Ipr_1,\ldots,Ipr_w$. With the two subgroups of a uniformization limit group, 
$Unif$, $<f,p,q,a>$ and $<\tilde f,p,\tilde q,a>$, 
we can naturally associate with $Unif$ two (possibly identical) of the duo limit groups, $Tduo_1,\ldots,Tduo_m$. By 
construction with $<f,p,q,a>$ we can associate the duo limit group $Tduo$. With $<\tilde f,p,\tilde q,a>$ it is possible to
associate another (possibly the same) duo limit group from the collection, $Tduo_1,\ldots,Tduo_m$, 
that we denote, $\tilde Tduo$.

\smallskip
Let $Tduo$ be one of the duo limit groups, $Tduo_1,\ldots,Tduo_m$, and let
$DE$ be the distinguished vertex in $Tduo$.  By the construction of $Tduo$,  there is an associated map from 
one of the 
rigid and solid limit groups, $Ipr_1,\ldots,Ipr_w$, into it, that we denoted $<f,p,q,a>$.
$Tduo$ is a duo limit group, and by our assumptions it terminates in a rigid limit group and the
two
graded completions that are associated with it contain no non-cyclic abelian vertex group
in any of their levels. Hence, $Tduo$, can be presented as an amalgamated product:
$Tduo = <d_1,p>*_{<d_0>} \, <d_2,q>$. If both subgroups, $<p>$ and $<q>$, are non-trivial in $Tduo$, then
the subgroup $<f,p,q,a>$ intersects non-trivially some conjugates of the distinguished vertex group in $Tduo$.  
We denoted by $H^1,\ldots,H^e$ the conjugacy classes of these intersections.

Let $\hat d_0'$ be a specialization of $d_0$ (i.e., a specialization of a fixed generating set of $DE$)
that does not belong to the finitely many equivalence
classes of $E(p,q)$ that are singled out in theorem 3.1. Suppose that  there are infinitely many
conjugacy classes of specializations of the subgroups $H^1,\ldots,H^e$ that are associated with the 
equivalence class of $\hat d_0'$, so that these
specializations of $H^1,\ldots,H^e$ can be extended to conjugates of rigid specializations of $DE$ 
(the terminal limit group in $Tduo$)
that can be extended to test sequences of $Tduo$ that restrict to valid proofs that
the corresponding couples: $(p_n,q_n)$, and $(p_n,q_0)$, belong to the same equivalence
class of $E(p,q)$, and to distinct couples of specializations: $\{(p_n,q_n)\}$. 

By part (i) of theorem 3.8 there exists a Duo limit group  $Sduo$, 
that is constructed from $Tduo$, 
in which not all the images of the subgroups 
$H^1,\ldots,H^e$ can be conjugated into the distinguished vertex in $Sduo$, 
and so that $\hat d_0'$ can be extended to a test sequence of $Sduo$ that restricts to
proofs that the couples: 
$(p_n,q_n)$, $(p_n,q_0(n))$, and $(q_0(n),\hat q_0')$ belong to the same equivalence
class of $E(p,q)$,  the corresponding specializations of the subgroups $H^1,\ldots,H^e$ belong to distinct
conjugacy classes, and the corresponding couples of specializations $\{(p_n,q_n)\}$ are distinct. 

By the construction of the uniformization  limit groups $Unif$, 
since this last conclusion holds for $Sduo$, it holds for at least one of the uniformization 
limit groups, $Unif$, that are associated with it.

Let $Unif$ be such a uniformization limit group, i.e., a uniformization limit group in which not all the subgroups.
$H^1,\ldots,H^e$, can be conjugated into the distinguished vertex group in $Unif$.
There are two maps
of the limit groups $Ipr_1,\ldots,Ipr_w$ into $Unif$, with images $<f,p,q,a>$ and $<\tilde f,p,\tilde q,a>$, that
are associated with two (possibly identical) duo limit groups, $Tduo$ and $\tilde Tduo$.

Recall that the duo limit groups, $Tduo_1,\ldots,Tduo_m$, encode all the extensions of test sequences
of the graded completions, $G_1,\ldots,G_t$, that form the Diophantine envelope of $E(p,q)$, to
rigid and almost shortest strictly solid specializations of the  rigid and solid
limit groups, $Ipr_1,\ldots,Ipr_w$, 
that do not factor through a free product in
which $<p>$ is contained in one factor and $<q>$ is contained in the second factor,
so that these extended test sequences restrict to 
valid proofs that the sequences of couples, $\{(p_n,q_n)\}$, are in the set $E(p,q)$,
and the  couples, $\{(p_n,q_n)\}$, are distinct.

Suppose that the second map of the limit groups, $Ipr_1,\ldots,Ipr_w$, into $Unif$, 
the one with image $<\tilde f,p,\tilde q,a>$,
is associated with $Tduo$ as well (i.e., $\tilde Tduo =Tduo$). Furthermore, suppose that the
images of the subgroups $H^1,\ldots,H^e$ in $<f,p,q,a>$ and in $<\tilde f,p,\tilde q,a>$ are conjugate. Then, by the
construction of the uniformization limit group $Unif$,  the two images of at least one of 
these subgroups can not be conjugated into the distinguished vertex group in $Unif$.

Suppose that the map from one of the limit groups, $Ipr_1,\ldots,Ipr_w$, into $Unif$, with image
$<\tilde f,p,\tilde q,a>$, is associated with a duo
limit group $Tduo_i$ which is not $Tduo$, or that it is associated with $Tduo$, but  the images of the subgroups 
$H^1,\ldots,H^e$ under the two maps from $Ipr_1,\ldots,Ipr_w$ to $Unif$,
with images, $<f,p,q,a>$ and $<\tilde f,p, \tilde q,a>$, are not conjugate in $Unif$.

By proposition 3.4, under our assumption that the terminal limit groups of the duo limit groups,
$Tduo_1,\ldots,Tduo_m$ are rigid, for each specialization $d_0'$ of a (finite) generating set $d_0$ 
of the distinguished vertex group $DE=<d_0>$ of $Tduo$,
there are at most boundedly many possible conjugacy classes
of specializations
of the subgroups, $H^1,\ldots,H^e$, that can extend $d_0'$, so that there exists a test sequence
of the duo family that is associated with
$d_0'$, that can be extended to shortest rigid and strictly solid specializations of one of the limit groups,
$Ipr_1,\ldots,Ipr_w$, that is associated with $Tduo$ and these specializations will have one of the given conjugacy
classes of specializations of the subgroups, $H^1,\ldots,H^e$.

Therefore, if the subgroup $<\tilde f,p,\tilde q,a>$ is not associated with $Tduo$, or it is associated with
$Tduo$, but the images of the subgroups $H^1,\ldots,H^e$ in $<f,p,q,a>$ and $<\tilde f,p,\tilde q,a>$ are not
pairwise conjugate, two of boundedly many such maps are already present,
and later on we will be able to apply the pigeon hole principle, to argue that after boundedly many steps 
two maps with pairwise conjugate subgroups $H_i^1,\ldots,H_i^{e_i}$ (that are associated with one
of the duo limit groups, $Tduo_1,\ldots,Tduo_m$) must be present. This will eventually guarantee the termination
of an iterative procedure that we present, that will finally give us the parameters for the equivalence classes of
the given equivalence relation $E(p,q)$.

\medskip
<Given the definable equivalence relation $E(p,q)$, we started its analysis with its Diophantine 
envelope, $G_1,\ldots,G_t$, and Duo envelope, $Duo_1,\ldots,Duo_r$ (theorems 1.5 and 1.7). 
We further associated with $E(p,q)$ the rigid
and solid limit groups, $Ipr_1,\ldots,Ipr_w$, so that their rigid and strictly solid specializations
(with respect to the parameter subgroup $<p,q,a>$) restrict to valid proofs that the couples $(p,q)$ are
in $E(p,q)$, and these specializations do not factor through a free product in which $<p>$ is contained in one
factor,  $<q>$ is contained in the second factor and $F_k$ is contained in one of the two factors, for all but finitely many equivalence classes of
$E(p,q)$ (theorem 3.1). Then we collected all the possible extensions of
test sequences of the graded completions,
$G_1,\ldots,G_t$, that form the Diophantine envelope of $E(p,q)$, to rigid and almost shortest strictly
solid specializations of the rigid and solid limit groups,
$Ipr_1,\ldots,Ipr_w$, 
and these were collected by the (finite) collection of duo limit groups, $Tduo_1,\ldots,Tduo_m$
(see propositions 3.3 and 3.4). 

Into each of the duo limit groups, $Tduo_1,\ldots,Tduo_m$, there is an associated  map of one of
the rigid and solid limit groups, $Ipr_1,\ldots,Ipr_w$. We denoted the image of the map from one of the 
rigid or solid limit
groups, $Ipr_1,\ldots,Ipr_w$, into a duo limit group, $Tduo$, by $<f,p,q,a>$. 
If both subgroups $<p>$ and $<q>$ are non-trivial in $Tduo$, then by proposition 3.3  the subgroup
$<f,p,q,a>$ intersects
some conjugates of the distinguished vertex group in $Tduo$ non-trivially. We denoted by $H^1,\ldots,H^e$
the conjugacy classes of these intersection subgroups. 
be non-trivial.

If the number of 
specializations of conjugacy classes of specializations of the subgroups,
$H^1,\ldots,H^e$ for a given equivalence class is finite then it is globally bounded (for all such equivalence classes). 
For the entire collection
of equivalence classes for which the number of conjugacy classes of specializations  of $H^1,\ldots,H^e$ is infinite, 
we have 
associated with the duo limit groups, $Tduo_1,\ldots,Tduo_m$, a finite collection of duo limit groups
$Sduo_1,\ldots,Sduo_h$.

With each of the duo limit groups, $Sduo$, we have associated a finite collection of uniformization limit groups,
that we denoted, $Unif_1,\ldots,Unif_d$.
Each of these uniformization limit groups admits a 
second 
map from one of the rigid and solid limit groups,
that we denote
$<\tilde f,p,\tilde q,a>$. By the universality of the collection of duo limit groups, $Tduo_1,\ldots,Tduo_m$, with each
of the subgroups, $<f,p,q,a>$ and $<\tilde f,p,\tilde q,a>$, it is possible to associate one of these limit groups. By construction
$Tduo$ is associated with  $<f,p,q,a>$, and we denoted $\tilde Tduo$, the duo limit group (from the collection
$Tduo_1,\ldots,Tduo_m$) that is associated with the subgroup $<\tilde f,p,\tilde q,a>$  ($Tduo$ and $\tilde Tduo$ may be the same
duo limit group).

If both subgroups $<p>$ and $<\tilde q>$ in $Unif$ are non-trivial, then the subgroup $<\tilde f,p,\tilde q,a>$
intersects some conjugates of the distinguished vertex group in $Unif$ non-trivially.
We denote by $\tilde H^1,\ldots,\tilde H^c$ the conjugacy classes of these subgroups of intersection.

To associate parameters with the various equivalence classes of the given equivalence relation $E(p,q)$,
and obtain separation of variables,
we repeat these constructions iteratively. The constructions that we perform in the
 second step of the iterative procedure,
depend on whether the two maps from the rigid and solid limit groups, $Ipr_1,\ldots,Ipr_w$, into
a uniformization  limit group $Unif$, with images $<f,p,q>$ and $<\tilde f,p,\tilde q,a>$, are associated with the same
duo limit group $Tduo$ (i.e., if $\tilde Tduo=Tduo$),
and if so whether the images of the subgroups $H^1,\ldots,H^e$ under the two associated maps,
are pairwise conjugate, or not.

Suppose that the two maps with images, $<f,p,q,a>$ and $<\tilde f,p,\tilde q,a>$, are associated with the same
duo limit group $Tduo$,
and the images of the subgroups, $<H^1,\ldots,H^e>$, under the two
associated maps, are pairwise conjugate.
Suppose
that there exists an equivalence class of $E(p,q)$, which is not one of the finitely many equivalence
classes that were singled out in theorem 3.1, for which there exist an infinite sequence of specializations
of the subgroups, $\tilde H^1,\ldots,\tilde H^c$, that are not pairwise conjugate,  that can be extended to test 
sequences of a uniformization limit group $Unif$ (i.e., sequences that restrict 
to test sequences of the 4 completions from which the uniformization limit group $Unif$ is composed), 
that restrict to 
pairwise non-conjugate specializations of the subgroups, $H^1,\ldots,H^e$, and so that the restrictions 
of these test sequences,
$\{(f_n,p_n,q_n,a)\}$ and $\{(\tilde f_n,p_n,\tilde q_n,a)\}$, prove that the couples,
$\{(p_n,q_n)\}$ and $\{(p_n,\tilde q_n)\}$, are in  $E(p,q)$, and the restrictions 
$\{(\hat f_n,q_0(n),\hat q_0,a)\}$ prove that the couples $\{(q_0(n),\hat q_0)\}$
are in $E(p,q)$. Furthermore these test sequences restrict to valid proofs that the couples $\{(p_n,q_0(n))\}$ 
are in $E(p,q)$ (recall that $q_0(n)$ is the restriction of the specializations $d_0(n)$ to
the elements $q_0$ and $\hat q_0$ is the restriction of the specialization $\hat d_0$), and the couples,
 $\{(p_n,q_n)\}$ and $\{(p_n,\tilde q_n)\}$,
are distinct.

We collect all these equivalence classes and their associated specializations of the subgroup $\tilde H$
in a finite collection of uniformization limit groups, in a similar way to the 
construction of the uniformization limit groups, $Unif_1,\ldots,Unif_d$. 
With these we associate a finite collection of
duo limit groups, in a similar way to the construction of the duo limit groups, $Cduo_1,\ldots,Cduo_v$.
Then we apply the sieve procedure (that is presented in [Se6]), to associate with this collection of duo
limit groups, a finite collection of duo limit groups, in a similar way to the construction of the
duo limit groups, $Sduo_1,\ldots,Sduo_h$. With each of these duo limit groups we associate a finite
collection of uniformization limit groups. 
We denote these uniformization
limit groups, $Unif^2_1,\ldots,Unif^2_{d^2}$. For presentation purposes we continue to assume that the
terminal limit groups of all the constructed duo limit groups are rigid (and not solid), and the graded completions
that are associated with them contain no abelian vertex groups in any of their levels.

Let $Unif^2$ be one of the duo limit groups,  
$Unif^2_1,\ldots,Unif^2_{d^2}$, in which at least one of the subgroups, $\tilde H^1,\ldots,\tilde H^c$,
can not be conjugated into the distinguished vertex group in $Unif^2$.
By construction, there are 3 maps from the rigid and solid limit groups, $Ipr_1,\ldots,Ipr_w$,
into $Unif^2$. Two of these maps with images, $<f,p,q,a>$ and $<\tilde f,p,\tilde q,a>$, 
are inherited from the associated 
uniformization limit group $Unif$, and are associated with the duo limit group $Tduo$ by our assumptions. 
The third map with image that we denote, $<f',p,q',a>$, is  
also associated with one of the duo limit
groups, $Tduo_1,\ldots,Tduo_m$. We denote the duo limit group with which $<f',p,q'>$ is associated, $Tduo'$.   

By our assumptions the two subgroups, $<f,p,q,a>$ and $<\tilde f,p,\tilde q,a>$, are associated with the same duo limit group $Tduo$,
and the images
of the subgroups, $H^1,\ldots,H^e$, under these two maps are pairwise conjugate. Suppose that the third map
from one of the subgroups, $Ipr_1,\ldots,Ipr_w$, into $Unif^2$ (with image $<f',p,q',a>$), is associated
with $Tduo$ as well (i.e., $Tduo'=Tduo$), 
and the images of
the subgroups, $H^1,\ldots,H^e$, under the third map are pairwise conjugate to their images under the first two maps.

The uniformization limit group, $Unif^2$, was constructed from a uniformization limit group,
$Unif$, and its associated duo limit group $Sduo$. Hence, by the universality
of the collection of uniformization limit groups that are associated with the duo limit group,
$Sduo$,  with the third map from one of the rigid
or solid limit groups, $Ipr_1,\ldots,Ipr_w$, into $Unif^2$, with image $<f',p,q',a>$, we can associate another
(possibly the same) uniformization limit group that is associated with the duo limit group, $Sduo$, that
we denote $Unif'$.

Suppose that the uniformization limit group, $Unif'$,  that is associated with
$Unif^2$, is $Unif$, the uniformization limit group from which $Unif^2$ was constructed.
Suppose further, that the images of the subgroups, $\tilde H^1,\ldots,\tilde H^c$, in the subgroup $<f',p,q',a>$,
are pairwise conjugate to the images of these subgroups in the subgroup, $<\tilde f,p,\tilde q,a>$.

If both subgroups, $<p>$ and $<q>$, are non-trivial, then the subgroup $<f',p, q',a>$ 
of the uniformization limit group $Unif^2$, intersects non-trivially some conjugates of
the distinguished vertex group in $Unif^2$. 
Let ${H^1}',\ldots,{H^b}'$ be conjugacy classes of these subgroups of intersection.

If the uniformization limit group that is associated with
the subgroup $<f',p, q',a>$ and the uniformization limit group $Unif^2$,
 is not $Unif$ (i.e., if $Unif'$ is not $Unif$), or if 
it is $Unif$, but the
images of the subgroups, $\tilde H^1,\ldots,\tilde H^c$, in the subgroup $<f',p,q',a>$,
are not pairwise conjugate to the images of these subgroups in the subgroup, $<\tilde f,p,\tilde q,a>$,
then as we argued for the uniformization limit group $Unif^2$, the two maps with images  
 $< \tilde f,p,\tilde q,a>$ and $<f',p,q',a>$,
occupies two of the boundedly many possibilities of such  maps (where the bound is uniform and does
not depend on the specific equivalence class of $E(p,q)$).

By our assumptions the two subgroups, $<f,p,q,a>$ and $<\tilde f,p,\tilde q,a>$, are associated with the same duo limit group $Tduo$,
and the images
of the subgroups, $H^1,\ldots,H^e$, under these two maps are pairwise conjugate. If the third map
from one of the subgroups, $Ipr_1,\ldots,Ipr_w$, into $Unif^2$ (with image $<f',p,q',a>$), is not associated
with $Tduo$ as well (i.e., $Tduo'$ is not $Tduo$), 
or if the images of
the subgroups, $H^1,\ldots,H^e$, under the third map are not pairwise conjugate to their images under the first two maps,
then by the same reasoning, the two maps with images $<f,p,q,a>$ and $<f',p,q',a>$, 
occupies two of the boundedly many possibilities of such  maps (where the bound is uniform and does
not depend on the specific equivalence class of $E(p,q)$).

\smallskip
Suppose that the two maps with images, $<f,p,q,a>$ and $<\tilde f,p,\tilde q,a>$, are not associated with the same
duo limit group $Tduo$, or that they are both associated with $Tduo$, and
the images of the subgroups, $<H^1,\ldots,H^e>$, under the two
associated maps, are not pairwise conjugate. In this case, the second map, the one with image 
$<\tilde f,p,\tilde q,a>$,
is associated with a duo limit group $\tilde Tduo$. If both subgroups, $<p>$ and $<q>$, in $\tilde Tduo$, are
non-trivial, then  
the image of the map from one of the rigid or solid limit groups,
$Ipr_1,\ldots,Ipr_w$, in $\tilde Tduo$, intersects non-trivially some conjugates of the distinguished vertex group
in $\tilde Tduo$. 
Let $\hat H^1,\ldots, \hat H^a$ be the conjugacy classes of these subgroups of intersection. 

Suppose
that there exists an equivalence class of $E(p,q)$, which is not one of the finitely many equivalence
classes that were singled out in theorem 3.1, for which there exist an infinite sequence of specializations
of the subgroups, $\hat H^1,\ldots,\hat H^a$, that are not pairwise conjugate,  that can be extended to test 
sequences of a uniformization limit group $Unif$ (i.e., sequences that restrict 
to test sequences of the 4 completions from which the uniformization limit group $Unif$ is composed), 
that restrict to 
pairwise non-conjugate specializations of the subgroups, $H^1,\ldots,H^e$, and so that if both 
$<f,p,q>$ and $<\tilde f,p,\tilde q,a>$ are associated with $Tduo$, then the  specializations of
the subgroups, $H^1,\ldots,H^e$, are not pairwise conjugate to those of $\hat H^1,\ldots,\hat H^{a=e}$, and the
restrictions 
of these test sequences,
$\{(f_n,p_n,q_n,a)\}$ and $\{(\tilde f_n,p_n,\tilde q_n,a)\}$, prove that the couples
$\{(p_n,q_n)\}$ and $\{(p_n,\tilde q_n)\}$ are in  $E(p,q)$, and the restrictions 
$\{(\hat f_n,q_0(n),\hat q_0,a)\}$ prove that the couples $\{(q_0(n),\hat q_0)\}$
are in $E(p,q)$. Furthermore these test sequences restrict to valid proofs that the couples $\{(p_n,q_0(n))\}$ 
are in $E(p,q)$ (recall that $q_0(n)$ is the restriction of the specializations $d_0(n)$ to
the elements $q_0$ and $\hat q_0$ is the restriction of the specialization $\hat d_0$), and the couples $\{(p_n,q_n)\}$ and $\{(p_n,\tilde q_n)\}$
are distinct.

We collect all these equivalence classes and their associated specializations of the subgroups, 
$\hat H^1,\ldots,\hat H^a$,
in a finite collection of duo limit groups, in a similar way to the 
construction of the duo limit groups, $Dduo_1,\ldots,Dduo_u$. With these we associate a finite collection of
duo limit groups, in a similar way to the construction of the duo limit groups, $Cduo_1,\ldots,Cduo_v$.
Then we apply the sieve procedure (that is presented in [Se6]), to associate with this collection of duo
limit groups, a finite collection of duo limit groups, in a similar way to the construction of the
duo limit groups, $Sduo_1,\ldots,Sduo_h$. With each of these duo limit groups we associate a finite
collection of uniformization limit groups, that we denote (once again):
$Unif^2_1,\ldots,Unif^2_{d^2}$. For presentation purposes we continue to assume that the
terminal limit groups of all the constructed duo limit groups are rigid (and not solid), and the graded completions
that are associated with them contain no abelian vertex groups in any of their levels.

Let $Unif^2$ be one of the uniformization limit groups,  
$Unif^2_1,\ldots,Unif^2_{d^2}$ in which at least one of the subgroups, $\hat H^1,\ldots,\hat H^a$, and
at least one of the subgroups, $H^1,\ldots,H^e$,
can not be conjugated into the distinguished vertex group in $Unif^2$.
By construction, there are 3 maps from the rigid and solid limit groups, $Ipr_1,\ldots,Ipr_w$,
into $Unif^2$. Two of these maps with images, $<f,p,q,a>$ and $<\tilde f,p,\tilde q,a>$, are inherited from 
the associated 
uniformization  limit group $Unif$, and are associated with the duo limit groups $Tduo$ and $\tilde Tduo$ in correspondence, 
by our assumptions. 
The third map with image that we denote, $<f',p,q',a>$, is associated with the construction of $Unif^2$,
and is also associated with one of the duo limit
groups, $Tduo_1,\ldots,Tduo_m$. We denote the duo limit group with which $<f',p,q'>$ is associated, $Tduo'$.

By our assumptions the two subgroups, $<f,p,q,a>$ and $<\tilde f,p,\tilde q,a>$, are associated with the  duo limit groups
 $Tduo$ and $\tilde Tduo$ in correspondence, and if they are both associated with $Tduo$, then 
the images
of the subgroups, $H^1,\ldots,H^e$, under these two maps are not pairwise conjugate. 
Suppose that the third map
from one of the subgroups, $Ipr_1,\ldots,Ipr_w$, into $DPduo^2$ (with image $<f',p,q',a>$), is associated
with either $Tduo$ or $\tilde Tduo$,
and the images of
the subgroups, $H^1,\ldots,H^e$, or the subgroups, $\hat H^1,\ldots,\hat H^a$, 
under the third map are pairwise conjugate to their images under the first or the second map in correspondence.

If both subgroups, $<p>$ and $<q>$, are non-trivial in $Unif^2$, then the subgroup $<f',p, q',a>$ 
intersects some conjugates of the distinguished vertex group in $Unif^2$ non-trivially.
Let ${H^1}',\ldots,{H^b}'$ be the conjugacy classes of these subgroups of intersection.

By our assumptions, the two subgroups, $<f,p,q,a>$ and $<\tilde f,p,\tilde q,a>$, are associated with the  duo limit groups
 $Tduo$ and $\tilde Tduo$, and if $Tduo=\tilde Tduo$ then
the images
of the subgroups, $H^1,\ldots,H^e$, under these two maps are not pairwise conjugate. If the third map
from one of the subgroups, $Ipr_1,\ldots,Ipr_w$, into $Unif^2$ (with image $<f',p,q',a>$), is not associated
with $Tduo$ or $\tilde Tduo$, 
or if the images of
the subgroups, $H^1,\ldots,H^e$, or $\hat H^1,\ldots,\hat H^a$,  under the third map are not pairwise conjugate to their images under the first or the second map,
then three maps with images $<f,p,q,a>$, $<\tilde f,p,\tilde q,a>$, and $<f',p,q',a>$, 
occupies 3 of the boundedly many possibilities of such  maps (where the bound is uniform and does
not depend on the specific equivalence class of $E(p,q)$).

\medskip
We continue iteratively. Suppose that there exists a uniformization limit group $Unif^2$, and 
an equivalence class of $E(p,q)$, which is not one of the finitely many classes that are singled out
in theorem 3.1, for which there exist 
infinitely many conjugacy classes of specializations of the subgroups ${H^j}'$ that are associated with the
image of the third map from $Ipr_1,\ldots,Ipr_w$, into $Unif^2$ with image $<f',p,q',a>$, that
can be extended to test sequences of the duo limit group  $Unif^2$ (i.e., sequences that restrict 
to test sequences of the 4 completions from which the duo limit group $DPduo^2$ is composed), so that restrictions 
of these test sequences to the subgroups, $H^j$ and $\tilde H^j$ are pairwise non-conjugate, and the restrictions
of the subgroups:
$\{(f_n,p_n,q_n,a)\}$, $\{(\tilde f_n,p_n,\tilde q_n.,a)\}$ and $\{(f',p,q',a)\}$, prove that the couples
$\{(p_n,q_n)\}$, $\{(p_n,\tilde q_n)\}$ and $\{(p_n,q_n')\}$, are in  $E(p,q)$, and the restrictions of these test sequences prove that the couples,
$\{(p_n,q_0(n)\}$, $\{(q_0(n),\hat q_0(n))\}$ and $\{(\hat q_0(n),q_0')\}$ are in $E(p,q)$. 
Then we repeat these constructions, and obtain new uniformization limit groups, 
$Unif^3_1,\ldots,Unif^3_{d^3}$ that admit
4 maps from the limit groups $Ipr_1,\ldots,Ipr_w$ into each of them.

To obtain a set of parameters for the equivalence classes of the given definable equivalence relation, $E(p,q)$,
we need to ensure a termination of this iterative procedure, that we'll leave us with only finitely many
uniformization limit groups, $Unif^i$, and so that for each equivalence class (apart from the finitely many that are
singled out in theorem 3.1) there will exist a uniformization limit group $Unif^i$ with only boundedly many possible
conjugacy classes of
values for the associated subgroups $H^i$ ($H^i$ are obtained as intersections between conjugates of 
the distinguished vertex group
in $Unif^i$ and an associated image of one of the limit groups $Ipr_1,\ldots,IPr_w$, $<f^i,p,q^i,a>$). 

\vglue 1pc
\proclaim{Theorem 3.9} The iterative procedure for the construction of the uniformization limit groups,
$Unif^i$, terminates after finitely many steps.
\endproclaim

\nfp Suppose that the iterative procedure does not terminate after finitely many steps. Since at each step
finitely many uniformization limit groups are constructed,
by Konig's lemma, if the procedure doesn't terminate there must exist an infinite path along it.

Each uniformization limit group along the infinite path is equipped with a map from one of the 
limit groups, $Ipr_1,\ldots,Ipr_w$, into it (we denote this image $<f_j,p,q_j,a>$). 
Hence, by passing to a subsequence of the uniformization 
limit groups along the infinite path, we may assume that they are all equipped with a map from the
same  limit group, $Ipr_i$. 

With each uniformization limit group from the chosen subsequence there is an
associated map from the rigid or solid limit group $Ipr_i$ into that uniformization limit group. By the
construction of the uniformization limit groups, $Unif^j$, by passing to a further subsequence
we may assume that the map from $Ipr_i$ extends to a map from a fixed uniformization limit group, $Unif^{j_1}$.
By passing  iteratively to further subsequences we obtain maps from fixed uniformization limit groups,
$Unif^{j_1},Unif^{j_2},\ldots$, into the uniformization limit groups from the corresponding
subsequences.

Now, we look at the sequence of images, $<f_j,p,q_j,a>$, of the limit group, $Ipr_i$, in the uniformization
limit groups, $Unif^j$, along the diagonal subsequence that is taken from the chosen subsequences of 
the infinite path.
Each uniformization limit group in the diagonal subsequence is constructed as a limit of homomorphisms into
a coefficient free group, $F_k$. With each uniformization limit group from the diagonal subsequence,
(that we still denote) $Unif^j$, we associate a homomorphism, $h_j: Unif^j \to F_k$, 
that restricts and lifts to a
homomorphism $s_j: Ipr_i \to F_k$. We choose the homomorphism $h_j$, so that it approximates the distances in
the limit action of $Unif^j$ on the limit $R^{n_j}$-tree, for larger and larger (finite) subsets of elements in
$Ipr_i$. 

To analyze the sequence of homomorphisms $\{s_j\}$, and obtain a contradiction to the existence of an
infinite path, we need the following theorem (theorem 1.5 in [Se3]), that gives a form of strong
accessibility for limit groups.

\proclaim{Theorem 3.10 ([Se3],1.3)} Let $G$ be a f.g.\ group, and let:
 $\{u_n \, | \, u_n : G \to F_k \}$ be a sequence of homomorphisms. Then there 
exist some
integer $m \geq 1$, and a subsequence
of the given sequence of homomorphisms, that converges into a faithful action of some
 limit quotient $L$
of $G$ on some $R$-tree. 

Furthermore, the elements of $L$ (hence, of $G$) have precisely $m$ rates of growth. i.e., we can divide the elements of $L$ into $m$ disjoint sets,
where up to conjugation the ratios of growths of two elements in the same set converges to some non-zero limit.
\endproclaim

By theorem 3.10, from the sequence of homomorphisms, $s_j:Ipr_i \to F_k$, it is possible to extract a subsequence 
(that we still denote $\{s_j\}$) that converges into a faithful action of a limit quotient $L$ of $Ipr_i$ on some $R$-tree, 
where the elements of $L$ can be divided into $m$ subsets with asymptotically comparable growth,
for some integer $m \geq 1$. 

By construction, the homomorphism from $Ipr_i$ into $Unif^j$ extends to homomorphisms from the
uniformization limit groups, $Unif^{j_1},Unif^{j_2},\ldots$ into $F_k$. Hence, the limit action of the image of
$Ipr_i$ in $Unif^j$ on the associated $R^{n_j}$-tree  contains at least $j$ levels of infinitesimals. Since the
homomorphisms $s_j:Ipr_i \to F_k$ were chosen to approximate these limit actions on larger and larger sets of elements
of $Ipr_i$, it can not be that the limit action that is obtained from the sequence of homomorphisms
$s_j: Ipr_i \to F_k$  contains only a finite (in fact, bounded) sequence of infinitesimals. Therefore, we
obtained a contradiction to the existence of an infinite path, and the procedure for the construction of
uniformization limit groups terminate after finitely many steps.

\line{\hss$\qed$}

Theorem 3.9 asserts that the iterative procedure for the construction of the uniformization limit groups,
$Unif^i$, terminates. Since the iterative procedure produces finitely many uniformization limit groups at
each step, until its termination it constructs finitely many uniformization limit groups, $Unif^i$, that we
denote $Unif_1,\ldots,Unif_v$ (we omit the notation for the step it was produced, since this will
not be important in the sequel). 
With each uniformization limit group, 
$Unif_i$, there is an associated map from one of the rigid and solid limit groups
$Ipr_1,\ldots,Ipr_w$ into $Unif_i$, that we denote,
$<f_i,p, q_i,a>$. Note that by our assumptions all
the terminal limit groups of the uniformization limit groups, $Unif_1,\ldots,Unif_v$, 
are rigid (and not solid). If the images of the subgroups $<p>$ and $<q>$ in 
$Unif_i$ are both non-trivial, then the subgroup $<f_i,p, q_i,a>$ intersects some conjugates of the 
distinguished vertex group in $Unif_i$ non-trivially.  
We set $H_i^1,\ldots,H_i^{e_i}$ to be the conjugacy classes of these subgroups of intersection.

\vglue 1pc
\proclaim{Theorem 3.11} Suppose that all the uniformization limit groups, $Unif_1,\ldots,Unif_v$, and the 
duo limit groups 
that were used for their construction, terminate in rigid limit groups, and the graded completions
that are associated with these groups contain no abelian vertex groups.

Then for every equivalence class of $E(p,q)$, which is not one of the finitely
many equivalence classes that are excluded in theorem 3.1, there exists a uniformization
 limit group, $Unif_i$, from
the finite collection, $Unif_1,\ldots,Unif_v$, so that there exists a (positive) bounded number of 
conjugacy classes of specializations of the subgroups, $H_i^1,\ldots,H_i^{e_i}$,  
for which (cf. theorem 3.8):
\roster
\item"{(1)}" there exist specializations in the given conjugacy classes of specializations of
the subgroups $H_i^1,\ldots,H_i^{e_i}$ that  can be extended to  rigid specializations of the distinguished
(terminal) rigid vertex group in the uniformization limit group $Unif_i$,
that can be further extended to test sequences of specializations that 
restrict to specializations of elements in the given
equivalence class of $E(p,q)$.

The test sequences of specializations that extend the corresponding rigid specializations of
the distinguished vertex group in the uniformization limit group $Unif_i$, 
restrict to valid proofs that the sequence
of couples $\{(p_n,q_n)\},\ldots,\{(p_n,q^i_n)\}$ are in the given equivalence class of $E(p,q)$, 
and to distinct sequence of couples 
$\{(p_n,q_n)\},\ldots,\{(p_n,q^i_n)\}$. 
Furthermore, with the uniformization limit group, $Unif_i$, there are finitely many associated maps from
the subgroups $Ipr_1,\ldots,Ipr_w$. With each such map, there are finitely many associated subgroups $H_{i,j}$
(that were associated with $Unif_i$ along the iterative procedure that constructs it). 
Then the test sequences that extend the rigid specializations of the distinguished vertex group in $Unif_i$,
restricts to non pairwise conjugate specializations of the subgroups $H_{i,j}$ (for each level $j$), except  
for the bottom level subgroups,   
$H_i^1,\ldots,H_i^{e_i}$.

\item"{(2)}" the boundedly many conjugacy classes of specializations of the subgroups,
$H_i^1,\ldots,H_i^{e_i}$, are the only conjugacy classes of specializations of these subgroups that satisfy
part (1) for the given equivalence class of $E(p,q)$. 
\endroster
Note that the bound on the number of  conjugacy classes of specializations of the
subgroups, 
$H_i^1,\ldots,H_i^{e_i}$,
is uniform and it does not depend on the
given equivalence class.
\endproclaim

\nfp By the construction of the first level uniformization limit groups, $Unif^1_j$, for each equivalence class
that is not one of the finitely many equivalence classes that are excluded in theorem 3.1, there exists 
a uniformization limit group, $Unif^1_j$, with (conjugacy classes of) specializations of the of the
subgroups, $H_1,\ldots,H_{e_1}$, that satisfy part (1). If for a given equivalence class there are
infinitely many such conjugacy classes, we pass the uniformization limit groups that were constructed in
the second level. By continuing iteratively, and by the termination of the procedure for the construction
of uniformization limit groups (theorem 3.9), for each given equivalence class (which was not excluded by theorem
3.1), we must reach a level in which there is a uniformization limit group with (conjugacy classes of)
specializations of the subgroups, $H^i_1,\ldots,H^i_{e_i}$, that satisfy the conclusion of the theorem.  

\line{\hss$\qed$}

Theorem 3.11 proves that for any equivalence class of $E(p,q)$ (except the finitely many equivalence
classes that are singled out in theorem 3.1), there exists some uniformization limit group, $Unif_i$, 
which is one of the
constructed uniformization limit groups, $Unif_1,\ldots,Unif_v$, for which the subgroups, 
$H_i^1,\ldots,H_i^{e_i}$, which are the conjugacy classes of intersecting subgroups between the subgroup,
$<f_i,p, q_i,a>$, and conjugates of the distinguished vertex group in the  uniformization limit group, $Unif_i$,
admit only boundedly many conjugacy classes of specializations (that can be extended to
test sequences of $Unif_i$ that satisfy part (1) in theorem 3.11).

Therefore, these boundedly 
many conjugacy classes of specializations of the subgroups, 
$H_i^1,\ldots,H_i^{e_i}$,
already enable us to construct a (definable) function from the
collection of equivalence classes of $E(p,q)$ into a power set of the coefficient group $F_k$, so that the function
maps each equivalence class of $E(p,q)$ into a (globally) bounded set. However, it is not guaranteed that the
class function that one can define in that way, separates between different classes of $E(p,q)$. 

The uniform bounds on the conjugacy classes of specializations of the subgroups, 
$H_i^1,\ldots,H_i^{e_i}$, does not yet give us the desired class function that we can associate with $E(p,q)$,
i.e., a class function with "bounded" image for each equivalence class. It does give us a $separation$
$of$ $variables$ that can be used as a step towards obtaining a desired class function. In order to obtain
this separation of variables we need to look once again at the  decomposition that we denote $\Lambda_i$,
which is the decomposition that is inherited by the subgroup $<f_i,p, q_i>$ 
from the uniformization limit group, $Unif_i$, from which the uniformization limit group, $Unif_i$,
was constructed. 

\vglue 1pc
\proclaim{Lemma 3.12} With the notation of theorem 3.11, $\Lambda_i$, the graph of groups decomposition that is 
inherited by the subgroup, $<f_i,p,q_i,a>$, from the (terminal) uniformization limit group, $Unif_i$, is either:
\roster
\item"{(1)}" $\Lambda_i$ is a trivial graph, i.e., a graph that contains a single vertex. 
In that case either the subgroup $<p>$
or the subgroup $<q_i>$ is contained in the distinguished vertex in $Unif_i$, and in particular,
it admits boundedly many values.

\item"{(2)}"  $\Lambda_i$ has at least two vertices,
 and (finitely many) edges between them. The subgroup $<p>$ is contained 
in one vertex group in $\Lambda_i$, and
the subgroup $<q_i>$ is contained in a different vertex group in $\Lambda_i$ and there exists an edge between 
the vertices that are stabilized by $<p>$ and $<q>$ in $\Lambda_i$. 

Only the edge group between the vertices that are stabilized between $<p>$ and $<q>$ can have a trivial stabilizer, and in that case
there must be more edges in $\Lambda_i$.
The subgroups,
$H_i^1,\ldots,H_i^{e_i}$, contain conjugates of all the edge groups in $\Lambda_i$ (except possibly the trivial
edge group.
\endroster
\endproclaim

\nfp Both parts follow because $\Lambda_i$ is the graph of groups that is induced by $<f_i,p,q>$ from the decomposition of the
uniformization limit group $Unif_i$. The edge groups in $\Lambda_i$ are by construction 
conjugates of $H_i^1,\ldots,H_i^{e_i}$. 

If $\Lambda_i$ contains more than a single trivial edge group, or a trivial edge group that
does not connect between the vertices that are stabilized by $<p>$ and $<q>$, then either $<f_i,p,q>$ admits a free decomposition in which
$<p>$ and $<q>$ are contained in the same factor, or a free product in which $<p>$ is contained in one factor, and $<q>$ in the second factor.
Both are not possible by the restriction of the homomorphisms from the groups $Ipr_1,\ldots,Ipr_w$ from which the uniformization limit groups were constructed.

\line{\hss$\qed$}

If we combine lemma 3.12 with theorem 3.11, for each equivalence class of $E(p,q)$ except the finitely
many equivalence classes that are singled out in theorem 3.1, there exists a uniformization limit group, 
$Unif_i$,
from the finite collection of the constructed uniformization limit groups, $Unif_1,\ldots,Unif_v$, 
for which:
\roster
\item"{(1)}" with the uniformization limit group, $Unif_i$, there is an associated map from 
one of the rigid and solid limit groups, 
$Ipr_1,\ldots,Ipr_w$, into $Unif_i$, with image, $<f_i,p,q_i,a>$.

\item"{(2)}" either one of the subgroups, $<p>$ or $<q_i>$ admits boundedly many values up to conjugacy, or
the subgroup $<f_i,p, q_i,a>$ inherits a graph of groups decomposition, $\Lambda_i$,  from the
presentation of $Unif_i$ as an amalgamated product. $\Lambda_i$ contains at least two vertices, where $<p>$ is contained
in one vertex group and $<q_i>$ in another second vertex group.

\item"{(3)}" the edge groups in $\Lambda_i$ are conjugates to some of the subgroups, 
$H_i^1,\ldots,H_i^{e_i}$, and these groups admit only boundedly many conjugacy classes of specializations, that
are associated with the given equivalence class, and satisfy the conditions that are presented in
part (1) of theorem 3.11.
\endroster

As we have already indicated, the bounded number of conjugacy classes of specializations of the subgroups,
$H_i^1,\ldots,H_i^{e_i}$, that are associated with a given equivalence class of $E(p,q)$, associated
a bounded set with each equivalence class. However, these bounded sets may not separate between different
equivalence classes. The graph of groups, 
$\Lambda_i$, that is 
inherited by the subgroup, $<f_i,p, q_i,a>$, from the uniformization limit group, $Unif_i$, 
can be viewed as a $separation$ $of$ $variables$
($\Lambda_i$ separates between the subgroups $<p>$ and $<q_i>$). 
This separation of variables is the goal of this section, and the key for associating parameters with equivalence
classes of the definable equivalence relation $E(p,q)$ in the next section, 
parameters that admit boundedly many values for each class,
and these values separate between the different classes.

\medskip
For presentation purposes, in all the constructions that were involved in obtaining the uniformization limit groups, 
$Unif_1,\ldots,Unif_v$, we assumed that the terminal limit group in all the associated duo limit groups are
rigid (and not solid), and that there is no abelian vertex group in all the abelian decompositions that are
associated with the various levels of the completions that are part of the constructed duo limit groups.
Before we continue to the next section, and use the separation of variables that we obtained to associate parameters
with equivalence classes, we generalize the constructions that we used, to omit these technical assumptions.


The construction of the rigid and solid limit groups, $Ipr_1,\ldots,Ipr_w$, and their properties that are
listed in theorem 3.1, do not depend on the structure of the graded completions, $G_1,\ldots,G_t$, that
form the Diophantine envelope of the given equivalence relation, $E(p,q)$. Hence, we can use them as in the special case (rigid
terminal groups of envelopes and uniformization limit groups)
that were analyzed before.

From the graded completions, $G_1,\ldots,G_t$, their test sequences and sequences of (rigid and almost
shortest strictly solid) homomorphisms
from the rigid and solid limit groups, $Ipr_1,\ldots,Ipr_w$, we constructed the duo limit groups,
$Tduo_1,\ldots,Tduo_m$, that can also serve as the duo envelope of $E(p,q)$.
By construction, one of the graded completions that is associated with each of
the duo limit groups, $Tduo_1,\ldots,Tduo_m$, 
is a closure of the graded completion, $G_j$, from which it was constructed.

With
each of the duo limit groups, $Tduo_1,\ldots,Tduo_m$ there is an associated subgroup, $<f,p,q>$, that denotes
the image of an associated map from one of the limit groups, $Ipr_1,\ldots,Ipr_w$, into it. Proposition 3.3
is stated for general duo limit groups $Tduo_i$. In order to generalize proposition 3.4 
to general duo limit groups, $Tduo_1,\ldots,Tduo_m$, we need the following observations.

Let $Tduo_i$ be one of the Duo limit groups, 
$Tduo_1,\ldots,Tduo_m$, and 
let the subgroup $<f,p,q,a>$, be the image in $Tduo_i$ of one of the rigid or solid limit groups, 
$Ipr_1,\ldots,Ipr_w$. Suppose that both subgroups, $<p>$ and $<q>$, are non-trivial in $Tduo_i$. 
$Tduo_i$ being a duo limit group, admits an amalgamated product decomposition:
$$Tduo_i=<d_1^i,p>*_{<d_0^i,e_1^i>} \, <d_0^i,e_1^i,e_2^i>*_{<d_0^i,e_2^i>} \, <d_2^i,q>.$$

Furthermore, the distinguished vertex group, $<d_0^i,e_1^i,e_2^i>$, admits a graph of groups decomposition
that we denote, $\Gamma^i_D$,
that is obtained from the graph of groups, $\Gamma^i_{<d_0>}$, that is associated with the terminal rigid or solid limit
group, $<d_0^i>$, so that to each rigid vertex group in $\Gamma^i_{<d_0>}$ one further connects
several (possibly none) free abelian vertex groups, which are subgroups  of  the subgroup $<e_1^i,e_2^i>$, 
along some
free abelian edge groups. 

As in the rigid case, with a duo limit group, $Tduo_i$, we also associate a decomposition (graph of groups), $\Delta_i$.
To construct the graph of groups $\Delta_i$, we start with the graph of groups $\Gamma_i^D$. From $\Gamma^i_D$ we take out all the QH vertex groups.
We call the fundamental group of each connected component $R^i_t$. Note that the graph of groups that is associated with each of the groups
$R^i_t$ has (possibly none) abelian edge groups,
and no $QH$ vertex groups.

In $R^i_t$ there are two canonical subgroups, $RP^i_t$ and $RQ^i_t$. $RP^i_t$ is a fundamental group of a connected component of a graph of groups that
is obtained from $\Gamma^i_{<d_0>}$ by taking out all the QH vertex groups from it, and adding abelian groups that are subgroups of $<e_1>$ and not
$<e_1,e_2>$. $RQ^i_t$ is similar to $RP^i_t$, but in this case we add to $\Gamma^i_{d_0}>$ the abelian groups $<e_2>$.

Starting with $\Gamma^i_D$ we collapse the subgraphs of groups that are associated with the subgroups $R^i_t$. i.e., we are left with a graph
of groups with vertices stabilized by the groups $R^i_t$ and the QH vertex groups in $\Gamma^i_D$. We start by constructing a graph of groups
 $\hat \Delta_i$ by adding to this (collapsed)
graph of groups two vertex groups and several edge groups. 

One of the vertex groups that we add to $\hat \Delta_i$ is stabilized by the completion $G_P=<d_1^i,p>$ in $Tduo_i$. 
The other vertex group is stabilized by the completion $G_Q=<d_2^i,q>$ in $Tduo_i$. The vertex that is stabilized by $G_P$ is connected to the vertices that are 
stabilized by each of the vertex groups $R^i_t$ by an edge that is stabilized by $RP^i_t$.  
The vertex that is stabilized by $G_Q$ is connected to the vertices that are 
stabilized by each of the vertex groups $R^i_t$ by an edge that is stabilized by $RQ^i_t$.  

To construct $\Delta_i$ from $\hat \Delta_i$, we  fold each of the edge groups that are stabilized by the edge groups $RP^i_t$ and $RQ^i_t$, and
replace them with the vertex groups $R^i_t$ that contain them.  
The obtained graph of groups that we denote $\Delta_i$, contains two vertex groups, that are stabilized by $<G_P,e_2>$ and $<G_Q,e_1>$, and possibly
some QH vertex groups from $\Gamma^i_D$. It contains cyclic edge groups that are connected to the QH vertex groups, 
and several edge groups that are stabilized by the edge groups, $R^i_t$. 
By construction, $<p>$ is a subgroup of one of the vertex groups in $\Delta_i$, and $<q>$ is a subgroup of the other vertex group in $\Delta_i$.

Recall that the subgroup $<f,p,q,a>$ is the image of one of the limit groups, $Ipr_1,\ldots,Ipr_w$, in the duo limit group $Tduo$, 
and it does not factor through a free product in which $<p>$ is contained in one factor
and $<q>$ is contained in  another factor.
Hence, if both subgroups $<p>$ and $<q>$ are non-trivial in $Tduo_i$, then $<f,p,q,a>$ 
intersects non-trivially some conjugates of the edge groups in the graph of groups $\Delta_i$. Furthermore, if it intersects a QH vertex group in
$\Delta_i$, either it intersects it in a finite index subgroup, or  the intersection is of infinite index, and it has to be
a free product of conjugates of (cyclic) subgroups
that are contained in  boundary subgroups of the QH vertex group, i.e., subgroups that are contained in conjugates of the edge groups in $\Delta_i$.

Let  $V_i^1,\ldots,V_i^f$ be  the conjugacy classes of intersections between the subgroup, $<f,p,q>$, and  the
edge groups in $\Delta_i$. 

\vglue 1pc
\proclaim{Lemma 3.13} Let $Tduo_i$ be one of the duo limit groups, 
$Tduo_1,\ldots,Tduo_m$, and suppose that both subgroups, $<p>$ and $<q>$, are non-trivial in
$Tduo_i$. Then for each value of the parameters, $q_0$,
%
%
there are at most boundedly many values
of $V^j_i$ that are associated with rigid or strictly solid values (homomorphisms) of the rigid or solid limit group $<d_0>$, 
up to conjugacy and the action of the modular groups that are associated with the edge groups in $R^i_t$ (that are subgroups of
the modular groups that is associated with the  graded abelian decomposition, $\Gamma^i_D$).
\endproclaim 

\nfp The lemma follows from the existence of a global bound on the number of rigid and strictly solid families of a rigid or a solid limit group
(Theorem 2.9 in [Se3]).

\line{\hss$\qed$}

Lemma 3.13 claims that with each value of $q_0$ there are at most boundedly many classes of values of
the groups $V^j_i$ (up to conjugacy and the modular groups of the edge groups $R^i_t$). 
However, as in the rigid case, it may be that for some of the equivalence classes there are infinitely many possible
classes of values of the groups $V^j_i$. i.e., for different values of the parameter $q_0$ we get different families of values. 
 
In this case of infinitely many families of values of the groups $V^j_i$ that are associated with an equivalence class we continue 
as we did in the rigid case. We  construct uniformization limit groups, and obtain further splittings of the groups $V^j_i$. With a
uniformization limit group that terminates with a solid limit group or has abelian vertex groups we associate a similar  decomposition
(graph of groups) $\Delta_i$ (that we still denote in the same way).

One of the groups $Ipr_1,\ldots,Ipr_w$ is mapped into each of the uniformization limit groups, and its image is (still) denoted, $<f,p,q,a>$.
The intersection subgroups of $<f,p,q,a>$ with the edge group in $\Delta_i$, that we still denote $V^j_i$,
 satisfy the conclusions of
lemma 3.13. As in the rigid case, when we use the uniformization limit groups, the conclusions of lemma 3.13 is valid globally, i.e., the boundedness
of the number of families of values of the intersection subgroups is uniform, and these families do not depend on the  parameters $q_0$ from the equivalence
classes, as it may happen in the case of
the duo
limit groups $Tduo$.

Lemma 3.13 enables one to generalize lemma 3.12, and its conclusions w.r.t.\ the graph of groups $\Lambda$ that the subgroup $<f,p,q>$ inherits from 
the decomposition of the terminal uniformization limit group, to the case in which the terminal limit group of the uniformization limit group may
be solid and the uniformization limit group may contain abelian vertex groups along its levels.




\vglue 1.5pc
\centerline{\bf{\S4. Equivalence Relations and their Parameters}}
\medskip
In the first section
of this paper
we have constructed the Diophantine envelope of a definable set (theorem 1.5), and then 
used it to construct the Duo envelope of a definable set (theorem 1.7). 

Recall that by its definition (see definition 1.1), a
Duo limit group $Duo$ admits an amalgamated product: 
$Duo=<d_1,p>*_{<d_0,e_1>} \, <d_0,e_1,e_2>*_{<d_0,e_2>} \, <d_2,q>$
where $<e_1>$ and $<e_2>$ are free abelian groups with pegs in $<d_0>$, i.e., free abelian groups that 
commute with non-trivial elements in $<d_0>$.
A specialization of the parameters $<d_0>$ of a Duo limit group gives us
a Duo family of it.

To analyze definable 
equivalence relations over a free (or a hyperbolic) group, our strategy is to further study
the parameters ($<d_0>$) that are associated with the Duo families that are associated with the Duo
limit groups that form the Duo envelope of a  definable equivalence relation. 

In the previous  section we modified and  analyzed   
the construction of the Duo envelopes that were presented in theorem 1.7, 
in the special case of a definable equivalence relation. We further
carefully studied the set of values of the parameters that are associated with the duo families that are
associated with each equivalence class.
This careful study, that uses what we called $uniformization$ limit groups that we associated with the
Duo envelope, enabled one to associate  
a "bounded" set of families of values of certain (edge) subgroups of the parameters that are associated 
with the Duo families of the Duo envelope, for each
equivalence class of a definable 
equivalence relation. 

The bounds that we achieved on the number of families of values of some subgroups 
that are associated with each equivalence class,
allowed us to obtain what we view as 
"separation of variables". This means that with the original subgroups of parameters, $<p>$ and $<q>$, we associate
a bigger subgroup, for which there exists a graph of groups decomposition, where $<p>$ is contained in one vertex group,
$<q>$ is contained in a second vertex group, and  the number of families of values of the edge groups in the graph of groups
is uniformly bounded for the classes in the given equivalence class $E(p,q)$.

However, the parameters that we associated with each equivalence class, i.e., the families of values of the edge groups in these graphs of groups, 
do not separate between equivalence classes in general. To obtain 
subgroups of parameters that have the same types of bounds as the ones that were constructed in the previous section, 
that do
separate between classes, we present a new iterative procedure that uses both the sieve procedure [Se6]
(that was used for quantifier elimination) together with the procedure for separation of variables that
was presented in the previous section (i.e., the iterative construction of uniformization limit groups). 

The combined procedure is a (new) sieve procedure that preserves the
separation of variables along its various steps, and its termination (that follows from the termination of
the sieve procedure and the procedure for the separation of variables), produces the desired subgroups
of parameters, that do separate between equivalence classes and admit boundedly many families of values 
for each equivalence class (where the bound on the number of families of values does not depend on the
specific equivalence class).

\medskip
Let $F_k=<a_1,\ldots,a_k>$ be a non-abelian free group, 
and let $E(p,q)$ be a definable equivalence relation over $F_k$. 
With the definable equivalence relation, $E(p,q)$, being a definable set, one associates 
using theorems 1.5 and 1.7, 
a Diophantine and a Duo envelopes. 
Let $G_1,\ldots,G_t$ be the Diophantine envelope of the given definable equivalence relation $E(p,q)$,
and let $Duo_1,\ldots,Duo_r$, be its Duo envelope. 



Theorem 3.1 associates with the given definable equivalence relation, $E(p,q)$, finitely many rigid and solid limit
 groups, $Ipr_1,\ldots,Ipr_w$, so that apart from finitely many equivalence classes, for each couple,
$(p,q) \in E(p,q)$, there exists a rigid or a strictly solid family of 
homomorphisms from at least one of the limit groups, $Ipr_1,\ldots,Ipr_w$,
to the coefficient group $F_k$, so that  the rigid homomorphisms or the strictly solid homomorphisms  
from the given strictly solid family
do not factor through a free product $A*B$ in which $<p> < A$ and $<q> < B$, and each of these homomorphisms
restricts to
a valid proof that $(p,q) \in E(p,q)$.

\medskip
To obtain separation of variables in the previous section,
we started with the Diophantine envelope of the given definable equivalence relation,
$G_1,\ldots,G_t$.  
With each graded completion $G_j$, $1 \leq j \leq t$,  we associated a finite collection  of duo
limit groups. First, we collected all the test sequences of the completion $G_j$, that can be extended
to rigid or strictly solid specializations of one of the rigid or solid limit groups, $Ipr_1,\ldots,Ipr_w$,
that restrict to valid proofs that the corresponding couples, $\{(p_n,q_n)\}$, are in the equivalence relation
$E(p,q)$. We further required that
these test sequences of specializations can not be factored through a free product in which $<p>$ is contained in 
one factor, and $<q>$ is
contained in the second factor. This collection of test sequences of the graded completions,
$G_1,\ldots,G_t$ that can be extended to specializations
of $Ipr_1,\ldots,Ipr_w$, can be collected in finitely many duo limit groups 
(using the techniques that were used for
collecting formal solutions in [Se2]). Then we used the sieve procedure [Se6]
to construct finitely many duo limit groups, $Tduo_1,\ldots,Tduo_m$, that still collect all these extended
test sequences of the graded completions, $G_1,\ldots,G_t$, and for which there exist
generic points (i.e., duo test sequences) that restrict to valid proofs that the restricted couples, 
$\{(p_n,q_n)\}$, are in the given equivalence relation $E(p,q)$.  

The collection of duo limit groups, $Tduo_1,\ldots,Tduo_m$, is the starting point for the iterative 
procedure for separation of variables. With them we associated  a collection of uniformization 
limit groups, $Unif_1,\ldots,Unif_v$. 
For each equivalence 
class of $E(p,q)$, apart from the finitely many equivalence classes that are singled out in theorem 3.1,
there exists at least one uniformization limit group, $Unif_i$, that enabled us to associate boundedly many families 
of values (of edge groups) with the equivalence class.

Recall that the image of  one of the subgroups, $Ipr_1,\ldots,Ipr_w$, that is associated with the uniformization
limit group, $Unif_i$,
that we denoted, $<f_i,p,q_i>$, inherits a graph of groups decomposition
$\Lambda_i$ from the graph of groups decomposition $\Delta_i$ of the ambient uniformization limit group,
$Unif_i$. For each equivalence class of $E(p,q)$ that is associated with
$Unif_i$, an edge group in $\Lambda_i$ admits boundedly many classes of values.
Furthermore, the subgroups
$<p>$ and $<q>$ are  both contained in vertex groups in $\Lambda_i$.  

The edge groups in the graphs of groups, $\Lambda_i$, enable one to associate parameters with each equivalence 
class of $E(p,q)$, where these parameters admit only boundedly many families of values for each equivalence class.  
The families of values of the edge groups are not guaranteed to separate between equivalence classes in general. 

To use the graphs of groups $\Lambda_i$ to
obtain parameters that do separate between equivalence classes, we use the graphs of groups $\Lambda_i$ as
a first step in an iterative procedure that combines the procedure for separation of variables (that was
presented in the previous section), with the sieve procedure for quantifier elimination that was presented in [Se6].

\medskip
For presentation purposes, we continue with what we did in the previous section, and start by  presenting the combined procedure
assuming that 
the graded closures that are associated with all the duo limit groups that were used in the construction of 
the uniformization limit groups, $Unif_i$, and the graphs of groups $\Lambda_i$, 
terminate in rigid limit groups and
do not contain abelian vertex groups in any of their levels. 
Later on we modify the procedure to omit these assumptions.

We continue with the uniformization limit groups, $Unif_i$, in parallel. Hence, for brevity, we denote the
uniformization limit group that we continue with, $Unif$. By construction, 
with each such uniformization limit group,
$Unif$, there is an associated subgroup, $<f,p,q>$, (which is the image of one of the rigid and solid limit groups,
$Ipr_1,\ldots,Ipr_w$), and graph of groups decomposition, $\Lambda$. Recall that under the assumption that
there are no abelian vertex groups in any of the levels of the graded closures that are associated with
the uniformization limit group, $Unif$, and that the uniformization limit group terminates with a rigid limit group,
the  edge group in $\Lambda$ 
that connects the vertex  that
is stabilized by $<p>$ with the vertex that is stabilized by $<q>$ may be trivial (only in case there are more edges that
connect these two vertices),
and that for each equivalence class
of $E(p,q)$ that is associated with $Unif$, the associated values of the the other edge groups in $\Lambda$ belong 
to boundedly
many conjugacy classes.

If both subgroups, $<p>$ and $<q>$, are contained in the same vertex group in $\Lambda$, then for each
equivalence class that is associated with $Unif$, there are only boundedly many associated values of
either the subgroup $<p>$ or the subgroup $<q>$, and these values that belong to the equivalence class
obviously determine the class. Hence, we can assume that the graph of groups $\Lambda$ contains at least
two vertex groups, and that $<p>$ is contained in one vertex group, and $<q>$ is contained in another
vertex group in $\Lambda$.

We start by associating with $\Lambda$ a graph of groups decomposition of a bigger group. 
In section 12 of [Se1] we presented the multi-graded Makanin-Razborov diagram.
Recall that this multi-graded diagram encodes all the homomorphisms of a given limit group into a free group,
if the specialization of a certain subgroup of the limit group is fixed, and the specializations of finitely many other
subgroups is fixed up to conjugacy. 


With each vertex group in $\Lambda$ 
$\Theta$, which does not contain $<p>$ 
we associate its multi-graded Makanin-Razborov diagram w.r.t.\ the edge groups that are
connected to the vertex group. Note that for each equivalence class in $E(p,q)$ there are only boundedly many values up to conjugacy that each edge may obtain.
Each such multi-graded diagram contains finitely many resolutions. We go over all the tuples of resolutions of these vertex groups, one from
each MR diagram of a vertex group in $\Lambda$. 

Given each such tuple of resolutions  we get a multi-graded resolution of the fundamental group of $\Lambda$, or a of a quotient of it. 
We go over this (finite) collection of 
multi-graded resolutions in parallel. We denote such a multi-graded resolution $MGRes$.



Given a multi-graded resolution $MGRes$, we associate with it the standard 
(multi-) graded closures following definition 1.25-1.30 in [Se5], that demonstrate how generic points in the (completion of the) resolution
(i.e., multi-graded test sequences of it) fail to be valid proof statements.

If there are test sequences of the multi-graded resolution $GRes$ that are not covered by the collection of the constructed multi-graded closures
then these test sequences are (possibly apart from a finite prefix) in the definable relation $E(p,q)$. 

Therefore, for these
collections of test sequences, the conjugacy classes of the values of the terminal levels of multi-graded resolutions, that are associated with 
all the edge groups in $\Lambda$, 
are not only associated with their equivalence classes, but they also separate them (i.e., define them). 
Since given the structure of 
the multi-graded resolution, the values that factor through the resolution are completely determined by the conjugacy classes of the terminal vertex groups,
and the finitely many cosets that are determined by the constructed multi-graded closures.
Hence, these conjugacy classes of the terminal levels of the multi-graded resolutions can serve as (definable) parameters of  their corresponding equivalence
classes.  

Some of the multi-graded closures of the multi-graded resolution $MGRes$ collect test sequences for which the restrictions of the test sequences
to the proof statement, $<f,p,q>$, are not valid proof statements because there are extra rigid or solid values that are not part of the proof
statement. We call these graded closures, Extra PS.

We continue to the next steps only with the Extra PS multi-graded closures.
%
If there are no Extra PS multi-graded closures, or if each extra PS multi-graded closure has a covering closure of Generic Collapse extra multi-graded closures
(Generic Collapse extra multi-graded closures are constructed from test sequences of the Extra PS multi-graded closures for which the restrictions
of the test sequence to the extra rigid and strictly solid values are not rigid or not strictly solid or they are in fact in the same families of 
rigid and strictly solid values that are part of the proof statement),
we are done. i.e., we reached a terminal point of the iterative procedure.

Otherwise, tuples $(p_0,q_0) \in E(p,q)$ that extend to specializations that factor through an Extra PS multi-graded closure, 
but for which test sequences in the same fiber restrict to proof statement that are not valid,
have to satisfy  one of finitely many Diophantine conditions that indicate that the extra
rigid or strictly solid specializations, that are part of the extra PS multi-graded closure,  are either flexible or are  the same as some rigid serialization or
in the same family of 
a strictly solid specialization that
are indicated by $<f,p,q>$ (the original proof statement). We call each of these finitely many Diophantine conditions, $collapse$ $forms$.

We continue with the (finitely many)  extra PS multi-graded closures. With each extra PS graded closure we continue with  each of its 
(finitely many) collapse forms in parallel.
Each collapse form is a Diophantine condition, that can be demonstrated by adding elements  that demonstrate or prove
the validity of the Diophantine condition
(see section 1 and 3 of [Se5] for more detailed
explanation of these Diophantine conditions, and the way that they are imposed). 

With the collection of specializations that factor through an extra PS multi-graded closure, $ExtraPS$, 
and extend to specializations of elements that demonstrate a fixed collapse form, we associate 
finitely many  multi-graded  limit groups.
Each such
multi-graded  limit group contains a subgroup that is generated by the subgroup $<f,p,q>$, that demonstrate the proof statement that
values of $<p,q>$ are indeed in the equivalence relation, $E(p,q)$.
We call these multi-graded limit groups
 $Extra$ $Collapse$ limit group and denote them $ExtraCollapse$.


Each Extra Collapse group contains an image of the multi-graded closure from which it was constructed, and additional elements for the
extra rigid and strictly solid values, and for their collapse form.
By the structure of the multi-graded closure from which $ExtraCollapse$ was constructed, in it the subgroups $<p>$ and $<q>$ were separated. This is
exactly what we achieved in the previous section. The additional elements that were added for the collapse forms may prevent separation of variable
in the Extra Collapse limit groups. 
Our aim is to get
separation of variables for these extra collapse limit groups, similar to the one that was obtained in the previous section for the subgroup $<f,p,q>$. 
i.e., a procedure to decompose the Extra collapse
limit group into a graph of groups, one that contains $<p>$ and one that contains $<q>$ and possibly finitely many other vertex groups, such that the values of the 
edge groups that are associated with 
each equivalence class in $E(p,q)$ belong to boundedly many conjugacy classes.

\smallskip
As we did in the previous section, to start a separation of variable procedure, we need the groups in question not to factor
through certain free products.  Hence, we collect all the specializations of the Extra Collapse limit groups, $ExtraCollapse$, 
for which:
\roster
\item"{(1)}"  the restriction to the
(image of the) subgroup $<f,p,q>$ form a proof that the couple $(p,q)$ is in $E(p,q)$.

\item"{(2)}"  the restriction to the specialization of the subgroup $<f,p,q>$ does not factor
through a free product of limit groups
in which the subgroup $<p>$ is in one factor and the subgroup $<q>$ is contained
in a second factor, and does not factor through a free product of limit groups 
in which the subgroup $<p,q>$ is contained in a factor.


\item"{(3)}" the ambient specialization (of $ExtraCollapse$) factors through a free product of limit
groups in which
the subgroup $<p,q>$  
is contained in one factor.
\endroster

The  collection of all these specializations (of $ExtraCollapse$) factor through 
finitely many limit groups. By looking at the actions of these limit groups on Bass-Serre trees corresponding
to the free products of limit groups
through which the specializations factor, and apply the shortening procedure for these
actions, 
we can replace these limit groups, by a collection of finitely many (quotient) limit groups,
$AGF_1,\ldots,AGF_a$, so that each of them admits a free product in which the subgroup $<f,p,q>$ is
contained in one factor. 
Therefore, we can replace the collection of these specializations of the Extra Collapse limit groups,
$ExtraCollapse$, by the limit groups that are
associated with the factors that contain the subgroup $<f,p,q>$ 
 in each of the limit groups, 
$AGF_1,\ldots,AGF_a$, and these factors still demonstrate that the Diophantine condition, that is associated
with the corresponding collapse form, does hold.

To analyze the  (non-generic) pairs $(p_0,q_0) \in E(p,q)$, that extend to a  specialization of the subgroup $<f,p,q>$, which
is a valid proof, and this proof extends to a specialization of one of the extra collapse limit groups, $ExtraCollapse$, we do 
 the following.


We continue by looking at the collection of specializations of $ExtraCollapse$ or of $AGF_1,\ldots,AGF_a$, for which the restriction to the subgroup $<f,p,q>$ 
does not factor through a free product in which $<p>$ is contained in one factor and $<q>$ in the other,
and the specialization of $ExtraCollapse$ or $AGF_1,\ldots,AGF_a$ does not factor through a free product in which $<p,q>$ is contained
in a factor.

We continue by applying the procedure for separation of variables, that was presented in the previous section, to the collection of specialization
that factor through $ExtraCollapse$ and do not factor through free products as we indicated. First, we construct with these specializations a graded
Diophantine envelope. i.e., a Diophantine envelope that is graded w.r.t.\ the variables $<q>$. Note that the edge groups in the graph
of groups $\Lambda$, from which the Extra multi-graded closures were constructed, admit only boundedly many values up to conjugacy for each equivalence
class of $E(p,q)$, so these edge groups must belong to the terminal vertex groups in the constructed Diophantine envelope.

Having a Diophantine envelope, we construct a Duo envelope, and uniformization limit groups, precisely as we did in the previous section. i.e., as we did
in constructing the graph of groups $\Lambda$.  
We continue construction uniformization limit groups iteratively, until the intersections of the image of the $ExtraCollapse$ with the terminal levels of the
uniformization limit groups admit only boundedly many values up to conjugacy for every equivalence class in $E(p,q)$.

At this point the image of $ExtraCollapse$ in the final  uniformization limit group
inherits a graph of groups from the graph of groups of the terminal uniformization limit group. The structure of this
graph of groups is similar to the graph of groups that $<f,p,q>$ inherited from the terminal uniformization limit group in the previous section.
It has one vertex group that contains $<p>$, one vertex group that contains $<q>$, and possibly finitely many additional vertex groups. At most
one edge groups in this inherited graph of groups may be trivial, and if there is an edge group with trivial stabilizer it has to be the edge between the vertices that are stabilized by $<p>$ and $<q>$, and in this case there must be additional edges with non-trivial stabilizers. The
edge groups in this inherited graph of groups do all obtain boundedly many values up to conjugacy for each equivalence class of $E(p,q)$. 
We denote this graph of groups $\Lambda_1$ as it was obtained in
a similar way to the construction of the graph of groups $\Lambda$ that was associated with the subgroup $<f,p,q>$. 

The graph of groups $\Lambda_1$ has a vertex group that contains $<p>$. Also, the completions of the other multi-graded resolutions from which the
multi-graded resolution $MGRes$ was constructed, are mapped into the limit group  $ExtCollapse$. We analyze the finite set of multi-graded resolutions 
that are associated with the other vertex groups in $\Lambda_1$ using the second step of the sieve procedure, as it appear in [Se6].

%

We continue iteratively. First we look at all the test sequences of the constructed multi-graded resolutions that restrict to proof statements 
$<f,p,q>$ are not valid proof statements. With these test sequences we associate finitely many graded closures, that enable us to decide
for what values of the edge groups in $\Lambda_1$ (up to conjugacy) there are test sequences that restrict to values $\{(p_n,q_n)$ that are
in $E(p,q)$. These are completely determined by the values of the edge groups in $\Lambda_1$, and the values of these edge groups are definable 
and are class functions, i.e, determined by the class in $E(p,q)$.

We continue with multi-graded fibers that do not have test sequences that restrict to valid proof statements. For these we continue with Extra PS
graded closures and $ExtraCollapse$ limit groups, precisely as we analyzed the multi-graded resolution that is associated with the graph of groups $\Lambda$.
At each step we construct a graph of groups $\Lambda_n$, obtain separation of variables, and apply the general step of the sieve procedure in [Se6]. 

\vglue 1pc
\proclaim{Theorem 4.1} In case the terminating limit groups in all the resolutions that are used in constructing the 
collapse uniformization limit groups are rigid, and there are no abelian vertex groups in any of the abelian decompositions that
are associated with the various levels of the collapse uniformization limit groups that are constructed along the
iterative procedure for the analysis of the parameters of equivalence classes, the iterative procedure, that combines
the sieve procedure with the procedure for separation of variables,
terminates after finitely many steps. 
\endproclaim

\nfp Follows from the termination of the sieve procedure (theorem 22 in [Se6]).

\line{\hss$\qed$}

\medskip
Suppose that the abelian decompositions that are associated with the various levels of the
uniformization limit groups,  that were constructed in the previous section,
do contain abelian vertex groups and their terminal limit groups are either rigid or solid. 

In this general case, we do what we did in the previous section. With each of the uniformization limit groups,
that were constructed in the procedure for separation of variables in the
previous section, there is (also) an associated subgroup, $<f,p,q>$, which is the image of one of
the rigid or solid limit groups, $Ipr_1,\ldots,Ipr_w$, that were constructed in theorem 3.1, and a graph
of groups decomposition, (also denoted) $\Lambda$, that is inherited from the decomposition $\Delta$ of the uniformization limit group,
as we described in the previous section (see lemma 3.13). 

The graph of groups, $\Lambda$, 
gives us a separation of variables in the general case as well. Note that in the general case, for each equivalence class of $E(p,q)$, there is no 
bound on the number of conjugacy classes of values of the edge groups in $\Lambda$ that are associated with the equivalence class, 
but there is a bound on the number of
the families of values of the edge groups up to conjugation and the modular groups of the terminal limit group of the uniformization limit group
(cf. lemma 3.13).

Recall that in lemma 3.13 we denoted the edge groups in $\Lambda$, that are the intersections of the group $<f,p,q>$ with the edge groups in 
$\Delta_i$ by $V^j_i$. If $V^j_i$ is the intersection of $<f,p,q>$ with a rigid vertex group or with an edge group
 in the terminal limit group of a duo  limit
group $Tduo$, or of a uniformization limit group $Unif$, then for each equivalence class $V^i_j$ admits boundedly many values up to conjugacy.
If $V^j_i$ is the intersection of $<f,p,q>$ with an edge group $R^t_i$ in the abelian decomposition $\Delta_i$, then it admits boundedly many values
up to conjugacy and the modular group that is associated with the abelian decomposition of $R^i_t$. 

Given $\Lambda$, we can associate the 
taut multi-graded abelian Makanin-Razborov diagram of the vertex group in $\Lambda$ w.r.t.\ its edge groups, i.e., w.r.t.\ the subgroups $V^j_i$.
We continue precisely as we did in case there were no abelian vertex groups in the constructed uniformization limit
groups. 

We construct a multi-graded resolution $MGRes$ from a collection of multi-graded resolutions of the vertex groups in $\Lambda$, in the same way that we did
it in the rigid case. 
Then we construct multi-graded closures and in particular Extra PS multi-graded closures, and  add elements to the demonstrate
the Diophantine condition that is associated with each possible collapse form, precisely as we did in case the terminal limit groups are rigid.
We further apply  the procedure for separation of variables using
collapse uniformization limit groups. Afterwards we apply the construction of quotient (multi-graded) resolutions in the sieve procedure. 
Altogether we get an iterative procedure that terminates after
finitely many steps precisely as in the case of rigid terminal groups and no abelian vertex groups along the constructed completions.

\medskip
The iterative procedure that constructs collapse uniformization limit groups and quotient multi-graded resolutions, 
enables one to associate parameters with the various equivalence classes of the
given definable equivalence relation, $E(p,q)$. The collection of these objects allows one to associate 
a collection of finitely
many elements in finitely many limit groups, with the equivalence relation $E(p,q)$, and for each
equivalence class in $E(p,q)$, these elements admit 
only (uniformly) boundedly many specializations, up to conjugation  and the modular groups of the edge groups $R^i_t$ in the abelian decompositions $\Delta_i$
that are associated with the collapse uniformization limit groups that are constructed along the iterative procedure.

To obtain a form of weak elimination of imaginaries,
after adding new (basic) sorts for the parameters that we associated with equivalence classes,
we still need to define these new sorts, and in particular prove that the parameters that we associated with equivalence classes
are definable.

\vglue 1pc
\proclaim{Theorem 4.2} Let $F$ be a (non-abelian) free group, and let $E(p,q)$ be a  definable equivalence relation over $F$.
Suppose that all the uniformization limit groups, and all the collapse uniformization limit groups that were constructed along the
procedure that analyzes $E(p,q)$ have terminal rigid limit groups, and no abelian vertex groups in any of their levels.
Then if we add a sort for conjugacy classes of elements, $E(p,q)$ can be weakly eliminated.

Suppose that $p$ and $q$ are $m$-tuples. There exist some integers $s$ and $t$ and a  definable multi-function:
$$f: \, F^m  \ \to \ F^s \times R_1 \times \ldots \times R_t$$ 
where each of the $R_i$ is the new sort for conjugacy classes of elements. 
The image of an element is uniformly bounded, 
and the multi-function
 is a class function, i.e., two elements in an equivalence class of $E(p,q)$ have the same image. The multi-function
$f$ separates between classes, i.e., the images of elements from distinct equivalence classes is distinct.  
\endproclaim

\nfp Let $E(p,q)$ be a  definable equivalence relation that satisfies the assumptions of the theorem. To prove that $E(p,q)$ can be  
weakly eliminated,
we need to construct a definable multi-function $f$ as described in the theorem. First, in theorem 3.1 we associated with
$E(p,q)$ finitely many rigid or solid limit groups, $Ipr_1,\ldots,Ipr_w$, so that for all but finitely many equivalence
classes, and for every pair $(p,q)$ in any of the remaining  equivalence classes,  
there is a rigid or a strictly solid homomorphism from one of these limit groups into the coefficient group $F_k$, that
restricts to a valid proof statement that the pair $(p,q)$  is in $E(p,q)$, and the homomorphism does not factor through a free product of
limit groups in which $<p>$ is contained in one factor, and $<q>$ is contained in a second factor. 

All our further 
constructions (of uniformization limit groups) are based on the existence of such homomorphisms, and hence, the finitely
many equivalence classes that were singled out by theorem 3.1 are excluded. Therefore, to construct the desired 
multi-function 
$f$, we need to show that the finite collection of equivalence classes that were singled out in theorem 3.1, is a definable
collection.

\vglue 1pc
\proclaim{Proposition  4.3} Let $E(p,q)$ be a  definable  equivalence relation in a free group $F$. 
The finite collection of equivalence classes that were singled out in theorem 3.1 is a  definable
collection. 
\endproclaim

\nfp Recall that in order to prove theorem 3.1 we have constructed finitely many ungraded limit groups,
$GFD_1,\ldots,GFD_d$, that admit a free product decomposition in which $<p>$ is contained in one factor,
and $<q>$ is contained in a second factor, $GFD_j=A_j*B_j$, $<p> \, < \, A_j$ and $<q> \, < \, B_j$, $1 \leq j \leq d$.

With each limit group $GFD_j$ we have associated its taut 
Makanin-Razborov diagram. Each (ungraded) resolution $Res$ in the diagram is composed from two resolutions, $Res_1$ of $A_j$ and $Res_2$ of $B_j$.
Equivalence classes of $E(p,q)$ that we excluded in theorem 3.1, were classes for which there exists a fixed value
$q_0$ of the variables $q$, and a test sequence
of  the resolution $Res_1$, such that the test sequence restricts to values $\{(p_n,q_0)\}$ that are in the equivalence relation $E(p,q)$.  
By lemma 3.2  such test sequences restrict to values $\{p_n\}$ that belong to finitely many equivalence classes of $E(p,q)$. 

The finite collection of equivalence classes of $E(p,q)$ 
that were excluded in theorem 3.1 are the classes that were excluded for each of the ungraded limit groups $GFD_j$ , and classes that were 
excluded for each of ungraded limit groups that were constructed along the finite iterative procedure that was used to prove theorem 3.1, where each of
these (finitely many) ungraded limit groups decomposes to a free product similar to that of the groups $GFD_j$, and the classes that were excluded for
each of these limit groups
have precisely the same properties as those that were excluded for the groups $GFD_j$.

To prove lemma 3.2 we constructed finitely many graded closures of the resolutions $Res_1$ (graded w.r.t\ the subgroup $<q>$), 
that are associated with the limit groups $GFD_j$. The classes that are excluded by lemma 3.2 are precisely the classes for which the graded closures that
collect all the obstructions for the proof statements that appear in $GFD_j$ to be valid proof statements, do not form a covering closure of the completion of
the resolution $Res_1$. 

The values of $q$ for which a collection of graded closures is not a covering closure is a definable set, in fact it is in the Boolean algebra of
AE sets. Hence, the finitely many equivalence classes that are excluded in theorem 3.1 are definable.



 
\line{\hss$\qed$}

Proposition 4.3 shows that the finite collection of equivalence classes that are excluded
in theorem 3.1 is definable. Therefore, to prove theorem 4.2, we need to construct a (definable)
function
with the properties that are listed in the statement of the theorem, that is defined
on the union of all the other
equivalence classes of $E(p,q)$.

Let $p_0$ be a specialization of the  variables $p$, that does not
belong to one of the finitely many equivalence classes that are excluded in theorem 3.1.
Then for  each pair $(p,q) \in E(p,q)$, that are in the same equivalence class as $p_0$, 
there exists a 
rigid or a strictly solid homomorphism $h$ from one of the rigid or strictly solid 
limit group, $Ipr_1,\ldots,Ipr_w$, into the coefficient group $F_k$,
that restricts to a valid proof that the given pair,
$(p,q)$ is in $E(p,q)$, and so that the homomorphism $h$, and all the homomorphisms
in its strictly solid family,  do not factor through a free product in  which
$<p>$ is contained in one factor, and $<q>$ is contained in a second factor. 

Based on the existence of such a homomorphism, we have associated (in section 3) with the equivalence class of $p_0$ at least
one uniformization limit group from the finite collection,
$Unif_1,\ldots,Unif_v$,  
so that the uniformization limit group 
satisfies the conclusions of theorem 3.11 and lemma 3.12 with respect to that equivalence class. 

From each uniformization limit group, $Unif_i$, the limit group $Ipr_{j(i)}$ that mapped into it  inherits a decomposition $\Lambda_i$.
$|Lambda_i$ gave us a separation of variables, i.e., $<p>$ is in one vertex group, and $<q>$ is in another vertex group, and the values that
the edge groups of $\Lambda_i$ obtain are a class function. i.e., for each equivalence class the edge groups obtain only boundedly many values up
to conjugation.

Note that the values of the edge groups in $\Lambda_i$ up to conjugation involves a simultaneous conjugation of finite tuples and not necessarily
conjugation of separate elements.
In a free group, a conjugacy class of a finite tuple of elements is determined by the conjugacy classes of finitely many elements. 
In a free group, if $u_1$ and $u_2$ don't
commute, and $v_1$ and $v_2$ don't commute, then if $u_1^{10}u_2^9=v_1^{10}v_1^9$ then $u_1=v_1$ and $u_2=v_2$. Hence, a non-commuting
 tuple, $u_1,u_2$, is conjugate
to a non-commuting tuple $v_1,v_2$, if and only if  $u_1^{10}u_2^9$ is conjugate to $v_1^{10}v_2^9$. If $u_1,u_2$ and $v_1,v_2$ are commuting pairs with no
non-trivial elements, then the two pairs are conjugate if and only if the individual elements are conjugate.

Similarly, two tuples, $u_1,u_2,u_3$ and $v_1,v_2,v_3$, 
with no pair of commuting
elements, are conjugate if and only if: $(u_1^{10}u_2^9)^{10}u_3^9$ is conjugate to 
$(v_1^{10}v_2^9)^{10}v_3^9$. Therefore, conjugation of two finite tuples can be definable reduced to separate 
conjugation of some finitely many words in the elements 
from the two tuples. 
 
Hence, the multi-function that associates the values of the edge groups in the  decompositions $\Lambda_i$ up to conjugacy is a
definable class multi-functions of the type that appears in the statement of theorem 4.2 (i.e., definable after adding new sorts for conjugacy classes).. 

With the  decomposition $\Lambda_i$ we associated finitely many multi-graded resolutions, and continued with a finite iterative procedure,
that finally constructed finitely many multi-graded resolutions. Each of these multi-graded resolutions contains a multi-graded
completion that contains the subgroup $<p>$, and a multi-graded completion that contains the subgroup $<q>$. The edge groups in that connect between
these two completions obtain boundedly many values for each equivalence class in $E(p,q)$. 

By the construction of the terminating iterative procedure, for each equivalence class in $E(p,q)$, that is not one of the finitely many
classes that were singled out in theorem 3.1, there must exist at least one multi-graded resolution that was constructed along the iterative procedure
that contains a multi-graded test sequence that is in the class.

With each multi-graded resolution that was constructed along the iterative procedure we have associated finitely many multi-graded closures
(according to definitions 1.25 to 1.30 in [Se5]) that collect the obstructions to the existence of test sequences that restrict to values $\{(p_n,q_n)\}$ 
that are in the relation $E(p,q)$. 

Every fiber of the multi-graded resolution (or its completion) that is not covered by these multi-graded
closures contain test sequences that are in $E(p,q)$, and the classes that contain test sequences in such a fiber are completely determined by the values
of the edge groups that connect between the completion that contains $<p>$ and the completion that contain $<q>$. These values are a class 
function,  for each class there are only boundedly many values of these edges up to conjugation, and these values are determined by a definable multi-function 
of the form that is indicated in theorem 4.2 (definable after adding new sorts for conjugacy classes). 

The collection of the definable multi-functions that are associated with all the multi-graded resolutions that were constructed along the
terminating iterative procedure that started with the graphs of groups $\Lambda_i$ (that was constructed from the uniformization limit groups)
give one definable multi-function in the form that is indicated in theorem 4.2 (after adding sorts for conjugacy classes), 
with values that are class functions, and values that separate between the classes of $E(p,q)$, and theorem 4.2 follows.

\line{\hss$\qed$}

Theorem 4.2 proves weak elimination of imaginaries in case the terminal limit groups of all the uniformization limit groups
that were constructed along the terminating iterative procedure  that we used to analyze definable equivalence relations are rigid, and
the uniformization limit groups do not contain abelian vertex groups in any of their levels. 

We already generalized the procedure for separation of variables to the case in which the terminal
limit groups of the uniformization and the collapse uniformization 
limit groups can be solid and the abelian decompositions that are associated with the various levels of the uniformization limit groups
may contain abelian vertex groups. In this general case we also obtain weak elimination of imaginaries but the new sorts that we need 
to add contain all the basic equivalence 
relations that are described in section 2.
 
\vglue 1pc
\proclaim{Theorem 4.4} Let $F_k$ be a (non-abelian) free group, and let $E(p,q)$ be a  definable equivalence relation over $F_k$.
If we add sorts for the imaginaries that are presented in section 2:
 conjugation, generalized double cosets (definition 2.7) and generalized conjugated double cosets (definition 2.9),  then $E(p,q)$ is weakly
eliminated. 

Suppose that $p$ and $q$ are $m$-tuples. There exist some integers $s$ and $t$ and a definable multi-function:
$$f: \, F^m  \ \to \ F^s \times R_1 \times \ldots \times R_t$$ 
where each of the $R_i$'s is a new sort for one of the 3 basic imaginaries 
(conjugation, generalized double cosets and generalized conjugated double cosets).
The image of an element is uniformly bounded, 
(and can be assumed to be of equal size), 
the multi-function
 is a class function, i.e., two elements in an equivalence class of $E(p,q)$ have the same image, and the multi-function
$f$ separates between classes, i.e., the images of elements from distinct equivalence classes is distinct.  
\endproclaim

\nfp Let $E(p,q)$ be a  definable equivalence relation. To prove that $E(p,q)$ can be weakly eliminated,
we need to construct a definable multi-function $f$ as described in the theorem. First, in theorem 3.1 we associated with
$E(p,q)$ finitely many rigid or solid limit groups, $Ipr_1,\ldots,Ipr_w$, so that for all but finitely many equivalence
classes, and for every pair $(p,q)$ in any of the remaining  equivalence classes,  
there is a rigid or a strictly solid homomorphism from one of these limit groups into the coefficient group $F_k$, that
restricts to a proof that the pair $(p,q)$  is in $E(p,q)$, and the homomorphism does not factor through a free product of
limit groups in which $<p>$ is contained in one factor, and $<q>$ is contained in a second factor. 

Proposition 4.3 shows that the finite collection of equivalence classes that are excluded
in theorem 3.1 is definable (the proof of proposition 4.3 is general and does not assume anything on the terminal limit groups of the
constructed completions).
Therefore, similarly to what we did in the proof of theorem 4.2, to prove theorem 4.4, we need to construct a (definable)
multi-function
with the properties that are listed in the statement of the theorem, that is defined
on the union of all the other
equivalence classes of $E(p,q)$.

We continue along the proof of theorem 4.2.
Let $p_0$ be a specialization of the  variables $p$, that does not
belong to one of the finitely many equivalence classes that are excluded in theorem 3.1.
Then for  each pair $(p,q) \in E(p,q)$, that are in the same equivalence class as $p_0$, 
there exists a 
rigid or a strictly solid homomorphism $h$ from one of the rigid or strictly solid 
limit group, $Ipr_1,\ldots,Ipr_w$, into the coefficient group $F_k$,
that restricts to a valid proof that the given pair,
$(p,q)$ is in $E(p,q)$, and so that the homomorphism $h$, and all the homomorphisms
in its strictly solid family,  do not factor through a free product in  which
$<p>$ is contained in one factor, and $<q>$ is contained in a second factor. 

Based on the existence of such a homomorphism, we have associated (in section 3) with the equivalence class of $p_0$ at least
one uniformization limit group from the finite collection,
$Unif_1,\ldots,Unif_v$,  
so that the uniformization limit group 
satisfies the conclusions of theorem 3.11 and lemma 3.13 with respect to that equivalence class. 

From each uniformization limit group, $Unif_i$, the limit group $Ipr_{j(i)}$ that mapped into it  inherits a decomposition $\Lambda_i$.
$\Lambda_i$ gave us a separation of variables, i.e., $<p>$ is in one vertex group, and $<q>$ is in another vertex group, and the values that
the edge groups of $\Lambda_i$ obtain are a class function. i.e., for each equivalence class the edge groups obtain only boundedly many values up
the basic equivalence relations that were defined in section 2:  
conjugation, generalized double cosets (definition 2.7) and generalized conjugated double cosets (definition 2.9).

As in the terminal rigid groups case, with the  decomposition $\Lambda_i$ we associated finitely many multi-graded resolutions, 
and continued with a finite iterative procedure,
that finally constructed finitely many multi-graded resolutions. Each of these multi-graded resolutions contains a multi-graded
completion that contains the subgroup $<p>$, and a multi-graded completion that contains the subgroup $<q>$. The edge groups in that connect between
these two completions obtain boundedly many values for each equivalence class in $E(p,q)$  up to the basic equivalence relations:
conjugation, generalized double cosets and generalized conjugated double cosets.

The rest of the proof is identical to the argument that was used to prove theorem 4.2 (in which we assumed that the terminal limit groups of the 
uniformization limit groups are rigid).

\line{\hss$\qed$}




By the results of [Se8], all the constructions that were associated with a (definable)  equivalence relation over a 
free group can be associated with a definable equivalence relation over a (non-abelian) torsion-free hyperbolic
group. Hence, torsion-free hyperbolic groups admit the same type of weak elimination of imaginaries as
a non-abelian free group.
 
\vglue 1pc
\proclaim{Theorem 4.5} Let $\Gamma$ be a non-elementary, torsion-free hyperbolic group, and let $E(p,q)$ be a 
definable equivalence relation over $\Gamma$. The conclusion of theorem 4.4 holds for $E(p,q)$.
If we add sorts for the imaginaries that are presented in section 2:
 conjugation, generalized double cosets (definition 2.7) and generalized conjugated double cosets (definition 2.9),  then $E(p,q)$ is weakly
eliminated  (weakly eliminated in the form that is stated in  theorem 4.4).
\endproclaim

\nfp By [Se8] the description of a definable set over a hyperbolic group is similar to the one over a free group.
In [Se8], the analysis of solutions to systems of equations, the construction of formal solutions, the analysis
of parametric equations, and in particular the uniform bounds on the number of rigid
and strictly solid families of solutions that are associated with a given value of the defining parameters, are
generalized to non-elementary, torsion-free hyperbolic groups. Furthermore, in [Se8] the sieve procedure
is generalized to torsion-free hyperbolic groups, as well as the
the analysis and the description of definable sets. Duo limit groups are defined and constructed over
torsion-free hyperbolic groups precisely as over free groups, and so are the Diophantine envelope and the duo envelope
(see section 1). 

Finally, theorem 1.5 in [Se3] that guarantees that given a f.g.\ group, and a sequence of homomorphisms from $G$
into a free group, there exists an integer $s$, and
a subsequence of the given homomorphisms that converges into a free action of some
limit quotient of $G$ on some $R^s$-tree, remains valid over torsion-free hyperbolic groups. This is the theorem that
is used to prove the termination of the iterative procedure for separation of variables, i.e., the iterative
procedure for the construction of uniformization limit groups.

Therefore, the procedure for separation of variables that was presented in the previous section
generalizes to torsion-free hyperbolic groups, and so is the modification of the sieve procedure that allows us
to run the sieve procedure while  preserving the separation of variables, that was described in this section. 
Hence, with a given
equivalence relation, $E(p,q)$, over a torsion-free hyperbolic group,
one can associate a finite collection of uniformization 
limit groups, precisely as over free groups. By the argument that was used to prove theorems 4.4 and 4.5
(that remains valid over torsion-free hyperbolic groups), this collection of  uniformization
limit groups, enables one to obtain weak elimination of imaginaries (over torsion-free hyperbolic groups),
when we add sorts for the  basic families of
imaginaries: 
 conjugation, generalized double cosets (definition 2.7) and generalized conjugated double cosets (definition 2.9), which implies theorem 4.5.

\line{\hss$\qed$}


\vskip 2em

\end